\documentclass[10pt]{article}
\usepackage{amsmath}
\usepackage{latexsym}
\usepackage{amssymb}
\usepackage{amsfonts}
\usepackage{amsthm}
\title{DEFORMATION  COHOMOLOGY  OF  \\  ALGEBRAIC  AND  GEOMETRIC  STRUCTURES}
\author{J.-F. Pommaret \\ CERMICS, Ecole des Ponts ParisTech,\\ 6/8 Av. Blaise Pascal, 77455 Marne-la-Vall\'ee Cedex 02, France \\
E-mail: jean-francois.pommaret@wanadoo.fr, pommaret@cermics.enpc.fr \\
URL: http://cermics.enpc.fr/~pommaret/home.html }
\date{  }
\begin{document}
\maketitle

\noindent
{\bf 1  INTRODUCTION}:\\

   In 1953 the physicists E. Inon\"{u} and E.P. Wigner (1963 Nobel prize) introduced the concept of {\it deformation of a Lie algebra} by considering the composition law $(u,v)\rightarrow (u+v)/(1+(uv/c^2))$ for speeds in special relativity (Poincar\' {e} group) when $c$ is the speed of light, claiming that the limit $c\rightarrow \infty$ or $1/c\rightarrow 0$ should produce the composition law $(u,v)\rightarrow u+v $ used in classical mechanics (Galil\' {e}e group) ([22]). However, this result is not correct indeed as $1/c\rightarrow 0$ has no meaning independently of the choice of length and time units. Hence one has to consider the dimensionless numbers $\bar{u}=u/c,\bar{v}=v/c$ in order to get $(\bar{u},\bar{v})\rightarrow (\bar{u}+\bar{v})/(1+\bar{u}\bar{v})$ with no longer any pertubation parameter involved ([32]). Nevertheless, this idea brought the birth of the theory of {\it deformation of algebraic structures} ([13],[19],[20],[21],[34],[35],[44],[45]), culminating in the use of the Chevalley-Eilenberg cohomology of Lie algebras ([9]) and one of the first applications of computer algebra in the seventies because a few counterexamples can only be found for Lie algebras of dimension $\geq 11$ ([2]). Finally, it must also be noticed that the main idea of general relativity is to deform the Minkowski metric $dx^2+dy^2+dz^2-c^2dt^2$ of space-time by means of the small dimensionless parameter $\phi /c^2$ where $\phi=GM/r$ is the gravitational potential at a distance $r$ of a central attractive mass $M$ with gravitational constant $G$ ([11],[33]).   \\
   
   Let $G$ be a Lie group with identity $e$ and Lie algebra ${\cal{G}}=T_e(G)$, the tangent space of $G$ at $e$. If $a=(a^{\tau})$ with $\tau =1,...,p$ are local coordinates on $G$, the {\it bracket} $[\cal{G},\cal{G}]\subset \cal{G}$ is defined by $p^2(p-1)/2$ {\it structure constants} $c=(c^{\tau}_{\rho\sigma}=-c^{\tau}_{\sigma\rho})$ satisfying the {\it Jacobi identities} $J(c)=0\Leftrightarrow c^{\lambda}_{\rho\sigma}c^{\mu}_{\lambda\tau}+c^{\lambda}_{\sigma\tau}c^{\mu}_{\alpha\rho}+c^{\lambda}_{\tau\rho}c^{\mu}_{\lambda\sigma}=0$. If now $c_t=c+tC+... $ is a deformation of $c$ satisfying $J(c_t)=0$, the idea is to study the vector space $\frac{\partial J(c)}{\partial c}C=0$ modulo the fact that $c$ behaves like a $3$-tensor under a change of basis of $\cal{G}$. The first condition implies $C\in Z_2(\cal{G})$ as a cocycle while the second implies $C\in B_2(\cal{G})$ as a coboundary, a result leading to introduce the Chevalley-Eilenberg cohomology groups $H_r({\cal{G}})=Z_r({\cal{G}})/B_r({\cal{G}})$ and to study $H_2(\cal{G})$ in particular.\\

  A few years later, a deformation theory was introduced for {\it structures on manifolds}, generally represented by {\it fields of geometric objects} like tensors and we may quote riemannian, symplectic, complex analytic or contact structures with works by H. Goldschmidt ([12],[14],[15]), R.E. Greene ([15]), V. Guillemin ([16],[17],[18]), K. Kodaira ([24],[25]), M. Kuranishi ([27]), L. Nirenberg ([24]), D.C. Spencer ([15],[25],[47],[48]) or S. Sternberg ([16],[17],[18]). The idea is to make the underlying geometric objects depending on a parameter while satisfying prescribed integrability conditions. The link betwen the two approaches, though often conjectured, has never been exhibited by the above authors and our aim is to sketch the solution of this problem that we have presented in many books, in particular to study the possibility to use computer algebra for such a purpose. It must be noticed that the general concept of natural bundle and geometric object is absent from the work of Spencer and coworkers ([27],[47],[48]]) though it has been discovered by E. Vessiot as early as in 1903 ([51]). It must also be noticed that the introduction of the book "Lie equations" ([26]) dealing with the "classical" structures on manifolds has nothing to do with the remaining of the book dealing with the nonlinear Spencer sequence and where all the vector bundles and nonlinear operators are different (See ([43]) for more details). This work is a natural continuation of symbolic computations done at RWTH-Aachen university by M. Barakat and A. Lorenz in 2008 ([3],[30]). \\
  
    The starting point is to refer to the three fundamental theorems of Lie, in particular the third which provides a way to realize an analytic Lie group $G$ from the knowledge of its Lie algebra $\cal{G}$ and the construction of the left or right invariant analytic $1$-forms on $G$ called {\it Maurer-Cartan forms} ${\omega}^{\tau}=({\omega}^{\tau}_{\sigma}(a)da^{\sigma})$ satisfying the well known {\it Maurer-Cartan equations} $d{\omega}^{\tau}+c^{\tau}_{\rho \sigma}{\omega}^{\rho}\wedge {\omega}^{\sigma}=0$ where $d$ is the exterior derivative, provided that $J(c)=0$.\\
    
    The general solution involves tools from differential geometry (Spencer operator, $\delta$-cohomology), differential algebra, algebraic geometry, algebraic analysis and homological algebra. However, the key argument is to acknowledge the fact that the {\it Vessiot structure equations} ([51],1903, still unknown today !) must be used in place of the {\it Cartan structure equations} ([4],[5]), along the following diagram describing the reference ([43]) which is the latest one existing on this subject while providing applications to engineering (continuum mechanics, electromagnetism) and mathematical (gauge theory, general relativity) physics.\\

\[   \begin{array}{rcccc}
   &   &  CARTAN   &   \longrightarrow   &  SPENCER  \\
      &   \nearrow  &   &   &   \\
 LIE   &   &   \updownarrow   &    &   \updownarrow   \\
      &   \searrow  &  &  &   \\
   &   &   VESSIOT   &   \longrightarrow   &  JANET   
   \end{array}     \]
\vspace*{1mm}  \\

  In a word, one has to replace Lie algebras by Lie algebroids with a bracket now depending on the Spencer operator and use the corresponding canonical {\it linear Janet sequence} in order to induce a new sequence locally described by finite dimensional vector spaces and linear maps, {\it even for structures providing infinite dimensional Lie algebroids} (contact and unimodular contact structures are good examples as we shall see). The cohomology of this sequence only depends on the "{\it structure constants} " appearing in the Vessiot structure equations (constant riemannian curvature is an example with only one constant having of course nothing to do with any Lie algebra !). Finally, the simplest case of a principal homogeneous space for $G$ (for example $G$ itself as before) gives back the Chevalley-Eilenberg cohomology.\\
  
  In order to motivate the reader and convince him about the novelty of the underlying concepts, we provide in a rather self-contained way and parallel manners the following  four striking examples which are the best nontrivial ones we know showing how and why computer algebra could be used and that will be revisited later on. We invite the reader at this stage to try to imagine any possible link that could exist between these examples.   \\

\noindent
{\bf EXAMPLE  1.1}: ({\it Principal homogeneous structure}) When $\Gamma$ is the Lie group of transformations made by the constant translations  $y^i=x^i+a^i$ for $i=1,...,n$ of a manifold $X$ with $dim(X)=n$, the characteristic object invariant by $\Gamma$ is a family $\omega=({\omega}^{\tau}={\omega}^{\tau}_idx^i)$ of $n$ $1$-forms with $det(\omega)\neq 0$ in such a way that $\Gamma =\{f\in aut(X){\mid}j_1(f)^{-1}(\omega)=\omega\}$ when $aut(X)$ denotes the pseudogroup of local diffeomorphisms of $X$, $j_q(f)$ denotes the derivatives of $f$ up to order $q$ and $j_1(f)$ acts in the usual way on covariant tensors. For any vector field $\xi\in T=T(X)$ the tangent bundle to $X$, introducing the standard Lie derivative ${\cal{L}}(\xi)$ of forms with respect to $\xi$, we may consider the $n^2$ 
{\it  first order Medolaghi equations}:  \\
\[   {\Omega}^{\tau}_i\equiv ({\cal{L}}(\xi)\omega)^{\tau}_i\equiv {\omega}^{\tau}_r(x){\partial}_i{\xi}^r+{\xi}^r{\partial}_r{\omega}^{\tau}_i(x)=0  \]
The particular situation is found with the special choice $\omega=(dx^i)$ that leads to the involutive system ${\partial}_i{\xi}^k=0$. Introducing the inverse matrix $\alpha=({\alpha}^i_{\tau})={\omega}^{-1}$, the above equations amount to the bracket relations $[\xi,{\alpha}_{\tau}]=0$ and, using crossed derivatives on the {\it solved form} ${\partial}_i{\xi}^k+{\xi}^r{\alpha}^k_{\tau}(x){\partial}_r{\omega}^{\tau}_i(x)=0$, we obtain the $n^2(n-1)/2$ zero order equations:  \\
\[   {\xi}^r{\partial}_r({\alpha}^i_{\rho}(x){\alpha}^j_{\sigma}(x)({\partial}_i{\omega}^{\tau}_j(x)-{\partial}_j{\omega}^{\tau}_i(x)))=0   \]
The {\it integrability conditions} (IC) are thus the $n^2(n-1)/2$ {\it Vessiot structure equations}:  \\
\[    {\partial}_i{\omega}^{\tau}_j(x)-{\partial}_j{\omega}^{\tau}_i(x)=c^{\tau}_{\rho\sigma}{\omega}^{\rho}_i(x){\omega}^{\sigma}_j(x)   \]
with $n^2(n-1)/2$ {\it structure constants} $c=(c^{\tau}_{\rho\sigma}=-c^{\tau}_{\sigma\rho})$. When $X=G$, these equations can be identified with the {\it Maurer-Cartan equations} (MC) existing in the theory of Lie groups, on the condition to change the sign of the structure constants involved because we have $[{\alpha}_{\rho},{\alpha}_{\sigma}]= - c^{\tau}_{\rho\sigma}{\alpha}_{\tau}$. Writing these equations in the form of the exterior system $d{\omega}^{\tau}=c^{\tau}_{\rho\sigma}{\omega}^{\rho}\wedge{\omega}^{\sigma}$ and closing this system by applying once more the exterior derivative $d$, we obtain the quadratic IC:   \\
\[   c^{\lambda}_{\mu\rho}c^{\mu}_{\sigma\tau}+c^{\lambda}_{\mu\sigma}c^{\mu}_{\tau\rho}+c^{\lambda}_{\mu\tau}c^{\mu}_{\rho\sigma}=0  \]
also called {\it Jacobi relations} $J(c)=0$. Finally, if another family $\bar{\omega}$ of forms is given, the {\it equivalence problem} 
$j_1(f)^{-1}(\omega)=\bar{\omega}$ cannot be solved even locally if $\bar{c}\neq c$.\\

\noindent
{\bf EXAMPLE  1.2}: ({\it Riemann structure}) If $\omega=({\omega}_{ij}={\omega}_{ji})$ is a metric on a manifold $X$ with $dim(X)=n$ such that $det(\omega)\neq 0$, the Lie pseudogroup of transformations preserving $\omega$ is $\Gamma=\{f\in aut(X){\mid}j_1(f)^{-1}(\omega)=\omega \}$ and is a Lie group with a maximum number of $n(n+1)/2$ parameters. A special metric could be the Euclidean metric when $n=1,2,3$ as in elasticity theory or the Minkowski metric when $n=4$ as in special relativity. The 
{\it  first order Medolaghi equations}:\\
\[ {\Omega}_{ij}\equiv ({\cal{L}}(\xi)\omega)_{ij}\equiv {\omega}_{rj}(x){\partial}_i{\xi}^r+{\omega}_{ir}(x){\partial}_j{\xi}^r+{\xi}^r{\partial}_r{\omega}_{ij}(x)=0 \]
are also called {\it Killing equations} for historical reasons. The main problem is that {\it this system is not involutive} unless we prolong it to order two by differentiating once the equations. For such a purpose, introducing ${\omega}^{-1}=({\omega}^{ij})$ as usual, we may define the {\it Christoffel symbols}:\\
\[ {\gamma}^k_{ij}(x)=\frac{1}{2}{\omega}^{kr}(x)({\partial}_i{\omega}_{rj}(x) +{\partial}_j  {\omega}_{ri}(x) -{\partial}_r{\omega}_{ij}(x))=
{\gamma}^k_{ji}(x) \]
This is a new geometric object of order $2$ providing an isomorphism $j_1(\omega)\simeq (\omega,\gamma)$ and allowing to obtain the {\it second order Medolaghi equations}:\\
\[ {\Gamma}^k_{ij}\equiv ({\cal{L}}(\xi)\gamma)^k_{ij}\equiv {\partial}_{ij}{\xi}^k+{\gamma}^k_{rj}(x){\partial}_i{\xi}^r+{\gamma}^k_{ir}(x){\partial}_j{\xi}^r-{\gamma}^r_{ij}(x){\partial}_r{\xi}^k+{\xi}^r{\partial}_r{\gamma}^k_{ij}(x)=0   \]
Surprisingly, the following expression:\\
\[ {\rho}^k_{lij}(x)\equiv {\partial}_i{\gamma}^k_{lj}(x)-{\partial}_j{\gamma}^k_{li}(x)+{\gamma}^r_{lj}(x){\gamma}^k_{ri}(x)-{\gamma}^r_{li}(x){\gamma}^k_{rj}(x)  \]
is still a first order geometric object and even a $4$-tensor with $n^2(n^2-1)/12$ independent components satisfying the purely algebraic relations :\\
 \[   {\rho}^k_{lij}+{\rho}^k_{ijl}+{\rho}^k_{jli}=0, \hspace{5mm} {\omega}_{rl}{\rho}^l_{kij}+{\omega}_{kr}{\rho}^r_{lij}=0  \]
Accordingly, the IC must express that the new first order equations $R^k_{lij}\equiv ({\cal{L}}(\xi)\rho)^k_{lij}=0$ are only linear combinations of the previous ones and we get the {\it Vessiot structure equations} with the only {\it structure constant} $c$:\\
    \[  {\rho}^k_{lij}(x)=c({\delta}^k_i{\omega}_{lj}(x)-{\delta}^k_j{\omega}_{li}(x))   \]
 describing the constant Riemannian curvature condition of Eisenhart ([10]). If $\bar{\omega}$ is another nondegenerate metric with structure constant 
 $\bar{c}$, the {\it equivalence problem} $j_1(f)^{-1}(\omega)=\bar{\omega}$ cannot be solved even locally if $\bar{c}\neq c$.\\ 

\noindent
{\bf EXAMPLE  1.3}: ({\it Contact structure}) We only treat the case $dim(X)=3$ as the case $dim(X)={2p+1}$ needs much more work ([37]). Let us consider the so-called {\it contact} $1$-form $\alpha=dx^1-x^3dx^2$ and consider the Lie pseudogroup $\Gamma\subset aut(X)$ of (local) transformations preserving $\alpha$ up to a function factor, that is $\Gamma=\{f\in aut(X){\mid}j_1(f)^{-1}(\alpha)= \rho \alpha\}$ where again $j_q(f)$ is a symbolic way for writing out the derivatives of $f$ up to order $q$ and $\alpha$ transforms like a $1$-covariant tensor. It may be tempting to look for a kind of "{\it object} " the invariance of which should characterize $\Gamma$. Introducing the exterior derivative $d\alpha=dx^2\wedge dx^3$ as a $2$-form, we obtain the volume $3$-form $\alpha\wedge d\alpha=dx^1\wedge dx^2\wedge dx^3$. As it is well known that the exterior derivative commutes with any diffeomorphism, we obtain sucessively:\\
\[      j_1(f)^{-1}(d\alpha)=d(j_1(f)^{-1}(\alpha))=d(\rho \alpha)=\rho d\alpha +d\rho \wedge \alpha  \Rightarrow j_1(f)^{-1}(\alpha \wedge d\alpha)={\rho}^2(\alpha\wedge d\alpha)   \]
As the volume $3$-form $\alpha\wedge d\alpha$ transforms through a division by the Jacobian determinant $\Delta=\partial (f^1,f^2,f^3)/\partial (x^1,x^2,x^3)\neq 0$ of the transformation $y=f(x)$ with inverse $x=f^{-1}(y)=g(y)$, {\it the desired object is thus no longer a} $1$-{\it form but a} $1$-{\it form density} $\omega=({\omega}_1,{\omega}_2,{\omega}_3)$ transforming like a $1$-form but up to a division by the square root of the Jacobian determinant. It follows that the infinitesimal contact transformations are vector fields $\xi\in T=T(X)$ the tangent bundle of $X$, satisfying the $3$ so-called 
{\it  first order Medolaghi equations}:ÊÊÊ\\
\[  {\Omega}_i\equiv ({\cal{L}}(\xi)\omega)_i\equiv {\omega}_r(x){\partial}_i{\xi}^r-(1/2){\omega}_i(x){\partial}_r{\xi}^r+{\xi}^r{\partial}_r{\omega}_i(x)=0   \]
When $\omega=(1,-x^3,0)$, we obtain the {\it special} involutive system:ÊÊ\\
\[  {\partial}_3{\xi}^3+{\partial}_2{\xi}^2+2x^3{\partial}_1{\xi}^2-{\partial}_1{\xi}^1=0, {\partial}_3{\xi}^1-x^3{\partial}_3{\xi}^2=0, {\partial}_2{\xi}^1-x^3{\partial}_2{\xi}^2+x^3{\partial}_1{\xi}^1-(x^3)^2{\partial}_1{\xi}^2-{\xi}^3=0   \]
with  $2$ equations of class $3$ and $1$ equation of class 2 obtained by exchanging $x^1$ and $x^3$ (see section 4 for the definitions) and thus only $1$ {\it compatibility conditions} (CC) for the second members.\\
For an arbitrary $\omega$, we may ask about the differential conditions on $\omega$ such that all the equations of order $r+1$ are only obtained by differentiating $r$ times the first order equations, exactly like in the special situation just considered where the system is involutive. We notice that, in a symbolic way, $\omega \wedge d\omega$ is now a scalar $c(x)$ providing the zero order equation ${\xi}^r{\partial}_rc(x)=0$ and the condition is $c(x)=c=cst$. The {\it integrability condition} (IC) is the {\it Vessiot structure equation}:   \\
\[  {\omega}_1({\partial}_2{\omega}_3-{\partial}_3{\omega}_2)+{\omega}_2({\partial}_3{\omega}_1-{\partial}_1{\omega}_3)+
{\omega}_3({\partial}_1{\omega}_2-{\partial}_2{\omega}_1)=c    \]
involving the only {\it structure constant} $c$.\\
For $\omega=(1,-x^3,0)$, we get $c=1$. If we choose $\bar{\omega}=(1,0,0)$ leading to $\bar{c}=0$, we may define $\bar{\Gamma}=\{f\in aut(X){\mid} j_1(f)^{-1}(\bar{\omega})=\bar{\omega}\}$ with infinitesimal transformations satisfying the involutive system:\\
\[  {\partial}_3{\xi}^3+{\partial}_2{\xi}^2-{\partial}_1{\xi}^1=0, {\partial}_3{\xi}^1=0, {\partial}_2{\xi}^1=0  \]
with again $2$ equations of class $3$ and $1$ equation of class $2$. The {\it equivalence problem} $j_1(f)^{-1}(\omega)=\bar{\omega}$ cannot be solved even locally because this system cannot have any invertible solution. Indeed, studying the system $j_1(g)^{-1}(\bar{\omega})=\omega$, we have to solve:  \\
\[       \frac{\partial g^1}{\partial y^2}+y^3\frac{\partial g^1}{\partial y^1}=0, \frac{\partial g^1}{\partial y^3}=0 \Rightarrow \frac{\partial g^1}{\partial y^1}=0, \frac{\partial g^1}{\partial y^2}=0, \frac{\partial g^1}{\partial y^3}=0  \]
by using crossed derivatives.\\

\noindent
{\bf EXAMPLE  1.4}: ({\it Unimodular contact structure})  With similar notations, let us again set $\alpha=dx^1-x^3dx^2\Rightarrow d\alpha=dx^2\wedge dx^3$  but let us now consider the new Lie pseudogroup  of transformations preserving $\alpha$ and thus $d\alpha$ too, that is preserving the {\it mixed} object $\omega=(\alpha,\beta)$ made up by a $1$-form $\alpha$ and a $2$-form $\beta$ with 
$\gamma = \alpha\wedge \beta\neq 0$ and $d\alpha=\beta\Rightarrow d\beta=0$. Then $\Gamma$ is a Lie subpseudogroup of the one just considered in the previous example and the corresponding infinitesimal transformations now satisfy the involutive system:\\
\[   {\partial}_1{\xi}^1=0, {\partial}_1{\xi}^2=0, {\partial}_1{\xi}^3=0, {\partial}_2{\xi}^1+x^3{\partial}_3{\xi}^3-{\xi}^3=0, {\partial}_2{\xi}^2+{\partial}_3{\xi}^3=0, {\partial}_3{\xi}^1-x^3{\partial}_3{\xi}^2=0  \]
with $3$ equations of class $3$, $2$ equations of class $2$ and $1$ equation of class $1$ if we exchange $x^1$ with $x^3$, a result leading now to $4$ CC. \\
More generally, when $\omega=(\alpha,\beta)$ where $\alpha$ is a $1$-form and $\beta$ is a $2$-form satifying $\alpha\wedge \beta\neq 0$, {\it we may study the same problem as before} for the {\it general} system ${\cal{L}}(\xi)\alpha=0, {\cal{L}}(\xi)\beta=0$ where ${\cal{L}}(\xi)$ is the standard Lie derivative of forms with respect to the vector field $\xi$, that is ${\cal{L}}(\xi)=i(\xi)d+di(\xi)$ if $i(\xi)$ is the interior multiplication of a form by the vector field $\xi$. We now provide details on the tedious computation involved as it is at this point that computer algebra may be used. With $\alpha=({\alpha}_i),\beta=({\beta}_{ij}=-{\beta}_{ji})$ we may suppose, with no loss of generality, that ${\alpha}_1\neq 0, {\beta}_{23}\neq 0$ in such a way that $\alpha\wedge \beta\neq 0 \Leftrightarrow {\alpha}_1{\beta}_{23}+{\alpha}_2{\beta}_{31}+{\alpha}_3{\beta}_{12}\neq 0$. We can then solve the general equations as in the special situation already considered with respect to the $6$ leading principal derivatives $pri=\{{\partial}_1{\xi}^1,{\partial}_1{\xi}^2,{\partial}_1{\xi}^3, {\partial}_2{\xi}^1, {\partial}_2{\xi}^2, {\partial}_3{\xi}^1\}$ as involution is an intrinsic local property and obtain:  \\
\[  {\partial}_1{\xi}^1 +\frac{{\alpha}_2}{{\alpha}_1}{\partial}_1{\xi}^2+\frac{{\alpha}_3}{{\alpha}_1}{\partial}_1{\xi}^3+ ... =0 \]
\[  {\partial}_1{\xi}^2 +\frac{{\beta}_{31}}{{\beta}_{32}}({\partial}_1{\xi}^1+{\partial}_3{\xi}^3)+\frac{{\beta}_{21}}{{\beta}_{32}}{\partial}_3{\xi}^2+ ...=0 \]
\[  {\partial}_1{\xi}^3 +\frac{{\beta}_{13}}{{\beta}_{32}}{\partial}_2{\xi}^3+\frac{{\beta}_{12}}{{\beta}_{32}}({\partial}_1{\xi}^1+{\partial}_2{\xi}^2)+ ... =0 \]
\[  {\partial}_2{\xi}^1 +\frac{{\alpha}_2}{{\alpha}_1}{\partial}_2{\xi}^2+\frac{{\alpha}_3}{{\alpha}_1}{\partial}_2{\xi}^3+ ... =0  \]
\[  {\partial}_2{\xi}^2 +{\partial}_3{\xi}^3+\frac{{\beta}_{13}}{{\beta}_{23}}{\partial}_2{\xi}^1+\frac{{\beta}_{21}}{{\beta}_{23}}{\partial}_3{\xi}^1+ ... =0 \]
\[  {\partial}_3{\xi}^1+\frac{{\alpha}_2}{{\alpha}_1}{\partial}_3{\xi}^2+\frac{{\alpha}_3}{{\alpha}_1}{\partial}_3{\xi}^3+ ... =0 \]
Solving wth respect to the 6 principal derivatives in a triangular way while introducing the $3$ parametric derivatives $par=\{{\partial}_2{\xi}^3, {\partial}_3{\xi}^2, {\partial}_3{\xi}^3\}$, we obtain for example:\\
\[  (1-\frac{{\beta}_{13}}{{\beta}_{23}}\frac{{\alpha}_2}{{\alpha}_1}){\partial}_2{\xi}^2-\frac{{\beta}_{21}}{{\beta}_{23}}\frac{{\alpha}_2}{{\alpha}_1}{\partial}_3{\xi}^2+(1-\frac{{\beta}_{21}}{{\beta}_{23}}\frac{{\alpha}_3}{{\alpha}_1}){\partial}_3{\xi}^3-\frac{{\beta}_{13}}{{\beta}_{23}}\frac{{\alpha}_3}{{\alpha}_1}{\partial}_2{\xi}^3 + ... =0                   \]
Setting now $\bar{\beta}=d\alpha$ and identifying the corresponding coefficients, then ${\cal{L}}(\xi)\bar{\beta}$ is a linear combination of ${\cal{L}}(\xi)\alpha$ and ${\cal{L}}(\xi)\beta$ if and only if we have for example: \\
\[  \frac{\frac{{\beta}_{13}}{{\beta}_{23}}\frac{{\alpha}_3}{{\alpha}_1}}{1-\frac{{\beta}_{13}}{{\beta}_{23}}\frac{{\alpha}_2}{{\alpha}_1}}= \frac{\frac{{\bar{\beta}}_{13}}{{\bar{\beta}}_{23}}\frac{{\alpha}_3}{{\alpha}_1}}{1-\frac{{\bar{\beta}}_{13}}{{\bar{\beta}}_{23}}\frac{{\alpha}_2}{{\alpha}_1}}\]
and thus:  \\
\[\frac{{\beta}_{13}}{{\beta}_{23}}=\frac{{\bar{\beta}}_{13}}{{\bar{\beta}}_{23}}\Rightarrow \frac{{\bar{\beta}}_{23}}{{\beta}_{23}}=\frac{{\bar{\beta}}_{13}}{{\beta}_{23}}=\frac{{\bar{\beta}}_{12}}{{\beta}_{12}}    \]
a result, not evident at first sight, showing that the $2$-form $d\alpha$ {\it must} be proportional to the $2$-form $\beta$, that is $d\alpha=c'(x)\beta$ and thus $\alpha\wedge d\alpha=c'(x)\alpha\wedge\beta$. As $\alpha\wedge\beta\neq 0$, we {\it must} have $c'(x)=c'=cst$ and thus $d\alpha=c'\beta$. Similarly, we get $d\beta=c''\alpha\wedge\beta$ and obtain finally the $4$ {\it Vessiot structure equations} $d\alpha=c'\beta,d\beta=c''\alpha\wedge\beta$ involving $2$ {\it structure constants} $c=(c',c'')$. Contrary to the previous situation (but like in the Riemann case !) we notice that we have now $2$ structure equations not containing any constant (called {\it first kind} by Vessiot) and $2$ structure equations with the same number of different constants (called {\it second kind} by Vessiot), namely $\alpha\wedge d\alpha=c'\alpha\wedge\beta, d\beta=c''\alpha\wedge\beta$.\\
Finally, closing this system by taking once more the exterior derivative, we get $0=d^2\alpha=c'd\beta=c'c''\alpha\wedge\beta$ and thus the unexpected purely algebraic {\it Jacobi condition} $c'c''=0$. For the special choice $\omega=(dx^1-x^3dx^2, dx^2\wedge dx^3)$ we get $c=(1,0)$, for the second special choice $\bar{\omega}=(dx^1,dx^2\wedge dx^3)$ we get $\bar{c}=(0,0)$ and for the third special choice $\bar{\bar{\omega}}=((1/x^1)dx^1, x^1dx^2\wedge dx^3)$ we get $\bar{\bar{c}}=(0,1)$ with similar comments as before for the possibility to solve the corresponding equivalence problems.\\

\noindent
{\bf REMARK  1.5}: Comparing the various Vessiot structure equations containing structure constants that we have just presented and that we recall below in a symbolic way, we notice that {\it these structure constants are absolutely on equal footing} though they have nothing to do with any Lie algebra. 
\[  \begin{array}{rcl}
\partial \omega - \partial \omega & = & c \hspace{2mm} \omega \hspace{1mm} \omega   \\
\partial \gamma -\partial \gamma + \gamma \gamma -\gamma \gamma & = & c \hspace{2mm}(\delta \omega-\delta \omega)  \\
\omega \wedge (\partial \omega - \partial \omega) & = & c 
\end{array}  \]
\[ \hspace{9mm} \left\{
\begin{array}{rcl}
\partial \alpha - \partial \alpha & = & c'\hspace{2mm} \beta  \\
\partial  \beta - \partial \beta & = & c'' \hspace{2mm}\alpha \wedge \beta
\end{array}  
\right.    \]
Accordingly, the fact that the ones appearing in the MC equations are related to a Lie algebra is a pure coincidence and we may even say that {\it the Cartan structure equations have nothing to do with the Vessiot structure equations}. Also, as their factors are either constant, linear or quadratic, {\it any identification of the quadratic terms appearing in the Riemann tensor with the quadratic terms appearing in the MC equations is definitively not correct} even though most of mathematical physics today is based on such a confusion ([43]). Meanwhile, we understand why the torsion is {\it automatically combined} with curvature in the Cartan structure equations but {\it totally absent} from the Vessiot structure equations, even though the underlying group (translations + rotations) is the same. In addition, despite the prophetic comments of the italian mathematician Ugo Amaldi in 1909 ([1]), we do believe that it has been a pity that Cartan deliberately ignored the work of Vessiot at the beginning of the last century and that the things did not improve afterwards in the eighties with Spencer and coworkers (Compare MR 720863 (85m:12004) and MR 954613 (90e:58166)).\\

In the second section of this paper, which is an extended version of a lecture given at the international conference SCA (Symbolic Computation and its Applications) 2012 held in Aachen (Aix-la-Chapelle), Germany, May 17-20, 2012, we shall recall the definition of the Chevalley-Eilenberg cohomology existing for Lie algebras and describe its use in the study of the {\it deformation theory of Lie algebras} ([28],[29],[45],[49],[50]). We insist on the fact that the challenge solved in the remaining of the paper is not to generalize this result to arbitrary Lie equations but rather to work out a general framework that will provide {\it exactly} the deformation cohomology of Lie algebras in the particular case of Example  1.\\

In the third section we study the "{\it Vessiot structure equations}" in the nonlinear framework and explain why they must contain "{\it structure constants}" satifying algebraic Jacobi-like conditions.\\

In the fourth section, we present the minimum amount of differential geometry (jet theory, Spencer operator, $\delta$-cohomology, differential sequences) needed in order to achieve the formal constructions done in the next sections.  \\

In the fifth and longest section we present for the first time a rather self-contained but complete general procedure for exhibiting a deformation cohomology for any system of transitive Lie equations, even of infinite type as in the case of the (unimodular) contact structure. A link with the difficult concept of a {\it truncated Lie algebra} will also be provided (Compare to [14], II).  \\

Finally, we conclude the paper with a few explicit computations based on the previous examples while showing out the possibility to use computer algebra techniques.\\
\vspace*{5mm}  \\
\noindent
{\bf 2  DEFORMATION THEORY OF LIE ALGEBRAS}: \\

   Let $V$ be finite dimensional vector space over a field $k$ containing the field $\mathbb{Q}$ as a subfield and set $V^*=hom_k(V,k)$. We shall denote the elements of $V$ by $X, Y, Z,...$ with components $X=(X^{\rho})$ for $\rho=1,...,dim(V)$.\\
   
\noindent
{\bf DEFINITION 2.1}: A {\it Lie algebra} ${\cal{G}}=(V,c)$ is an algebraic structure on $V$ defined over $k$ by a bilinear map $[ \hspace{2mm} ]:{\wedge}^2V\rightarrow V$ called {\it bracket} through the formula $([X,Y])^{\tau}=c^{\tau}_{\rho\sigma}X^{\rho}Y^{\sigma}$ where the {\it structure constants} 
$c\in {\wedge}^2V^*\otimes V$ are satisfying the {\it Jacobi relations}:\\
\[   J(c)=0 \hspace{2cm}  c^{\lambda}_{\rho\sigma}c^{\mu}_{\lambda\tau}+c^{\lambda}_{\sigma\tau}c^{\mu}_{\alpha\rho}+c^{\lambda}_{\tau\rho}c^{\mu}_{\lambda\sigma}=0  \]
a result leading to the {\it Jacobi identity} for the bracket, namely:   \\
\[       [X,[Y,Z]] + [Y,[Z,X]] + [Z,[X,Y]]  \equiv 0 \hspace{1cm}, \forall X,Y,Z\in V   \]
Considering a Lie algebra as a point $c$ of the algebraic set defined by the quadratic equations $J(c)=0$, we may state:   \\

\noindent
{\bf DEFINITION 2.2}: A {\it deformation} $c_t$ of $c$ is a curve passing through $c$ in this algebraic set, that is to say a set of points $c_t$ indexed by a parameter $t$ and such that $J(c_t)=0$ with $c_t=c+tC+...$. As a byproduct, the Lie algebra ${\cal{G}}_t=(V,c_t)$ is called a deformation of ${\cal{G}}=(V,c)$.  \\

Of course, a central problem in the theory of Lie algebras is to exhibit properties of $\cal{G}$ that do not depend on the basis chosen for $V$. In particular, if $a\in aut(V) \subset V^*\otimes V$, we may define an equivalence relation among the structure constants as follows:\\

\noindent
{\bf DEFINITION 2.3}: $\bar{c}\sim c \Leftrightarrow a\bar{c}(X,Y)=c(aX,aY), \forall X,Y\in V$ or equivalently $\bar{c}=aaa^{-1}c$ in a symbolic tensorial notation showing that $c$ transforms like a ($2$-covariant, $1$-contravariant)-tensor under the action of $a$. In particular, we get at once $\bar{c}\sim c \Leftrightarrow J(\bar{c})=0$ whenever $J(c)=0$.   \\

When $a_t\in aut(V)$ is such that $a_0=id_V$, then the set of points $c_t\sim c$ is a deformation of $c$ called {\it trivial} and we may state:\\

\noindent
{\bf DEFINITION 2.4}: $a\in aut(V)$ is called an {\it automorphism} of $\cal{G}$ if $\bar{c}=c$ and $A\in end(V)=V^*\otimes V$ is called a {\it derivation} of $\cal{G}$ if $A[X,Y]=[AX,Y]+[X,AY], \forall X,Y\in V$. In particular, if $a_t=a_0+tA+ ...\in aut({\cal{G}})$ is such that $a_0=id_V$, then $A={\frac{da_t}{dt}}{\mid}_{t=0}$ is a derivation of $\cal{G}$. Moreover, if we define the {\it adjoint action} of $V$ on $V$ by $ad(X)Y=[X,Y]$, then it follows from the Jacobi identity that $ad(X)$ is a derivation of $\cal{G}$ called 
{\it inner derivation}.\\

\noindent
{\bf DEFINITION 2.5}: A Lie algebra $\cal{G}$ is said to be {\it rigid} if it cannot admit a deformation which is not trivial.  \\

As the algebraic set of structure constants may have very bad local properties, it may be interesting to study infinitesimal deformations when $t$ is a small parameter, that is $t\ll 1$. In particular, we should have $\frac{\partial J(c)}{\partial c}C=0$ for any deformation while a trivial deformation should lead to:  \\
\[  ({\bar{C}})^{\tau}_{\rho\sigma}=(A)^{\mu}_{\rho}(C)^{\tau}_{\mu\sigma}+(A)^{\mu}_{\sigma}(C)^{\tau}_{\rho\mu}-(A)^{\tau}_{\mu}(C)^{\mu}_{\rho\sigma}   \]
with no way to unify all these technical formulas.\\

   By chance, such an infinitesimal study can be made easy by means of the {\it Chevalley-Eilenberg cohomology} that we now describe.. For this, by analogy with the exterior derivative, let us define an application $d:{\wedge}^rV^*\otimes V \rightarrow {\wedge}^{r+1}V^*\otimes V$ depending on 
   $c$ by the formula:ÊÊ\\
\[  \begin{array}{rl}
df(X_1, ... , X_{r+1})=  &\sum_{i<j} (-1)^{i+j}f([X_i,X_j],X_1, ... ,{\hat{X}}_i, ... ,{\hat{X}}_j, ... ,X_{r+1}) \\
    &+ \sum_i (-1)^{i+1}[X_i, f(X_1, ... ,{\hat{X}}_i, ... ,X_{r+1})]    
\end{array}    \]
where a "hat" indicates an omission. Using the Jacobi identity for the bracket and a straightforward but tedious computation left to the reader as an exercise, one can prove that $d\circ d=0$. Hence we may define in the usual way {\it coboundaries} $B_r({\cal{G}})$ as images of $d$, {\it cocycles} $Z_r({\cal{G}})$ as kernels of $d$ and {\it cohomology groups} $H_r({\cal{G}})=Z_r({\cal{G}})/B_r({\cal{G}})$ in such a way that $B_r({\cal{G}})\subseteq Z_r({\cal{G}})\subseteq {\wedge}^rV^*\otimes V$.   \\

\noindent
{\bf LEMMA 2.6}: $C=\frac{dc_t}{dt}{\mid}_{t=0}\in Z_2(\cal{G})$ and $c_t\sim c \Leftrightarrow C=\frac{dc_t}{dt}{\mid}_{t=0}\in B_2(\cal{G}) $ . Accordingly, a sufficient condition of rigidity is $H_2({\cal{G}})=0 $.\\

\noindent
{\bf LEMMA 2.7}: $A$ is a derivation $\Rightarrow A\in Z_1({\cal{G}}) $ and $A$ is an inner derivation $\Rightarrow A\in B_1({\cal{G}}) $. Accordingly, the vector space of derivations of ${\cal{G}}$ modulo the vector space of inner derivations is $H_1({\cal{G}})$.  \\

Up to the moment we have only been looking at infinitesimal deformations. However, a Lie algebra ${\cal{G}}$ may be rigid on the finite level even if $H_2({\cal{G}})\neq 0$, that is even if it can be deformed on the infinitesimal level ([28]) and we may state:  \\

\noindent
{\bf DEFINITION 2.8}: An element $C\in Z_2({\cal{G}})$ is said to be {\it integrable} if there exists a deformation $c_t$ of $c$ such that $\frac{dc_t}{dt}{\mid}_{t=0}=C$. It is said to be {\it formally integrable} if there exists a formal power series $c_t=\sum_{\nu=0}^{\infty}\frac{t^{\nu}}{\nu !}C_{\nu}$ with $C_0=c, C_1=C$.ÊÊ\\ 

\noindent
{\bf EXAMPLE  2.9}: Even if Example 1.4 has nothing to do at first sight with any Lie algebra as we already said, we may adapt the previous arguments to $c'c''=0$. Indeed, we must have:\\
\[  0=(c'+tC'_1+\frac{t^2 }{2} C'_2+ ... )(c''+tC''_1+\frac{t^2}{2}C''_2+ ...)=t(c''C'_1+c'C''_1)+\frac{t^2}{2}(c''C'_2+2C'_1C''_1+c'C''_2) + ...    \]  
As cocycles are defined by the condition $c''C'+c'C''=0$, it follows that the cocycle $(1,1)$ cannot be integrable at $c=(0,0)$ as we should have $C'C''=0$.\\

In order to study the formal integrability of a cocycle $C\in Z_2({\cal{G}})$, we shall use a trick by introducing two parameters $s$ and $t$. Then $c_{s+t}$ becomes a deformation of $c_t$ in $s$ and we have:  \\
\[   c_{s+t}=\sum_{\nu=0}^{\infty}\frac{(s+t)^{\nu}}{\nu !}C_{\nu}=\sum_{\nu=0}^{\infty}\sum_{\lambda+\mu=\nu}\frac{s^{\lambda}}{\lambda !}\frac{t^{\mu}}{\mu !}C_{\nu}=c_t+s(\sum_{\nu=0}^{\infty}\frac{t^{\nu}}{\nu !}C_{\nu+1}) + ...  \]
where one may notice the change $C \rightarrow \sum_{\nu=0}^{\infty}\frac{t^{\nu}}{\nu !}C_{\nu+1}$.  \\

\noindent
{\bf PROPOSITION 2.10}: A sufficient condition for the formal integrability of any cocycle is $H_3({\cal{G}})=0$.Ê \\
\noindent
{\it Proof}: By definition we have $d=d(c)$ and we may set $d(t)=d(c_t)$ in order to obtain:  \\
\[    dC_1=0  \Rightarrow  d(t)\sum_{\nu=0}^{\infty}\frac{t^{\nu}}{\nu !}C_{\nu+1}=0   \]
Differentiating with respect to $t$ while using the Leibnitz formula, we get:  \\
\[    \sum_{\lambda,\mu\geq 0,\lambda+\mu=\nu }\frac{\nu !}{\lambda !\mu !} (\frac{d^{\lambda} d(t)}{dt^{\lambda}}){\mid}_{t=0}C_{\mu+1}=0  \]
or, in an equivalent way:\\
\[dC_{\nu+1}+\sum_{\lambda>0,\mu\geq 0,\lambda+\mu=\nu}\frac{\nu !}{\lambda !\mu !}(\frac{d^{\lambda}d(t)}{dt^{\lambda}}){\mid}_{t=0}C_{\mu+1}=0\]
As the first term on the left belongs to $B_3({\cal{G}})$, it just remains to prove that the sum on the right , which only depends on $C_1, ... ,C_{\nu}$, is in $Z_3({\cal{G}})$. Otherwise, as shown by the previous example, the study of these inductively related conditions may sometimes bring {\it obstructions} to the deformation at a certain order.\\
  Now we have $d\circ d\equiv 0 \Rightarrow d(t)\circ d(t)\equiv 0$ and thus again:\\
  \[     \sum_{\lambda,\mu\geq 0, \lambda+\mu=\nu}\frac{\nu !}{\lambda !\mu !}(\frac{d^{\lambda}d(t)}{dt^{\lambda}}){\mid}_{t=0}\circ (\frac{d^{\mu}d(t)}{dt^{\mu}}){\mid}_{t=0}\equiv 0   \]
  Accordingly, we obtain successively by using this identity (care to the minus sign):\\
 \[ - d\circ \sum_{\lambda >0,\mu\geq 0, \lambda+\mu=\nu}\frac{\nu !}{\lambda !\mu !}(\frac{d^{\lambda}d(t)}{dt^{\lambda}}){\mid}_{t=0}C_{\mu+1}\equiv  
 -\sum_{\lambda >0,\mu\geq 0, \lambda+\mu=\nu}\frac{\nu !}{\lambda !\mu !} d\circ (\frac{d^{\lambda}d(t)}{dt^{\lambda}}){\mid}_{t=0}C_{\mu+1}   \]
 \[  \equiv \sum_{\lambda >0,\mu\geq 0, \lambda+\mu=\nu}\frac{\nu !}{\lambda !\mu !}[ \sum_{\alpha >0,\beta\geq 0,\alpha+\beta=\lambda}\frac{\lambda !}{\alpha !\beta !}(\frac{d^{\alpha}d(t)}{dt^{\alpha}}){\mid}_{t=0}\circ (\frac{d^{\beta}d(t)}{dt^{\beta}}){\mid}_{t=0}C_{\mu+1} ]  \]
 \[  \equiv \sum_{\alpha=1}^{\nu}\frac{\nu !}{\alpha !(\nu-\alpha)!}(\frac{d^{\alpha}d(t)}{dt^{\alpha}}){\mid}_{t=0}\circ [\sum_{\mu\geq 0,\beta\geq 0,\beta+\mu=\nu-\alpha}\frac{(\nu-\alpha)!}{\beta !\mu !}(\frac{d^{\beta}d(t)}{dt^{\beta}}){\mid}_{t=0}C_{\mu+1}] \equiv 0  \]
 because we may take into account the integrability conditions up to order $\nu-1$ as $\alpha>0$.  \\
\hspace*{13cm}     Q.E.D.  \\

Studying the successive integrability conditions, we get:\\
$\nu=0$     \hspace{4cm}  $dC_1=0 $  \\
$\nu=1$     \hspace{4cm}  $ dC_2+\frac{{\partial}^2J(c)}{\partial c\partial c}C_1 C_1=0 $ \\
and we have the following result not evident at first sight:\\

\noindent
{\bf COROLLARY 2.11}: The Hessian of the Jacobi conditions provides a quadratic map:\\
\[  H_2({\cal{G}})  \rightarrow  H_3({\cal{G}}): C \rightarrow  \frac{{\partial}^2J(c)}{\partial c\partial c}C C   \]

\noindent
{\it Proof}: The previous proposition proves that this map takes $C\in Z_2({\cal{G}})$ to $Z_3({\cal{G}})$ and we just need to prove that it also takes $C\in B_2({\cal{G}})$ to $B_3({\cal{G}})$. For this, let $A\in V^*\otimes V$ be such that $C=dA$. We have:\\
\[ -(\frac{dd(t)}{dt}){\mid}_{t=0}C=-(\frac{dd(t)}{dt}){\mid}_{t=0}\circ dA=d\circ (\frac{dd(t)}{dt}){\mid}_{t=0}A  \]
as we wished.\\
\hspace*{13cm}    Q.E.D.   \\
\vspace*{5mm}  \\
\noindent
{\bf 3  VESSIOT STRUCTURE EQUATIONS }:   \\

If $X$ is a manifold, we denote as usual by $T=T(X)$ the {\it tangent bundle} of $X$, by $T^*=T^*(X)$ the {\it cotangent bundle}, by ${\wedge}^rT^*$ the {\it bundle of r-forms} and by $S_qT^*$ the {\it bundle of q-symmetric tensors}. More generally, let $\cal{E}$ be a {\it fibered manifold}, that is a manifold with local coordinates $(x^i,y^k)$ for $i=1,...,n$ and $k=1,...,m$ simply denoted by $(x,y)$, {\it projection} $\pi:{\cal{E}}\rightarrow X:(x,y)\rightarrow (x)$ and changes of local coordinates $\bar{x}=\varphi(x), \bar{y}=\psi(x,y)$. If $\cal{E}$ and $\cal{F}$ are two fibered manifolds over $X$ with respective local coordinates $(x,y)$ and $(x,z)$, we denote by ${\cal{E}}{\times}_X{\cal{F}}$ the {\it fibered product} of $\cal{E}$ and $\cal{F}$ over $X$ as the new fibered manifold over $X$ with local coordinates $(x,y,z)$. We denote by $f:X\rightarrow {\cal{E}}: (x)\rightarrow (x,y=f(x))$ a global {\it section} of $\cal{E}$, that is a map such that $\pi\circ f=id_X$ but local sections over an open set $U\subset X$ may also be considered when needed. Under a change of coordinates, a section transforms like $\bar{f}(\varphi(x))=\psi(x,f(x))$ and the derivatives transform like:\\
\[   \frac{\partial{\bar{f}}^l}{\partial{\bar{x}}^r}(\varphi(x)){\partial}_i{\varphi}^r(x)=\frac{\partial{\psi}^l}{\partial x^i}(x,f(x))+\frac{\partial {\psi}^l}{\partial y^k}(x,f(x)){\partial}_if^k(x)  \]
We may introduce new coordinates $(x^i,y^k,y^k_i)$ transforming like:\\
\[ {\bar{y}}^l_r{\partial}_i{\varphi}^r(x)=\frac{\partial{\psi}^l}{\partial x^i}(x,y)+\frac{\partial {\psi}^l}{\partial y^k}(x,y)y^k_i  \]
We shall denote by $J_q({\cal{E}})$ the {\it q-jet bundle} of $\cal{E}$ with local coordinates $(x^i, y^k, y^k_i, y^k_{ij},...)=(x,y_q)$ called {\it jet coordinates} and sections $f_q:(x)\rightarrow (x,f^k(x), f^k_i(x), f^k_{ij}(x), ...)=(x,f_q(x))$ transforming like the sections $j_q(f):(x) \rightarrow (x,f^k(x), {\partial}_if^k(x), {\partial}_{ij}f^k(x), ...)=(x,j_q(f)(x))$ where both $f_q$ and $j_q(f)$ are over the section $f$ of $\cal{E}$. Of course $J_q({\cal{E}})$ is a fibered manifold over $X$ with projection ${\pi}_q$ while $J_{q+r}({\cal{E}})$ is a fibered manifold over $J_q({\cal{E}})$ with projection ${\pi}^{q+r}_q, \forall r\geq 0$.\\

\noindent
{\bf DEFINITION 3.1}: A {\it system} of order $q$ on $\cal{E}$ is a fibered submanifold ${\cal{R}}_q\subset J_q({\cal{E}})$ and a {\it solution} of ${\cal{R}}_q$ is a section $f$ of $\cal{E}$ such that $j_q(f)$ is a section of ${\cal{R}}_q$.\\

\noindent
{\bf DEFINITION 3.2}: When the changes of coordinates have the linear form $\bar{x}=\varphi(x),\bar{y}= A(x)y$, we say that $\cal{E}$ is a {\it vector bundle} over $X$ and denote for simplicity a vector bundle and its set of sections by the same capital letter $E$. When the changes of coordinates have the form $\bar{x}=\varphi(x),\bar{y}=A(x)y+B(x)$ we say that $\cal{E}$ is an {\it affine bundle} over $X$ and we define the {\it associated vector bundle} $E$ over $X$ by the local coordinates $(x,v)$ changing like $\bar{x}=\varphi(x),\bar{v}=A(x)v$.\\

\noindent
{\bf DEFINITION 3.3}: If the tangent bundle $T({\cal{E}})$ has local coordinates $(x,y,u,v)$ changing like ${\bar{u}}^j={\partial}_i{\varphi}^j(x)u^i, {\bar{v}}^l=\frac{\partial {\psi}^l}{\partial x^i}(x,y)u^i+\frac{\partial {\psi}^l}{\partial y^k}(x,y)v^k$, we may introduce the {\it vertical bundle} $V({\cal{E}})\subset T({\cal{E}})$ as a vector bundle over $\cal{E}$ with local coordinates $(x,y,v)$ obtained by setting $u=0$ and changes ${\bar{v}}^l=\frac{\partial {\psi}^l}{\partial y^k}(x,y)v^k$. Of course, when $\cal{E}$ is an affine bundle over $X$ with associated vector bundle $E$ over $X$, we have $V({\cal{E}})={\cal{E}}\times_XE$.\\

For a later use, if $\cal{E}$ is a fibered manifold over $X$ and $f$ is a section of $\cal{E}$, we denote by $f^{-1}(V({\cal{E}}))$ the {\it reciprocal image} of $V({\cal{E}})$ by $f$ as the vector bundle over $X$ obtained when replacing $(x,y,v)$ by $(x,f(x),v) $ in each chart. A similar construction may also be done for any affine bundle over ${\cal{E}}$.  \\

We now recall a few basic geometric concepts that will be constantly used through this paper. First of all, if $\xi,\eta\in T$, we define their {\it bracket} $[\xi,\eta]\in T$ by the local formula $([\xi,\eta])^i(x)={\xi}^r(x){\partial}_r{\eta}^i(x)-{\eta}^s(x){\partial}_s{\xi}^i(x)$ leading to the {\it Jacobi identity} $[\xi,[\eta,\zeta]]+[\eta,[\zeta,\xi]]+[\zeta,[\xi,\eta]]=0, \forall \xi,\eta,\zeta \in T$ allowing to define a {\it Lie algebra} and to the useful formula $[T(f)(\xi),T(f)(\eta)]=T(f)([\xi,\eta])$ where $T(f):T(X)\rightarrow T(Y)$ is the tangent mapping of a map $f:X\rightarrow Y$.\\

When $I=\{ i_1< ... < i_r\}$ is a multi-index, we may set $dx^I=dx^{i_1}\wedge ... \wedge dx^{i_r}$ for describing ${\wedge}^rT^*$ and introduce the {\it exterior derivative} $d:{\wedge}^rT^*\rightarrow {\wedge}^{r+1}T^*:\omega={\omega}_Idx^I \rightarrow d\omega={\partial}_i{\omega}_Idx^i\wedge dx^I$ with $d^2=d\circ d\equiv 0$ in the {\it Poincar\'{e} sequence}:\\
\[  {\wedge}^0T^* \stackrel{d}{\longrightarrow} {\wedge}^1T^* \stackrel{d}{\longrightarrow} {\wedge}^2T^* \stackrel{d}{\longrightarrow} ... \stackrel{d}{\longrightarrow} {\wedge}^nT^* \longrightarrow 0  \]

The {\it Lie derivative} of an $r$-form with respect to a vector field $\xi\in T$ is the linear first order operator ${\cal{L}}(\xi)$ linearly depending on $j_1(\xi)$ and uniquely defined by the following three properties:\\
1) ${\cal{L}}(\xi)f=\xi.f={\xi}^i{\partial}_if, \forall f\in {\wedge}^0T^*=C^{\infty}(X)$.\\
2) ${\cal{L}}(\xi)d=d{\cal{L}}(\xi)$.\\
3) ${\cal{L}}(\xi)(\alpha\wedge \beta)=({\cal{L}}(\xi)\alpha)\wedge \beta+\alpha\wedge ({\cal{L}}(\xi) \beta), \forall \alpha,\beta \in \wedge T^*$.\\
It can be proved that ${\cal{L}}(\xi)=i(\xi)d+di(\xi)$ where $i(\xi)$ is the {\it interior multiplication} $(i(\xi)\omega)_{i_1...i_r}={\xi}^i{\omega}_{ii_1...i_r}$ and that $[{\cal{L}}(\xi),{\cal{L}}(\eta)]={\cal{L}}(\xi)\circ {\cal{L}}(\eta)-{\cal{L}}(\eta)\circ {\cal{L}}(\xi)={\cal{L}}([\xi,\eta]), \forall \xi,\eta\in T$.\\

We now turn to group theory and start with two basic definitions:\\

Let $G$ be a {\it Lie group}, that is another manifold with local coordinates $a=(a^1, ... , a^p)$ called {\it parameters}, a {\it composition} $G\times G \rightarrow G: (a,b)\rightarrow ab$, an {\it inverse} $G \rightarrow G: a \rightarrow a^{-1}$ and an {\it identity} $e\in G$ satisfying:\\
\[(ab)c=a(bc)=abc,\hspace{1cm} aa^{-1}=a^{-1}a=e,\hspace{1cm} ae=ea=a,\hspace{1cm} \forall a,b,c \in G \]

\noindent
{\bf DEFINITION 3.4}: $G$ is said to {\it act} on $X$ if there is a map $X\times G \rightarrow X: (x,a) \rightarrow y=ax=f(x,a)$ such that $(ab)x=a(bx)=abx, \forall a,b\in G, \forall x\in X$ and we shall say that we have a {\it Lie group of transformations} of $X$. In order to simplify the notations, we shall use global notations even if only local actions are existing. The set $G_x=\{a\in G\mid ax=x\}$ is called the {\it isotropy subgroup} of $G$ at $x\in X$ and the action is said to be {\it effective} if $ax=x, \forall x\in X\Rightarrow a=e$. \\ 

\noindent
{\bf DEFINITION 3.5}: A {\it Lie pseudogroup of transformations} $\Gamma\subset aut(X)$ is a group of transformations solutions of a system of OD or PD equations such that, if $y=f(x)$ and $z=g(y)$ are two solutions, called {\it finite transformations}, that can be composed, then $z=g\circ f(x)=h(x)$ and $x=f^{-1}(y)=g(y)$ are also solutions while $y=x$ is a solution and we shall set $id_q=j_q(id)$.  \\

It becomes clear from Examples 1.1 and 1.2 that Lie groups of transformations are particular cases of Lie pseudogroups of transformations as the system defining the finite transformations can be obtained by eliminating the parameters among the equations $y_q=j_q(f)(x,a)$ when $q$ is large enough. The underlying system may be nonlinear and of high order as we have seen. We shall speak of an {\it algebraic pseudogroup} when the system is defined by {\it differential polynomials} that is polynomials in the derivatives. Looking for transformations "close" to the identity, that is setting $y=x+t\xi(x)+...$ when $t\ll 1$ is a small constant parameter and passing to the limit $t\rightarrow 0$, we may linearize the above nonlinear {\it system of finite Lie equations} in order to obtain a linear {\it system of infinitesimal Lie equations} of the same order for vector fields. Such a system has the property that, if $\xi,\eta$ are two solutions, then $[\xi,\eta]$ is also a solution. Accordingly, the set $\Theta\subset T$ of solutions of this new system satisfies $[\Theta,\Theta]\subset \Theta$ and can therefore be considered as the Lie algebra of $\Gamma$.\\

We now turn to the theory proposed by Vessiot in 1903 ([51]) and sketch in a few successive steps the main results that we have obtained in many books ([36],[37],[38],[39]). We invite the reader to follow the procedure on each of the examples provided for this purpose in the introduction.  \\

1) If ${\cal{E}}=X\times X$, we shall denote by ${\Pi}_q={\Pi}_q(X,X)$ the open subfibered manifold of $J_q(X\times X)$ defined independently of the coordinate system by $det(y^k_i)\neq 0$ with {\it source projection} ${\alpha}_q:{\Pi}_q\rightarrow X:(x,y_q)\rightarrow (x)$ and {\it target projection} ${\beta}_q:{\Pi}_q\rightarrow X:(x,y_q)\rightarrow (y)$. We shall sometimes introduce a copy $Y$ of $X$ with local coordinates $(y)$ in order to avoid any confusion between the source and the target manifolds. Let us start with a Lie pseudogroup $\Gamma\subset aut(X)$ defined by a system ${\cal{R}}_q\subset {\Pi}_q$ of order $q$. In all the sequel we shall suppose that the system is {\it involutive} (see next section) and that $\Gamma$ is {\it transitive} that is $\forall x,y\in X, \exists f\in \Gamma, y=f(x)$ or, equivalently, the map $({\alpha}_q,{\beta}_q):{\cal{R}}_q\rightarrow X\times X:(x,y_q)\rightarrow (x,y)$ is surjective.\\

2) The Lie algebra $\Theta\subset T$ of infinitesimal transformations is then obtained by linearization, setting $y=x+t\xi(x)+...$ and passing to the limit $t\rightarrow 0$ in order to obtain the linear involutive system $R_q= id^{-1}_q(V({\cal{R}}_q))\subset J_q(T)$ by reciprocal image with $\Theta=\{\xi \in T{\mid}j_q(\xi) \in R_q\}$ while taking into account the fact that $T=id^{-1}(V(X\times X))$. From now on we shall suppose that $R_q$ is {\it transitive}, that is to say the canonical projection ${\pi}^q_0:J_q(T) \rightarrow T$ induces an epimorphism ${\pi}^q_0:R_q \rightarrow T$ with kernel $R^0_q\subset R_q$ and we have the useful short exact sequence $0\longrightarrow R^0_q \longrightarrow R_q \stackrel{{\pi}^q_0}{\longrightarrow} T \longrightarrow 0 $.  \\

3) Passing from source to target, we may {\it prolong} the vertical infinitesimal transformations $\eta={\eta}^k(y)\frac{\partial}{\partial y^k}$ to the jet coordinates up to order $q$ in order to obtain:\\
\[   {\eta}^k(y)\frac{\partial}{\partial y^k}+\frac{\partial {\eta}^k}{\partial y^r}y^r_i\frac{\partial}{\partial y^k_i}+(\frac{{\partial}^2{\eta}^k}{\partial y^r\partial y^s}y^r_iy^s_j+\frac{\partial {\eta}^k}{\partial y^r}y^r_{ij})\frac{\partial}{\partial y^k_{ij}}+...    \]
where we have replaced $j_q(f)(x)$ by $y_q$, each component beeing the "formal" derivative of the previous one obtained by introducing $d_i={\partial}_i+y^k_{\mu +1_i}\frac{\partial}{\partial y^k_{\mu}}$.\\

4) As $[\Theta,\Theta]\subset \Theta$, we may use the Frobenius theorem in order to find a generating fundamental set of {\it differential invariants} $\{{\Phi}^{\tau}(y_q)\}$ up to order $q$ which are such that ${\Phi}^{\tau}({\bar{y}}_q)={\Phi}^{\tau}(y_q)$ by using the chain rule for derivatives whenever $\bar{y}=g(y)\in \Gamma$ acting now on $Y$. Of course, in actual practice {\it one must use sections of} $R_q$ {\it instead of solutions} but it is only in section 5 through Definition 5.2 that we shall see why the use of the Spencer operator will be crucial for this purpose. Specializing the ${\Phi}^{\tau}$ at $id_q(x)$ we obtain the {\it Lie form} ${\Phi}^{\tau}(y_q)={\omega}^{\tau}(x)$ of ${\cal{R}}_q$. Finally, if ${\Phi}^{\tau}$ is any differential invariant at the order $q$, then $d_i{\Phi}^{\tau}$ is a differential invariant at order $q+1, \forall i=1,...,n$.\\

5) The main discovery of Vessiot, fifty years in advance, has been to notice that the prolongation at order $q$  of any horizontal vector field $\xi={\xi}^i(x)\frac{\partial}{\partial x^i}$ commutes with the prolongation at order $q$ of any vertical vector field $\eta={\eta}^k(y)\frac{\partial}{\partial y^k}$, exchanging therefore the differential invariants. Keeping in mind the well known property of the Jacobian determinant while passing to the finite point of view, any (local) transformation $y=f(x)$ can be lifted to a (local) transformation of the differential invariants between themselves of the form $u\rightarrow \lambda(u,j_q(f)(x))$ allowing to introduce a {\it natural bundle} $\cal{F}$ over $X$ by patching changes of coordinates $\bar{x}=\varphi(x), \bar{u}=\lambda(u,j_q(\varphi)(x))$. A section $\omega$ of $\cal{F}$ is called a {\it geometric object} or {\it structure} on $X$ and transforms like ${\bar{\omega}}(f(x))=\lambda(\omega(x),j_q(f)(x))$ or simply $\bar{\omega}=j_q(f)(\omega)$. This is a way to generalize vectors and tensors ($q=1$) or even connections ($q=2$). As a byproduct we have $\Gamma=\{f\in aut(X){\mid} {\Phi}_{\omega}(j_q(f))\equiv j_q(f)^{-1}(\omega)=\omega\}$ as a new way to write out the Lie form and we may say that $\Gamma$ {\it preserves} $\omega$. Replacing $j_q(f)$ by $f_q$, we also obtain ${\cal{R}}_q=\{f_q\in {\Pi}_q{\mid} f_q^{-1}(\omega)=\omega\}$. Coming back to the infinitesimal point of view and setting $f_t=exp(t\xi)\in aut(X), \forall \xi\in T$, we may define the {\it ordinary Lie derivative} with value in ${\omega}^{-1}(V({\cal{F}}))$ by the formula :\\
\[   {\cal{D}}\xi={\cal{D}}_{\omega}\xi={\cal{L}}(\xi)\omega=\frac{d}{dt}j_q(f_t)^{-1}(\omega){\mid}_{t=0} \Rightarrow \Theta=\{\xi\in T{\mid}{\cal{L}}(\xi)\omega=0\}      \]
We have $x\rightarrow \bar{x}=x+t\xi(x)+...\Rightarrow u^{\tau}\rightarrow {\bar{u}}^{\tau}=u^{\tau}+t{\partial}_{\mu}{\xi}^kL^{\tau\mu}_k(u)+...$ where $\mu=({\mu}_1,...,{\mu}_n)$ is a multi-index and we 
may write down the system of infinitesimal Lie equations in the {\it Medolaghi form}:\\
\[     {\Omega}^{\tau}\equiv ({\cal{L}}(\xi)\omega)^{\tau}\equiv -L^{\tau\mu}_k(\omega(x)){\partial}_{\mu}{\xi}^k+{\xi}^r{\partial}_r{\omega}^{\tau}(x)=0    \]
as a way to state the invariance of the section $\omega$ of ${\cal{F}}$, that is $u^{\tau}-{\omega}^{\tau}(x)=0\Rightarrow {\bar{u}}^{\tau}-{\bar{\omega}}^{\tau}(\bar{x})=0$.  \\
Finally, replacing $j_q(\xi)$ by a section ${\xi}_q\in J_q(T)$ over $\xi\in T$, we may define $R_q\subset J_q(T)$ {\it on sections} by the linear (non-differential) equations:ÊÊ\\
\[  {\Omega}^{\tau}\equiv (L({\xi}_q)\omega)^{\tau}\equiv - L^{\tau\mu}_k(\omega (x)){\xi}^k_{\mu}+ {\xi}^r{\partial}_r{\omega}^{\tau}(x)=0  \]
and obtain the first prolongation $R_{q+1}\subset J_{q+1}(T)$ by adding:   \\
\[ {\Omega}^{\tau}_i\equiv (L({\xi}_{q+1})j_1(\omega))^{\tau}_i\equiv -L^{\tau\mu}_k(\omega (x)){\xi}^k_{\mu+1_i}-\frac{\partial L^{\tau\mu}_k(\omega (x))}{\partial u^{\sigma}}{\partial}_i{\omega}^{\sigma}(x){\xi}^k_{\mu}+{\partial}_r{\omega}^{\tau}(x){\xi}^r_i+{\xi}^r{\partial}_r({\partial}_i{\omega}^{\tau}(x))=0   \]

6) By analogy with "special" and "general" relativity, we shall call the given section {\it special} and any other arbitrary section {\it general}. The problem is now to study the formal properties of the linear system just obtained with coefficients only depending on $j_1(\omega)$, exactly like we did in the examples of the introduction. In particular, if any expression involving $\omega$ and its derivatives is a scalar object, it must reduce to a constant because $\Gamma$ is assumed to be transitive and thus cannot be defined by any zero order equation. Now one can prove that the CC for $\bar{\omega}$, thus for $\omega$ too, only depend on the $\Phi$ and take the quasi-linear symbolic form $v\equiv I(u_1)\equiv A(u)u_x+B(u)=0$ with $u_1=(u,u_x)$, allowing to define an affine subfibered manifold ${\cal{B}}_1\subset J_1({\cal{F}})$ over $\cal{F}$. Indeed, if $\sum A(y_q)d_x\Phi$ is a minimum sum of formal derivatives of differential invariants of order $q$ not containing any jet coordinate of strict order $q+1$, we may suppose by division that the first $A$ in the sum is equal to $1$. Applying the prolonged distribution of vector fields introduced in the step 3 at order $q+1$, we obtain a new sum with less terms and a contradiction unless all the $A$ are again differential invariants at order $q$ and thus functions of the $\Phi$ because of the Frobenius theorem. A similar comment can be done for the $B$. Now, if one has two sections $\omega$ and $\bar{\omega}$ of $\cal{F}$, the {\it equivalence problem} is to look for $f\in aut(X)$ such that $j_q(f)^{-1}(\omega)=\bar{\omega}$. When the two sections satisfy the same CC, the problem is sometimes locally possible (Lie groups of transformations, Darboux problem in analytical mechanics,...) but sometimes not ([36], p 333).\\

7) Instead of the CC for the equivalence problem, let us look for the {\it integrability conditions} (IC) for the system of infinitesimal Lie equations and suppose that, for the given section, all the equations of order $q+r$ are obtained by differentiating $r$ times only the equations of order $q$, then it was claimed by Vessiot ([51] with no proof, see [36], p 313,[39], p 207-212) that such a property is held if and only if there is an equivariant section $c:{\cal{F}}\rightarrow {\cal{F}}_1:(x,u)\rightarrow (x,u,v=c(u))$ where ${\cal{F}}_1=J_1({\cal{F}})/{\cal{B}}_1$ is a natural vector bundle over $\cal{F}$ with local coordinates $(x,u,v)$. Moreover, any such equivariant section depends linearly on a finite number of constants $c$ called {\it structure constants} and the IC for the {\it Vessiot structure equations} $I(u_1)=c(u)$ are of a polynomial form $J(c)=0$. It is important to notice that {\it the form of the Vessiot structure equations is invariant under any change of coordinate system}. In actual practice, this study can be divided into two parts according to the following commutative and exact diagram: \\
\[  \begin{array}{rcccccl}
  & 0 & &  0 & & &  \\
  & \downarrow & & \downarrow & & &   \\
  0\longrightarrow  & g_{q+1} &  =  &g_{q+1} &\longrightarrow & 0 &  \\
         & \downarrow & & \downarrow &  & \downarrow &   \\
  0\longrightarrow &  R^0_{q+1} & \longrightarrow & R_{q+1}  & \longrightarrow &  T  &     \\
     &  \downarrow &  & \downarrow &  & \parallel  &    \\
       0 \longrightarrow &  R^0_q & \longrightarrow & R_q  & \longrightarrow &  T  &  \longrightarrow 0  \\
   &  &  &  &  &  \downarrow &  \\
   &  &  &  &  &   0  &   
   \end{array}    \]
 \noindent
 Indeed, chasing in this diagram, we discover that ${\pi}^{q+1}_q:R_{q+1}\rightarrow R_q $ is an epimorphism if and only if ${\pi}^{q+1}_q:R^0_{q+1} \rightarrow R^0_q$ is an epimorphism {\it and} ${\pi}^{q+1}_0:R_{q+1} \rightarrow T$ is also an epimorphism. Looking at the form of the corresponding Medolaghi equations $L({\xi}_q)\omega=0$ and $L({\xi}_{q+1})j_1(\omega)=0$, these two conditions respectively bring the Vessiot structure equations of {\it first kind} $I_{\ast}(u_1)=0$ not invoving any structure constant and the Vessiot stucture equations of {\it second kind} $I_{\ast\ast}(u_1)=c$ involving the same number of different structure constants. Such a study, only depending now on linear algebraic techniques, can be achieved by means of computer algebra ([3],[30]).  \\
Finally, looking at the formal integrability of the system ${\cal{B}}_1 \subset J_1({\cal{F}})$ defined by the equations $A(u)u_x+B(u)=0$ and their first prolongation $A(u)u_{xx}+\frac{\partial A(u)}{\partial u}u_xu_x+\frac{\partial B(u)}{\partial u}u_x=0$, the symbols only depend on $A(u)$ and we may obtain equations of the form $a(u)u_xu_x+b(u)u_x=0$ by eliminating the jets of order $2$. Using local coordinates $(x,u,v=A(u)u_x+B(u))$ for ${\cal{F}}_1$, and substituting $(u,v)$ in place of $(u,u_x)$, we obtain equations of the form $\alpha (u)vv+\beta (u)v+\gamma (u)=0$. As we may suppose that $c=0$ for the special section, we finally get equations of the form $\alpha (u) I I + \beta (u) I=0$ and it only remains to set  $I(u_1)=c(u)$ in order to get polynomial Jacobi conditions of degree $\leq 2$ which may not depend on $u$ anymore because these equations are invariant in form under any change of coordinates.  \\

\noindent
{\bf REMARK  3.6}: When $q=1$, a close examination of the Medolaghi equations and their first prolongation shows at once that we can choose $v=A(u)u_x$ and we get homogeneous Jacobi conditions of degree $2$. Such a result explains the difference existing between Examples 1.1, 1.3, 1.4 ($q=1$) and Example 1.2 ($q=2$).  \\
\vspace*{5mm}   \\
\noindent
{\bf 4  LINEAR AND NONLINEAR JANET SEQUENCES}:   \\

Let $\mu=({\mu}_1,...,{\mu}_n)$ be a multi-index with {\it length} ${\mid}\mu{\mid}={\mu}_1+...+{\mu}_n$, {\it class} $i$ if ${\mu}_1=...={\mu}_{i-1}=0,{\mu}_i\neq 0$ and $\mu +1_i=({\mu}_1,...,{\mu}_{i-1},{\mu}_i +1, {\mu}_{i+1},...,{\mu}_n)$. We set $y_q=\{y^k_{\mu}{\mid} 1\leq k\leq m, 0\leq {\mid}\mu{\mid}\leq q\}$ with $y^k_{\mu}=y^k$ when ${\mid}\mu{\mid}=0$. If $E$ is a vector bundle over $X$ with local coordinates $(x^i,y^k)$ for $i=1,...,n$ and $k=1,...,m$, we denote by $J_q(E)$ the $q$-{\it jet bundle} of $E$ with local coordinates simply denoted by $(x,y_q)$ and {\it sections} $f_q:(x)\rightarrow (x,f^k(x), f^k_i(x), f^k_{ij}(x),...)$ transforming like the section $j_q(f):(x)\rightarrow (x,f^k(x),{\partial}_if^k(x),{\partial}_{ij}f^k(x),...)$ when $f$ is an arbitrary section of $E$. Then both $f_q\in J_q(E)$ and $j_q(f)\in J_q(E)$ are over $f\in E$ and the {\it Spencer operator} just allows to distinguish them by introducing a kind of "{\it difference}" through the operator $D:J_{q+1}(E)\rightarrow T^*\otimes J_q(E): f_{q+1}\rightarrow j_1(f_q)-f_{q+1}$ with local components $({\partial}_if^k(x)-f^k_i(x), {\partial}_if^k_j(x)-f^k_{ij}(x),...) $ and more generally $(Df_{q+1})^k_{\mu,i}(x)={\partial}_if^k_{\mu}(x)-f^k_{\mu+1_i}(x)$. In a symbolic way, {\it when changes of coordinates are not involved}, it is sometimes useful to write down the components of $D$ in the form $d_i={\partial}_i-{\delta}_i$ and the restriction of $D$ to the kernel $S_{q+1}T^*\otimes E$ of the canonical projection ${\pi}^{q+1}_q:J_{q+1}(E)\rightarrow J_q(E)$ is {\it minus} the {\it Spencer map} $\delta=dx^i\wedge {\delta}_i:S_{q+1}T^*\otimes E\rightarrow T^*\otimes S_qT^*\otimes E$. The kernel of $D$ is made by sections such that $f_{q+1}=j_1(f_q)=j_2(f_{q-1})=...=j_{q+1}(f)$. Finally, if $R_q\subset J_q(E)$ is a {\it system} of order $q$ on $E$ locally defined by linear equations ${\Phi}^{\tau}(x,y_q)\equiv a^{\tau\mu}_k(x)y^k_{\mu}=0$ and local coordinates $(x,z)$ for the parametric jets up to order $q$, the $r$-{\it prolongation} $R_{q+r}={\rho}_r(R_q)=J_r(R_q)\cap J_{q+r}(E)\subset J_r(J_q(E))$ is locally defined when $r=1$ by the linear equations ${\Phi}^{\tau}(x,y_q)=0, d_i{\Phi}^{\tau}(x,y_{q+1})\equiv a^{\tau\mu}_k(x)y^k_{\mu+1_i}+{\partial}_ia^{\tau\mu}_k(x)y^k_{\mu}=0$ and has {\it symbol} $g_{q+r}=R_{q+r}\cap S_{q+r}T^*\otimes E\subset J_{q+r}(E)$ if one looks at the {\it top order terms}. If $f_{q+1}\in R_{q+1}$ is over $f_q\in R_q$, differentiating the identity $a^{\tau\mu}_k(x)f^k_{\mu}(x)\equiv 0$ with respect to $x^i$ and substracting the identity $a^{\tau\mu}_k(x)f^k_{\mu+1_i}(x)+{\partial}_ia^{\tau\mu}_k(x)f^k_{\mu}(x)\equiv 0$, we obtain the identity $a^{\tau\mu}_k(x)({\partial}_if^k_{\mu}(x)-f^k_{\mu+1_i}(x))\equiv 0$ and thus the restriction $D:R_{q+1}\rightarrow T^*\otimes R_q$. This first order operator induces, up to sign, the purely algebraic monomorphism $ 0 \rightarrow g_{q+1} \stackrel{\delta}{\rightarrow} T^*\otimes g_q$ on the symbol level ([36],[40],[49]). \\
    
\noindent
{\bf DEFINITION} 4.1: $R_q$ is said to be {\it formally integrable} when the restriction ${\pi}^{q+r+1}_{q+r}:R_{q+r+1}\rightarrow R_{q+r} $ is an epimorphism $\forall r\geq 0$ or, equivalently, when all the equations of order $q+r$ are obtained by $r$ prolongations only, $\forall r\geq 0$. In that case, $R_{q+1}\subset J_1(R_q)$ is a canonical equivalent formally integrable first order system on $R_q$ with no zero order equations, called the {\it Spencer form}.\\

\noindent
{\bf DEFINITION} 4.2: $R_q$ is said to be {\it involutive} when it is formally integrable and the symbol $g_q$ is {\it involutive}, that is all the sequences $... \stackrel{\delta}{\rightarrow} {\wedge}^sT^*\otimes g_{q+r}\stackrel{\delta}{\rightarrow}...$ are exact $\forall 0\leq s\leq n, \forall r\geq 0$. Equivalently, using a linear change of local coordinates if necessary, we may {\it successively} solve the maximum number ${\beta}^n_q, {\beta}^{n-1}_q, ... , {\beta}^1_q$ of equations with respect to the principal jet coordinates of strict order $q$ and class $n,n-1,...,1$ in order to introduce the {\it characters} ${\alpha}^i_q=m\frac{(q+n-i-1)!}{(q-1)!((n-i)!}-{\beta}^i_q$ for $i=1, ..., n$ with ${\alpha}^n_q=\alpha$. Then $R_q$ is involutive if $R_{q+1}$ is obtained by only prolonging the ${\beta}^i_q$ equations of class $i$ with respect to $d_1,...,d_i$ for $i=1,...,n$. In that case $dim(g_{q+1})={\alpha}^1_q+...+{\alpha}^n_q$ and one can exhibit the {\it Hilbert polynomial} $dim(R_{q+r})$ in $r$ with leading term $(\alpha/n!)r^n$ when $\alpha \neq 0$. Such a prolongation procedure allows to compute {\it in a unique way} the principal ($pri$) jets from the parametric ($par$) other ones ([23],[46]). \\

\noindent
{\bf REMARK 4.3 }: This definition may also be applied to nonlinear systems as well while using a generic linearization by means of vertical bundles. For example, with $m=1$, $n=2$ and $q=2$, the nonlinear system $y_{22}-\frac{1}{3}(y_{11})^3=0, y_{12}-\frac{1}{2}(y_{11})^2=0$ is involutive but the nonlinear system $y_{22}-\frac{1}{2}(y_{11})^2=0, y_{12}-y_{11}=0$ is not involutive.  \\

    When $R_q$ is involutive, the linear differential operator ${\cal{D}}:E\stackrel{j_q}{\rightarrow} J_q(E)\stackrel{\Phi}{\rightarrow} J_q(E)/R_q=F_0$ of order $q$ with space of solutions $\Theta\subset E$ is said to be {\it involutive} and one has the canonical {\it linear Janet sequence} ([39], p 144):\\
\[  0 \longrightarrow  \Theta \longrightarrow T \stackrel{\cal{D}}{\longrightarrow} F_0 \stackrel{{\cal{D}}_1}{\longrightarrow}F_1 \stackrel{{\cal{D}}_2}{\longrightarrow} ... \stackrel{{\cal{D}}_n}{\longrightarrow} F_n \longrightarrow 0   \]
where each other operator is first order involutive and generates the {\it compatibility conditions} (CC) of the preceding one. As the Janet sequence can be cut at any place, {\it the numbering of the Janet bundles has nothing to do with that of the Poincar\'{e} sequence}, contrary to what many 
people believe.\\

    Equivalently, we have the involutive {\it first Spencer operator} $D_1:C_0=R_q\stackrel{j_1}{\rightarrow}J_1(R_q)\rightarrow J_1(R_q)/R_{q+1}\simeq T^*\otimes R_q/\delta (g_{q+1})=C_1$ of order one induced by $D:R_{q+1}\rightarrow T^*\otimes R_q$. Introducing the {\it Spencer bundles} $C_r={\wedge}^rT^*\otimes R_q/{\delta}({\wedge}^{r-1}T^*\otimes g_{q+1})$, the first order involutive ($r+1$)-{\it Spencer operator} $D_{r+1}:C_r\rightarrow C_{r+1}$ is induced by $D:{\wedge}^rT^*\otimes R_{q+1}\rightarrow {\wedge}^{r+1}T^*\otimes R_q:\alpha\otimes {\xi}_{q+1}\rightarrow d\alpha\otimes {\xi}_q+(-1)^r\alpha\wedge D{\xi}_{q+1}$ and we obtain the canonical {\it linear Spencer sequence} ([39], p 150):\\
\[    0 \longrightarrow \Theta \stackrel{j_q}{\longrightarrow} C_0 \stackrel{D_1}{\longrightarrow} C_1 \stackrel{D_2}{\longrightarrow} C_2 \stackrel{D_3}{\longrightarrow} ... \stackrel{D_n}{\longrightarrow} C_n\longrightarrow 0  \]
\noindent
as the Janet sequence for the first order involutive system $R_{q+1}\subset J_1(R_q)$. Introducing the other Spencer bundles $C_r(E)={\wedge}^rT^*\otimes J_q(E)/\delta({\wedge}^{r-1}T^*\otimes S_{q+1}T^*\otimes  E)$ with $C_r \subset C_r(E)$, the linear Spencer sequence is induced by the {\it linear hybrid sequence}: \\
\[  0 \longrightarrow E \stackrel{j_q}{\longrightarrow} C_0(E) \stackrel{D_1}{\longrightarrow} C_1(E) \stackrel{D_2}{\longrightarrow} C_2 \stackrel{D_3}{\longrightarrow} ... \stackrel{D_n}{\longrightarrow } C_n \longrightarrow 0 \]
\noindent
which is at the same time the Janet sequence for $j_q$ and the Spencer sequence for $J_{q+1}(E)\subset J_1(J_q(E))$ ([39], p 153):\\

We have the following commutative and exact diagram allowing to relate the Spencer bundles $C_r$ and $C_r(E)$ to the {\it Janet bundles} $F_r={\wedge}^rT^*\otimes F_0/\delta({\wedge}^{r-1}T^*\otimes h_1)$ if we start with the short exact sequence $ 0 \rightarrow g_{q+1} \rightarrow S_{q+1}T^*\otimes E \rightarrow h_1  \rightarrow 0  $ where $h_1\subset T^*\otimes F_0$:ÊÊÊ\\

\[  \begin{array}{cccccl}
  0  &   &  0  &  &  0  &   \\
  \downarrow &  &  \downarrow  &   &  \downarrow  &   \\
  {\wedge}^{r-1}T^*\otimes g_{q+1}  &  \stackrel{\delta}{\longrightarrow} &{\wedge}^rT^*\otimes R_q  &\longrightarrow & C_r  &  \longrightarrow 0  \\
  \downarrow &  &  \downarrow  &   &  \downarrow  &   \\
  {\wedge}^{r-1}T^*\otimes S_{q+1}T^*\otimes E &  \stackrel{\delta}{\longrightarrow} & {\wedge}^rT^*\otimes J_q(E) & \longrightarrow & C_r(E) & \longrightarrow 0  \\
  \downarrow  &    &  \hspace{4mm}\downarrow \Phi &  & \hspace{5mm}\downarrow {\Phi}_r  &   \\
  {\wedge}^{r-1}T^*\otimes h_1  &  \stackrel{\delta}{\longrightarrow} & {\wedge}^rT^*\otimes F_0 & \longrightarrow  & F_r &  \longrightarrow  0  \\
  \downarrow  &  & \downarrow &  & \downarrow   &   \\
  0  &  &  0  &  &  0  &   
  \end{array}   \]
\noindent
In this diagram, only depending on the linear differential operator ${\cal{D}}=\Phi\circ j_q$, the epimorhisms ${\Phi}_r:C_r(E)\rightarrow F_r$ for $0\leq r \leq n$ are induced by the canonical projection $\Phi={\Phi}_0:C_0(E)=J_q(E)\rightarrow J_q(E)/R_q=F_0$ if we start with the knowledge of $R_q\subset J_q(E)$ or from the knowledge of an epimorphism $\Phi:J_q(E)\rightarrow F_0$ if we set $R_q=ker(\Phi)$. It follows that the hybrid sequence projects onto the Janet sequence and that the kernel of this projection is the Spencer sequence. Also, chasing in the diagram, we may finally define the Janet bundles, up to an isomorphism, by the formula:\\
\[    F_r={\wedge}^rT^*\otimes J_q(E)/({\wedge}^rT^*\otimes R_q + \delta ({\wedge}^{r-1}T^*\otimes S_{q+1}T^*\otimes E))   \]
\noindent
that will be crucially used in the next section dealing with the deformation theory of Lie equations when $E=T$, $R_q\subset J_q(T)$ is a transitive  involutive system of infinitesimal Lie equations of order $q$ and the corresponding operator $\cal{D}$ is a Lie operator.  \\

\noindent
{\bf DEFINITION  4.4}: The Janet sequence is said to be {\it locally exact at} $F_r$ if any local section of $F_r$ killed by ${\cal{D}}_{r+1}$ is the image by ${\cal{D}}_r$ of a local section of $F_{r-1}$ over a convenient open subset. It is called {\it locally exact} if it is locally exact at each $F_r$ for  $0 \leq r \leq n $. The Poincar\' {e} sequence is locally exact but counterexamples may exist ([36], p 202).  \\
\vspace*{5mm}  \\
\noindent
{\bf  5  DEFORMATION THEORY OF LIE EQUATIONS}:  \\

Having in mind the application of computer algebra to the local theory of Lie pseudogroups, we want first of all to insist on two points which have never been emphasized up to our knowledge.  \\

In order to motivate the {\it first point}, we sketch it on an example.  \\

\noindent
{\bf EXAMPLE  5.1}:  Let $\alpha={\alpha}_i(x)dx^i$ be a $1$-form as in Example 1.1 or 1.4. If we look for infinitesimal transformations preserving $\alpha$, we have to cancel the Lie derivative as follows:  \\
\[  ({\cal{L}}(\xi)\alpha)_i \equiv {\alpha}_r(x){\partial}_i{\xi}^r+{\xi}^r{\partial}_r{\alpha}_i(x)=0   \]
As a byproduct, we have a well defined {\it Lie operator} ${\cal{D}}:T \rightarrow T^*:\xi \rightarrow {\cal{D}}\xi={\cal{L}}(\xi)\alpha$ such that ${\cal{D}}\xi=0, {\cal{D}}\eta=0 \Rightarrow {\cal{D}}[\xi,\eta] =0$ because of the well known property of the Lie derivative $[{\cal{L}}(\xi),{\cal{L}}(\eta)]={\cal{L}}(\xi)\circ {\cal{L}}(\eta)-{\cal{L}}(\eta)\circ {\cal{L}}(\xi)={\cal{L}}([\xi,\eta])$ and such a property can be extended to tensors or even any geometric object. Accordingly, it is usual to say that, if we have two solutions of the system, then their bracket is again a solution. However such a result is coming from mathematics and {\it cannot } be recognized by means of computer algebra, contrary to what is sometimes claimed. Surprisingly, the underlying reason has to do with formal integrability. Indeed, if we study the first derivative of the bracket $[\xi, \eta]$, it involves in fact the {\it second derivatives} of $\xi$  and $\eta$ and sometimes things may change a lot. For example, if $d\alpha \neq 0$, as ${\cal{L}}(\xi)d\alpha=di(\xi)d\alpha=d{\cal{L}}(\xi)\alpha = 0 $, then the first order equations brought by ${\cal{L}}(\xi)d\alpha=0$ may not be linear combinations of the first order equations brought by ${\cal{L}}(\xi)\alpha =0$. In particular, if $n=2$ and $\alpha = x^2dx^1$, then $d\alpha = - dx^1\wedge dx^2$ and the new first order equation ${\partial}_1{\xi}^1+{\partial}_2{\xi}^2=0$, which is {\it automatically } satisfied by any solution of the system $R_1\subset J_1(T)$ defined by ${\cal{L}}(\xi)\alpha=0$, is not a linear combination of the equations defining $R_1$, that is $R_1$ is not involutive as it is not even formally integrable because $R^{(1)}_1={\pi}^2_1(R_2)\subset R_1$ with a strict inclusion. {\it Hence it becomes a challenge to define a kind of bracket for sections of $R_1\subset J_1(T)$ and not for solutions as usual}, {\it that is to say independently of formal integrability}. This idea, which is a crucial one indeed as it will lead to the concept of Lie algebroid, is to replace the {\it classical Lie derivative} ${\cal{L}}(\xi)$ for any $\xi\in T$ by a {\it formal Lie derivative} $L({\xi}_1)$ for any ${\xi}_1\in J_1(T)$ over $\xi\in T$ in such a way that ${\cal{L}}(\xi)=L(j_1(\xi))$ and to compute the bracket $[L({\xi}_1),L({\eta}_1)]=L({\xi}_1)\circ L({\eta}_1)-L({\eta}_1)\circ L({\xi}_1)$ in the operator sense in order to be sure that the new bracket on $J_1(T)$ will satisfy the desired Jacobi identity. We obtain successively:  \\
\[ (L({\xi}_1)\alpha)_i\equiv {\alpha}_r(x){\xi}^r_i+{\xi}^r{\partial}_r{\alpha}_i(x)=0,\hspace{1cm}
   (L({\eta}_1)\alpha)_j\equiv {\alpha}_s(x){\eta}^s_j+{\eta}^s{\partial}_s{\alpha}_j(x)=0      \]
\[ L([{\xi}_1,{\eta}_1])\alpha\equiv {\alpha}_k(x)({\xi}^r{\partial}_r{\eta}^k_i+{\xi}^r_i{\eta}^k_r-{\eta}^s_i{\xi}^k_s-{\eta}^s{\partial}_s{\xi}^k_i)+([\xi,\eta])^k{\partial}_k{\alpha}_i(x)=0     \]
as a way to {\it define}:   \\
\[   ([{\xi}_1,{\eta}_1])^k=([\xi,\eta])^k={\xi}^r{\partial}_r{\eta}^k-{\eta}^s{\partial}_s{\xi}^k   \]
\[   ([{\xi}_1,{\eta}_1])^k_i = {\xi}^r{\partial}_r{\eta}^k_i+{\xi}^r_i{\eta}^k_r-{\eta}^s_i{\xi}^k_s-{\eta}^s{\partial}_s{\xi}^k_i   \]
{\it The induced property} $[R_1,R_1]\subset R_1$ {\it can therefore be checked linearly on sections and no longer on solutions}. In particular, we may exhibit the section $\{{\xi}^1=0,{\xi}^2=0,{\xi}^1_1=0,{\xi}^1_2=0,{\xi}^2_1=0,{\xi}^2_2=1\}$ of $R_1$, even if ${\xi}^1_1+{\xi}^2_2\neq 0$. \\

Now, using the {\it algebraic bracket} $\{ j_{q+1}(\xi),j_{q+1}(\eta)\}=j_q([\xi,\eta]), \forall \xi,\eta\in T$, we may  obtain by bilinearity a {\it differential bracket} on $J_q(T)$ extending the bracket on $T$:\\
\[   [{\xi}_q,{\eta}_q]=\{{\xi}_{q+1},{\eta}_{q+1}\}+i(\xi)D{\eta}_{q+1}-i(\eta)D{\xi}_{q+1},\hspace{1cm} \forall {\xi}_q,{\eta}_q\in J_q(T) \]
which does not depend on the respective lifts ${\xi}_{q+1}$ and ${\eta}_{q+1}$ of ${\xi}_q$ and ${\eta}_q$ in $J_{q+1}(T)$. Applying $j_q$ to the Jacobi identity for the ordinary bracket, we obtain:  \\
\[  \{{\xi}_{q+1},\{{\eta}_{q+2},{\zeta}_{q+2}\}\}+\{{\eta}_{q+1},\{{\zeta}_{q+2},{\xi}_{q+2}\}\}+\{{\zeta}_{q+1},\{{\xi}_{q+2},{\eta}_{q+2}\}\}\equiv 0 \hspace{3mm}  \forall {\xi}_{q+2},{\eta}_{q+2}, {\zeta}_{q+2}\in J_{q+2}(T)   \]
As we shall see later on, this bracket on sections satisfies the Jacobi identity and the following definition is the only one that can be tested by means of computer algebra:\\

\noindent
{\bf DEFINITION  5.2}: We say that a vector subbundle $R_q\subset J_q(T)$ is a {\it system of infinitesimal Lie equations} or a {\it Lie algebroid} if $[R_q,R_q]\subset R_q$, that is to say $[{\xi}_q,{\eta}_q]\in R_q, \forall {\xi}_q,{\eta}_q\in R_q$. The kernel $R^0_q$ of the projection ${\pi}^q_0:R_q\rightarrow T$ is the {\it isotropy Lie algebra bundle} of ${\cal{R}}^0_q={id}^{-1}({\cal{R}}_q)$ and $[R^0_q,R^0_q]\subset R^0_q$ does not contain derivatives, being thus defined fiber by fiber.   \\

Of course, another difficulty to overcome in this new setting, is that we have no longer an identity like $d{\cal{L}}(\xi)\alpha-{\cal{L}}(\xi)d\alpha=0$ but it is easy to check in local coordinates that:  \\
\[      (dL({\xi}_1)\alpha-L({\xi}_1)d\alpha)_{ij}={\alpha}_r(x)({\partial}_i{\xi}^r_j-{\partial}_j{\xi}^r_i)+({\partial}_i{\xi}^r-{\xi}^r_i){\partial}_r{\alpha}_j(x)-({\partial}_j{\xi}^r-{\xi}^r_j){\partial}_r{\alpha}_i (x)   \]
and the Spencer operator allows to factorize the formula if we notice that:ÊÊ\\
\[     ({\partial}_i{\xi}^r_j - {\partial}_j{\xi}^r_i)=({\partial}_i{\xi}^r_j-{\xi}^r_{ij})-({\partial}_j{\xi}^r_i-{\xi}^r_{ij})  \]
We finally obtain:  \\
\[      i({\zeta}_{(1)})i({\zeta}_{(2)})(dL({\xi}_1)\alpha-L({\xi}_1)d\alpha)=i({\zeta}_{(2)})L(i({\zeta}_{(1)})D{\xi}_2)\alpha-i({\zeta}_1)L(i({\zeta}_{(2)}D{\xi}_2)\alpha    \]
and more generally:   \\

\noindent
{\bf LEMMA 5.3}: When $\alpha \in {\wedge}^{r-1}T^*$ we have the formula:   \\
\[ i({\zeta}_{(1)}) ... i({\zeta}_{(r)})(dL({\xi}_1)\alpha - L({\xi}_1)d\alpha)=\sum^{r}_{s=1}(-1)^{s+1}i({\zeta}_{(1)}) ... i({\hat{\zeta}}_{(s)}) ... i({\zeta}_{(r)})L(i({\zeta}_{(s)}D{\xi}_2)\alpha   \]
which does not depend on the lift ${\xi}_2\in J_2(T)$ of ${\xi}_1\in J_1(T)$.  \\

In order to understand the {\it second point}, we have to revisit the work of Vessiot. Indeed, if we have a geometric object, that is a section $\omega$ of a natural bundle ${\cal{F}}$ of order $q$, then we may consider the system ${\cal{R}}_q=\{f_q\in {\Pi}_q\mid f_q^{-1}(\omega)=\omega\}$ of finite Lie equations and the corresponding linearized system $R_q=\{{\xi}_q\in J_q(T) \mid L({\xi}_q)\omega=0\}$ of infinitesimal Lie equations, both with the particular way to write them out, namely the Lie form and the Medolaghi form as it becomes clear from Example 1.1 to 1.4. As a byproduct, when constructing the Janet sequence, we can write $F_0=J_q(T)/R_q$ but we can also use the isomorphic definition $F_0={\omega}^{-1}(V({\cal{F}}))$ {\it depending on whether we want to pay attention to the system or to the object}. The main idea of deformation theory will be to begin with the second point of view and finish with the first. Starting with a system $R_q\subset J_q(T)$, we shall suppose that $R_q$ is transitive with a short exact sequence $0 \rightarrow R^0_q \rightarrow R_q \rightarrow T \rightarrow 0$ and, whatever is the definition of $F_0$, introduce an epimorphism $ \Phi:J_q(T) \rightarrow F_0$ while considering the following commutative and exact diagram:   \\
\[   \begin{array}{rcccccl}
   & 0 && 0 && 0 &    \\
   & \downarrow & & \downarrow & & \downarrow  &  \\
   0 \longrightarrow & R^0_q & \longrightarrow & J^0_q(T) & \longrightarrow & F_0 & \longrightarrow 0  \\
   & \downarrow &  &  \downarrow &  & \parallel  &   \\
   0 \longrightarrow & R_q   & \longrightarrow & J_q(T)  & \stackrel{\Phi}{\longrightarrow} & F_0 &\longrightarrow 0   \\
    & \downarrow &  &  \downarrow &  & \downarrow  &   \\
0\longrightarrow & T & =  & T & \longrightarrow & 0 &    \\
   &   \downarrow  &  &  \downarrow &  &  &   \\
   &         0   &&    0   & & &
   \end{array}    \]

The next definition will also be crucial for our purpose and generalizes the standard definition:    \\
\[  {\cal{L}}(\xi)\omega=\frac{d}{dt}j_q(exp\hspace{2mm} t\xi)^{-1}(\omega){\mid}_{t=0}.  \]

\noindent
{\bf DEFINITION  5.4}: We say that a vector bundle $F$ is {\it associated} with $R_q$ if there exists a first order differential operator $L({\xi}_q):F \rightarrow  F$ called {\it formal Lie derivative} and such that:  \\
1)  $L({\xi}_q+{\eta}_q)=L({\xi}_q)+L({\eta}_q)     \hspace{2cm} \forall {\xi}_q,{\eta}_q\in R_q$.  \\
2)  $L(f{\xi}_q)=fL({\xi}_q)  \hspace{2cm}   \forall {\xi}_q\in R_q, \forall f\in C^{\infty}(X)$.   \\
3)  $[L({\xi}_q),L({\eta}_q)]=L({\xi}_q)\circ L({\eta}_q)-L({\eta}_q)\circ L({\xi}_q)=L([{\xi}_q,{\eta}_q])  \hspace{2cm}  \forall {\xi}_q,{\eta}_q\in R_q$.   \\
4)  $ L({\xi}_q)(f\eta)=fL({\xi}_q)\eta + ({\xi}.f)\eta \hspace{15mm} \forall {\xi}_q\in R_q,\forall f\in C^{\infty}(X), \forall \eta \in F$ where ${\xi}.f=i(\xi)df$. \\

As a byproduct, if $E$ and $F$ are associated with $R_q$, we may set on $E\otimes F$:Ê\\
\[    L({\xi}_q)(\eta \otimes \zeta)=L({\xi}_q)\eta\otimes \zeta+\eta\otimes L({\xi}_q)\zeta \hspace{2cm} \forall {\xi}_q\in R_q,\forall \eta \in E, \forall \zeta \in F   \] 

\noindent
{\bf REMARK  5.5}: If $\Theta \subset T$ denotes the solutions of $R_q$, then ${\cal{L}}(\xi)=L(j_q(\xi))$ is simply called the {\it classical Lie derivative} but cannot be used in actual practice as we already said because $\Theta$ may be infinite dimensional as in Example 1.3 and 1.4. We obtain at once:  \\
1)  ${\cal{L}}(\xi+\eta)={\cal{L}}(\xi)+{\cal{L}}(\eta)  \hspace{2cm} \forall \xi,\eta \in \Theta$.  \\
3)  $[{\cal{L}}(\xi),{\cal{L}}(\eta)]={\cal{L}}([\xi,\eta])  \hspace{2cm} \forall \xi,\eta \in \Theta$.  \\
4)  ${\cal{L}}(\xi)(f\eta)=f{\cal{L}}(\xi)\eta+({\xi}.f) \eta \hspace{2cm} \forall \xi\in \Theta, \forall f\in C^{\infty}(X), \forall \eta\in F$.  \\

The extension to tensor products is well known.  \\

The following technical proposition and its corollary will be of constant use later on:  \\

\noindent
{\bf PROPOSITION 5.6}: We have:\\
\[  i(\zeta)D\{{\xi}_{q+1},{\eta}_{q+1}\}=\{i(\zeta)D{\xi}_{q+1},{\eta}_q\}+\{{\xi}_q,i(\zeta)D{\eta}_{q+1}\}  \]

\noindent
{\it Proof}: We have:   \\
 \[  (\{{\xi}_{q+1},{\eta}_{q+1}\})^k_{\nu}=\sum_{\lambda+\mu=\nu}({\xi}^r_{\lambda}{\eta}^k_{\mu+1_r}-{\eta}^s_{\lambda}{\xi}^k_{\mu+1_s}) \]
Now, caring only about ${\xi}_{q+1}$, we get: \\
\[  {\partial}_i(\{{\xi}_{q+1},{\eta}_{q+1}\})^k_{\nu}-(\{{\xi}_{q+1},{\eta}_{q+1}\})^k_{\nu+1_i}=\sum_{\lambda+\mu=\nu}({\partial}_i{\xi}^r_{\lambda}-{\xi}^r_{\lambda+1_i}){\eta}^k_{\mu+1_r}-({\partial}_i{\xi}^k_{\mu+1_s}-{\xi}^k_{\mu+1_s+1_i}){\eta}^s_{\lambda}+ ...  \]
and the Proposition follows by bilinearity.  \\
\hspace*{12cm}  Q.E.D.  \\

The proof of the following proposition is similar and left to the reader as an exercise:  \\

\noindent
{\bf PROPOSITION  5.7}: We have the formula:  \\
\[  i(\zeta)D[{\xi}_{q+1},{\eta}_{q+1}]=[i(\zeta)D{\xi}_{q+1},{\eta}_q]+[{\xi}_q,i(\zeta)D{\eta}_{q+1}]+i(L({\eta}_1)\zeta)D{\xi}_{q+1}-i(L({\xi}_1)\zeta)D{\eta}_{q+1}   \]

\noindent
{\bf COROLLARY  5.8}:  If $R_q\subset J_q(T)$ is such that $[R_q,R_q]\subset R_q$, then $R_{q+1}\subset J_{q+1}(T)$ satisfies $[R_{q+1},R_{q+1}]\subset R_{q+1}$  even if $R_q$ is not formally integrable.  \\

\noindent
{\bf EXAMPLE  5.9}: $T$ and $T^*$ both with any tensor bundle are associated with $J_1(T)$. The case of $T^*$ has been treated at the beginning of this section while for $T$ we may define $L({\xi}_1)\eta=[\xi,\eta]+i(\eta)D{\xi}_1=\{{\xi}_1,j_1(\eta)\}$. We have indeed ${\xi}^r{\partial}_r{\eta}^k-{\eta}^s{\partial}_s{\xi}^k+{\eta}^s({\partial}_s{\xi}^k-{\xi}^k_s)=-{\eta}^s{\xi}^k_s+{\xi}^r{\partial}_r{\eta}^k$ and the four properties of the formal Lie derivative can be checked directly as we did for $T^*$. Of course, we find back ${\cal{L}}(\xi)\eta=[\xi,\eta], \forall \xi,\eta \in T$.  \\

More generally, we have in a coherent way:   \\

\noindent
{\bf PROPOSITION  5.10}: $J_q(T)$ is associated with $J_{q+1}(T)$ if we define:  \\
\[     L({\xi}_{q+1}){\eta}_q=\{{\xi}_{q+1},{\eta}_{q+1}\}+i(\xi)D{\eta}_{q+1}=[{\xi}_q,{\eta}_q]+
i(\eta)D{\xi}_{q+1}        \]
and thus $R_q$ is associated with $R_{q+1}$.  \

\noindent
{\it Proof}: It is easy to check the properties 1, 2, 4 and it only remains to prove property 3 as follows.\\
\[  \begin{array}{rcl}
[L({\xi}_{q+1}),L({\eta}_{q+1})]{\zeta}_q & = & L({\xi}_{q+1})(\{{\eta}_{q+1},{\zeta}_{q+1}\}+i(\eta)D{\zeta}_{q+1})-L({\eta}_{q+1})(\{{\xi}_{q+1},{\zeta}_{q+1}\}+i(\xi)D{\zeta}_{q+1})   \\
       &  =  & \{{\xi}_{q+1},\{{\eta}_{q+2},{\zeta}_{q+2}\}\}-\{{\eta}_{q+1},\{{\xi}_{q+2},{\zeta}_{q+2}\}\}  \\
   &   & +\{{\xi}_{q+1},i(\eta)D{\zeta}_{q+2}\}-\{{\eta}_{q+1},i(\xi)D{\zeta}_{q+2}\}   \\
   &   & +i(\xi)D\{{\eta}_{q+2},{\zeta}_{q+2}\}-i(\eta)D\{{\xi}_{q+2},{\zeta}_{q+2}\}   \\
   &   & +i(\xi)D(i(\eta)D{\zeta}_{q+2})-i(\eta)D(i(\xi)D{\zeta}_{q+2})    \\
   &= & \{\{{\xi}_{q+2},{\eta}_{q+2}\},{\zeta}_{q+1}\}+\{i(\xi)D{\eta}_{q+2},{\zeta}_{q+1}\}-\{i(\eta)D{\xi}_{q+2},{\zeta}_{q+1}\}      \\
    &   &+i([\xi,\eta])D{\zeta}_{q+1}   \\
   &= & \{[{\xi}_{q+1},{\eta}_{q+1}],{\zeta}_{q+1}\}+i([\xi,\eta])D{\zeta}_{q+1}
   \end{array}   \]
by using successively the Jacobi identity for the algebraic bracket and the last proposition.Ê\\
\hspace*{12cm}   Q.E.D.   \\

\noindent
{\bf COROLLARY  5.11}: The differential bracket satisfies the Jacobi identity :   \\
\[  [{\xi}_q,[{\eta}_q,{\zeta}_q]]+[{\eta}_q,[{\zeta}_q,{\xi}_q]]+[{\zeta}_q,[{\xi}_q,{\eta}_q]]\equiv 0  \hspace{2cm}   \forall {\xi}_q,{\eta}_q,{\zeta}_q \in J_q(T)   \]

\noindent
{\bf PROPOSITION  5.12}: We have the formula:   \\
\[ i(\zeta)(DL({\xi}_{q+2}){\eta}_{q+1}-L({\xi}_{q+1})D{\eta}_{q+1})=L(i(\zeta)D{\xi}_{q+2}){\eta}_q \]

\noindent
{\it Proof}: Using Proposition 5.6, we have:  \\
\[  \begin{array}{rcl}
   i(\zeta)DL({\xi}_{q+2}){\eta}_{q+1}&=& i(\zeta)D\{{\xi}_{q+2},{\eta}_{q+2}\}+i(\zeta)Di(\xi)D{\eta}_{q+2} \\           & =  & \{i(\zeta)D{\xi}_{q+2},{\eta}_{q+1}\}+\{{\xi}_{q+1},i(\zeta)D{\eta}_{q+2}\}+i(\zeta)Di(\xi)D{\eta}_{q+2} 
  \end{array}        \]
and we must substract:   \\
\[\begin{array}{rcl}
i(\zeta)L({\xi}_{q+1})D{\eta}_{q+1} & = & L({\xi}_{q+1})(i(\zeta)D{\eta}_{q+1})-i(L({\xi}_1)\zeta)D{\eta}_{q+1}    \\
        & =   & \{{\xi}_{q+1},i(\zeta)D{\eta}_{q+2}\}+i(\xi)Di(\zeta)D{\eta}_{q+2}-i(L({\xi}_1)\zeta)D{\eta}_{q+1}
\end{array}     \]
in order to obtain for the difference:  \\
\[ L(i(\zeta)D{\xi}_{q+2}){\eta}_q-i(i(\zeta)D{\xi}_1)D{\eta}_{q+1}+i(L({\xi}_1)\zeta)D{\eta}_{q+1}+i(\zeta)Di(\xi)D{\eta}_{q+2}-i(\xi)Di(\zeta)D{\eta}_{q+2}    \]
Finally, the last four terms vanish because $L({\xi}_1)\zeta-i(\zeta)D{\xi}_1=[\xi,\zeta] $ and:  \\
\[   i(\zeta)Di(\xi)D{\eta}_{q+2}-i(\xi)Di(\zeta)D{\eta}_{q+2}=- i([\xi,\zeta])D{\eta}_{q+1}.   \]
\hspace*{12cm}        Q.E.D.   \\

Combining this proposition and Lemma 5.3, we obtain:  \\

\noindent
{\bf PROPOSITION  5.13}: When $A^{r-1}_{q+1}\in {\wedge}^{r-1}T^*\otimes J_{q+1}(T)$, we have the formula:  \\
\[ i({\zeta}_{(1)}) ... i({\zeta}_{(r)})(DL({\xi}_{q+2})-L({\xi}_{q+1})D)A^{r-1}_{q+1}=
   \sum_{s=1}^{r}(-1)^{s+1}i({\zeta}_{(1)}) ... i({\hat{\zeta}}_{(s)}) ... i({\zeta}_{(r)})L(i({\zeta}_{(s)})D{\xi}_{q+2})A^{r-1}_q   \]
where $A^{r-1}_q\in {\wedge}^{r-1}T^*\otimes J_q(T)$ is the projection of $A^{r-1}_{q+1}$.   \\
   
\noindent
{\it Proof}: With $\alpha\in {\wedge}^{r-1}T^*$ and ${\eta}_{q+1}\in J_{q+1}(T)$, we obtain successively:\\
\[  \begin{array}{rcl}
DL({\xi}_{q+2})(\alpha\otimes {\eta}_{q+1})& = & D(L({\xi}_1)\alpha\otimes {\eta}_{q+1}+\alpha\otimes L({\xi}_{q+2}){\eta}_{q+1})  \\
     &  =  & dL({\xi}_1)\alpha\otimes {\eta}_q+(-1)^{r-1}(L({\xi}_1)\alpha)\wedge D{\eta}_{q+1}   \\
     &  & +d\alpha\otimes L({\xi}_{q+1}){\eta}_q+(-1)^{r-1}\alpha \wedge DL({\xi}_{q+2}){\eta}_{q+1}
     \end{array}   \]
\[   \begin{array}{rcl}
L({\xi}_{q+1})D(\alpha\otimes {\eta}_{q+1}) & =& L({\xi}_{q+1})(d\alpha\otimes {\eta}_q+(-1)^{r-1}\alpha\wedge D{\eta}_{q+1})    \\
 & = & L({\xi}_1)d\alpha\otimes {\eta}_q+d\alpha\otimes L({\xi}_{q+1}){\eta}_q   \\
   &  &+(-1)^{r-1}L({\xi}_1)\alpha\wedge D{\eta}_{q+1}+(-1)^{r-1}\alpha\wedge L({\xi}_{q+1})D{\eta}_{q+1}
 \end{array}   \]
 and obtain y substraction: \\
 \[ (DL({\xi}_{q+2})-L({\xi}_{q+1})D)(\alpha\otimes {\eta}_{q+1})=(dL({\xi}_1)-L({\xi}_1)d)\alpha\otimes {\eta}_q+(-1)^{r-1}\alpha\wedge (DL({\xi}_{q+2})-L({\xi}_{q+1})D){\eta}_{q+1}  \]
 and the proposition follows by skewlinearity.\\
 \hspace*{12cm}    Q.E.D.   \\

\noindent
{\bf PROPOSITION  5.14}: We have the formula:   \\
\[   L({\xi}_q)\{{\eta}_q,{\zeta}_q\}=\{L({\xi}_{q+1}){\eta}_q,{\zeta}_q\}+\{{\eta}_q,L({\xi}_{q+1}){\zeta}_q\}  \]
which does not depend on the lift ${\xi}_{q+1}\in J_{q+1}(T)$ of ${\xi}_q\in J_q(T)$.\\

\noindent
{\it Proof}: Using the Jacobi identity for the algebraic bracket and proposition 5.6, we obtain:      \\
\[  \begin{array}{rcl}
L({\xi}_q)\{{\eta}_q,{\zeta}_q\} &   = &\{{\xi}_q, \{{\eta}_{q+1},{\zeta}_{q+1}\}+i(\xi)D\{{\eta}_{q+1},
{\zeta}_{q+1}\} \\
    & = & \{\{{\xi}_{q+1},{\eta}_{q+1}\},\zeta_q\}+\{{\eta}_q,\{{\xi}_{q+1},{\zeta}_{q+1}\} \}  \\
     &   &+\{i(\xi)D{\eta}_{q+1},{\zeta}_q\}+\{ {\eta}_q, i(\xi)D{\zeta}_{q+1}  \}    \\
     & =  & \{L({\xi}_{q+1}){\eta}_q, {\zeta}_q\}+\{{\eta}_q,L({\xi}_{q+1}){\zeta}_q\}
     \end{array}    \]
\hspace{12cm}       Q.E.D.  \\

Finally, using Proposition 5.12, we obtain at once:  \\

\noindent
{\bf COROLLARY  5.15}: We have the formula:   \\
\[  \begin{array}{rcl}
L({\xi}_{q+1})[{\eta}_q,{\zeta}_q] & =  & [L({\xi}_{q+1}){\eta}_q,{\zeta}_q]+[{\eta}_q,L({\xi}_{q+1}){\zeta}_q]\\
     &     & +L(i(\zeta)D{\xi}_{q+2}){\eta}_q-L(i(\eta) D{\xi}_{q+2}){\zeta}_q   
     \end{array}   \]
which does not depend on the lift ${\xi}_{q+2}\in J_{q+2}(T)$ of ${\xi}_{q+1}\in J_{q+1}(T)$.  \\

Before going ahead, let us stop for a moment and wonder how we could proceed for generalizing the deformation theory of Lie algebras by using the Vessiot structure equations even though we know that the structure constants have nothing to do in general with any Lie algebra. Of course we could start similarly from the Jacobi relations but, if we do want to exhibit a kind of cohomology, we should be able to define a trivial deformation, that is the analogue of a change of basis of the underlying vector space $V$ of the Lie algebra ${\cal{G}}$ in such a natural way that it could induce a change of the structure constants which is surely not of a tensorial nature anymore.  \\
  The following "trick", already known to Vessiot in 1903 ( [51], p 445), is still ignored today. For this, assuming that the natural bundle ${\cal{F}}$ is known, let us consider two sections $\omega$ and $\bar{\omega}$ giving rise respectively to the systems $R_q$ and ${\bar{R}}_q$ of infinitesimal Lie equations:  \\
 \[        R_q            \hspace{4cm}  {\Omega}^{\tau}\equiv - L^{\tau\mu}_k(\omega (x)){\xi}^k_{\mu}+
 {\xi}^r{\partial}_r{\omega}^{\tau}(x)=0      \]
\[       {\bar{R}}_q   \hspace{4cm}  {\bar{\Omega}}^{\tau}\equiv - L^{\tau\mu}_k(\bar{\omega} (x)){\xi}^k_{\mu}+{\xi}^r{\partial}_r{\bar{\omega}}^{\tau}(x)=0      \]
and define the following equivalence relation:   \\

\noindent
{\bf DEFINITION  5.16}:   $\bar{\omega} \sim  \omega  \Leftrightarrow  {\bar{R}}_q=R_q  $ \\

The study of such an equivalence relation is not evident at all and we improve earlier presentations (compare to [36], p 336). First of all, having in mind what we did for Example 1.4, we shall use a solved form of the system obtained by choosing principal jets or, equivalently, choosing a square submatrix $M=(M(u))$ of rank $dim(F_0)$ in the matrix $L=(L(u))$ defining $R^0_q$ with $dim(J^0_q(T))$ columns which describe the infinitesimal generators of prolongations of changes of coordinates on $X$ acting on the fibers of ${\cal{F}}$ and $dim(F_0)=m$ rows. Of course, a major problem will be to obtain intrinsic results not depending on this choice. The columns of $L$ are thus made by vector fields $L^{\mu}_k={L}^{\tau\mu}_k(u) \frac{\partial}{\partial u^{\tau}}$ that we can therefore separate into two parts, namely the vectors ${\L}_{\sigma}={M}^{\tau}_{\sigma}(u)\frac{\partial}{\partial u^{\tau}}$ for $\sigma=1,...,dim(F_0)$ and the vectors $L_{m+r}={\cal{E}}^{\sigma}_{m+r}(u)L_{\sigma}$ for $r=1,...,dim(R^0_q)$ obtained by introducing the {\it stationary functions} ${\cal{E}}(u)$, also called {\it Grassmann determinants}, while describing the matrix $M^{-1}L=(id_{F_0},{\cal{E}}(u))$. We are therefore led to look for transformations 
$\bar{u}=g(u)$ of the fibers of ${\cal{F}}$ such that:    \\
\[  (M^{-1})^{\sigma}_{\tau}(\bar{u})du^{\tau}=(M^{-1})^{\sigma}_{\tau}(u)du^{\tau}, 
\hspace{1cm} {\cal{E}}^{\sigma}_{m+r}(\bar{u})={\cal{E}}^{\sigma}_{\tau}(u)      \]
In order to study such a system and to prove that it is defining a Lie pseudogroup of transformations, let us notice that the first conditions are equivalent to saying that the transformations $\bar{u}=g(u)$ preserve the vector fields $L_{\sigma}$ and also the vector fields $L_{m+r}$ according to the second conditions. It follows that the transformations $\bar{u}=g(u)$ preserves the vector fields $L^{\mu}_k$, a property thus not depending on the choice of the principal jets. In addition, we have:   \\

\noindent
{\bf PROPOSITION  5.17}: The Lie pseudogroup of transformations of the fibers of ${\cal{F}}$ that we have exhibited is in fact a Lie group of transformations, namely the {\it reciprocal} of the lie group of transformations describing the natural structure of ${\cal{F}}$.  \\

\noindent
{\it Proof}: The defining system is finite type with a zero first order symbol. If $W={W}^{\tau}(u)\frac{\partial}{\partial u^{\tau}}$ is an infinitesimal transformation, we obtain therefore the Lie operator 
$[W,L^{\mu}_k]=0, \forall 1\leq {\mid}\mu{\mid}\leq q, \forall k=1,...,m$. Indeed, if $W_1$ and $W_2$ are two solutions, then $[W_1,W_2]$ is also a solution because of the Jacobi identity for the bracket and there are at most $dim(F_0)$ linearly independent such vector fields denoted by $W_{\alpha}$. It follows that $[W_1,W_2]={\rho}^{\alpha}_{12}(u)W_{\alpha}$ and we deduce from the Jacobi identity again that $L^{\mu}_k\cdot {\rho}^{\alpha}_{12}(u)=0 \Rightarrow {\rho}^{\alpha}_{12}(u)={\rho}^{\alpha}_{12}=cst$ because $rk(L^{\mu}_k)=dim(F_0)$. Accordingly, the $W_{\alpha}$ are the infinitesimal generators of a Lie group of transformations of the fibers of ${\cal{F}}$ and the effective action does not depend on the coordinate system. \\
\hspace*{12cm}     Q.E.D.  \\

\noindent
{\bf DEFINITION  5.18}: These finite transformations will be called {\it label transformations} and will be noted $\bar{u}=g(u,a)$ where the number of parameters $a$ is $\leq dim(F_0)$.  \\

If $R_q$ is formally integrable/involutive, then ${\bar{R}}_q=R_q$ is also formally integrable/involutive and thus $I(j_1(\omega))=c(\omega)\Leftrightarrow I(j_1(\bar{\omega}))=\bar{c}(\bar{\omega})$ with eventually different structure constants.  \\

\noindent
{\bf COROLLARY  5.19}: Any finite label transformation $\bar{u}=g(u,a)$ induces a finite transformation 
$\bar{c}=h(c,a)$ of the structure constants which is not effective in general and we may set $\bar{\omega}\sim \omega \Rightarrow \bar{c}\sim c$.   \\

It now remains to exhibit a deformation cohomology coherent with the above results.\\

\noindent
{\bf DEFINITION  5.20}: When $F$ is a vector bundle associated with $R_q$, we may define $\Upsilon=\Upsilon(F)=\{\eta\in F{\mid} L({\xi}_q)\eta=0, \forall {\xi}_q\in R_q\}$ and the sub-vector bundle $E=\{\eta\in F{\mid}L({\xi}^0_q)\eta=0, \forall {\xi}^0_q\in R^0_q\}\subseteq F$ in such a way that $\Upsilon \subset E\subseteq F$.\\

In order to look for $\Upsilon$ in general, we shall decompose this study into two parts, exactly as we did in section 3.7, by using a splitting of the short exact sequence $0\rightarrow R^0_q \rightarrow R_q \stackrel{{\pi}^q_0}{\rightarrow} T \rightarrow 0$ called $R_q$-{\it connection}, namely a map ${\chi}_q:T\rightarrow R_q$ such that ${\pi}^q_0\circ {\chi}_q=id_T$, in order to have $R_q\simeq R^0_q\oplus {\chi}_q(T)$. Such a procedure does not depend on the choice of ${\chi}_q$ because, if ${\bar{\chi}}_q$ is another $R_q$-connection, then $({\bar{\chi}}_q-{\chi}_q)(T)\in R^0_q$. An $R_q$-connection may also be considered as a section ${\chi}_q\in T^*\otimes R_q$ over $id_T\in T^*\otimes T$ and ${\bar{\chi}}_q-{\chi}_q\in T^*\otimes R^0_q$ in this case. It follows that we have equivalently $\Upsilon=\{\eta\in E{\mid}L({\chi}_q(\xi))\eta=0, \forall \xi\in T\}$ and we may define a first order operator $\nabla :E\rightarrow T^*\otimes E$ with zero symbol, called {\it covariant derivative}, by the formula $(\nabla\cdot \eta)(\xi)={\nabla}_{\xi}\eta=L({\chi}_q(\xi))\eta$. We may now extend $\nabla$ to a first order operator $\nabla= {\wedge}^rT^*\otimes E\rightarrow {\wedge}^{r+1}T^*\otimes E$ by the formula:  \\
\[   \nabla(\alpha\otimes \eta)=d\alpha\otimes \eta+ (-1)^r\alpha \wedge \nabla \eta, \hspace{1cm}\forall \alpha \in {\wedge}^rT^*, \forall \eta \in E   \]

\noindent
{\bf LEMMA  5.21}: With ${\nabla}^2=\nabla\circ \nabla$, we have:   \\
\[    ({\nabla}^2\eta)(\xi,{\bar{\xi}})=L([{\chi}_q(\xi),{\chi}_q({\bar{\xi}})]-{\chi}_q([\xi,{\bar{\xi}}]))\eta=0 \Rightarrow {\nabla}^2=0 \]

\noindent
{\it Proof}: We have:   \\
\[  \begin{array}{rcl}
{\nabla}^2(\alpha\otimes \eta)& = & \nabla (d\alpha\otimes \eta+(-1)^r\alpha\wedge\nabla \eta)  \\
  &  =  & d^2\alpha+(-1)^{r+1}d\alpha \wedge\nabla\eta+(-1)^rd\alpha\wedge\nabla\eta+(-1)^r\alpha\wedge{\nabla}^2\eta  \\
  &  =  & (-1)^r\alpha\wedge {\nabla}^2\eta
  \end{array}   \]
  Setting $\nabla\eta=dx^i{\nabla}_i\eta$, we get:   \\
  \[ {\nabla}^2\eta=\nabla(dx^i{\nabla}_i\eta)=dx^i\wedge dx^j{\nabla}_j{\nabla}_i\eta=-\frac{1}{2}dx^i\wedge dx^j({\nabla}_i{\nabla}_j-{\nabla}_j{\nabla}_i)\eta  \]
and we have just to use the fact that ${\nabla}_i=L({\chi}_q({\partial}_i))$ with 
$[{\partial}_i,{\partial}_j]=0$. Indeed, we have:   \\
\[  [{\chi}_q(\xi),{\chi}_q({\bar{\xi}})]-{\chi}_q([\xi,{\bar{\xi}}])=\{{\chi}_{q+1}(\xi),{\chi}_{q+1}({\bar{\xi}}\}+i(\xi)D{\chi}_{q+1}({\bar{\xi}})-i({\bar{\xi}})D{\chi}_{q+1}(\xi)-{\chi}_q([\xi,{\bar{\xi}}])  \]
The first term in the right member is linear in $\xi$ and ${\bar{\xi}}$ while the sum of the others 
becomes also linear in $\xi$ and ${\bar{\xi}}$ because:         \\
\[ {\xi}^i({\partial}_i({\chi}^k_{\mu,j}{\bar{\xi}}^j)-{\chi}^k_{\mu+1_i,j}{\bar{\xi}}^j)-{\bar{\xi}}^j({\partial}_j({\chi}^k_{\mu,i}{\xi}^i)-{\chi}^k_{\mu+1_j,i}{\xi}^i)-{\chi}^k_{\mu,r}({\xi}^i{\partial}_i{\bar{\xi}}^r-{\bar{\xi}}^j{\partial}_j{\xi}^r)\]
   \[=({\partial}_i{\chi}^k_{\mu,j}-{\partial}_j{\chi}^k_{\mu,i}+{\chi}^k_{\mu+1_j,i}-{\chi}^k_{\mu+1_i,j}){\xi}^i{\bar{\xi}}^j   \]
The proposition follows from the fact that ${\chi}_0=id_T$ and $E$ is $R^0_q$-invariant. \\
\hspace*{12cm}    Q.E.D.   \\

Hence, we obtain by linearity the $\nabla$-{\it sequence}: \\
\[   0 \longrightarrow \Upsilon \longrightarrow E \stackrel{\nabla}{\longrightarrow}T^*\otimes E \stackrel{\nabla}{\longrightarrow} {\wedge}^2T^*\otimes E \stackrel{\nabla}{\longrightarrow} ... \stackrel{\nabla}{\longrightarrow} {\wedge}^nT^*\otimes E \longrightarrow 0   \]
which does not depend on the choice of the connection and is made by first order involutive operators. It follows that $\Upsilon$ can be locally described by a linear combination with constant coefficients of certain sections of $E\subset F$ and we may therefore set $dim(\Upsilon)=dim(E)\leq dim(F)$. The use of computer algebra will essentially be to compute these dimensions by using linear algebra combined with homological algebra techniques.\\

We now provide a few definitions:   \\

\noindent
{\bf DEFINITION  5.22}: When $\Theta$ is given, we may define:  \\
{\it Centralizer}   \hspace{3cm}   $C(\Theta)=\{\eta\in T{\mid}[\xi,\eta]=0, \forall \xi\in \Theta \}$.  \\
{\it Center}    \hspace{37mm}         $Z(\Theta)=\{\eta\in \Theta {\mid}[\xi,\eta]=0, \forall \xi\in \Theta \}$.  \\
{\it Normalizer}\hspace{31mm}   $N(\Theta)=\{\eta\in T {\mid} [\xi,\eta]\subset \Theta, \forall \xi\in \Theta\}$. \\

It is essential to notice that these definitions are not very useful at all in actual practice when $\Theta$ is infinite dimensional.   \\

\noindent
{\bf PROPOSITION  5.23}: $C(\Theta)=\Upsilon (T) $.   \\

\noindent
{\it Proof}: If $R_1={\pi}^q_1(R_q)\subset J_1(T)$, it follows from Example 5.9 that $\Upsilon (T)=\{ \eta\in T{\mid}L({\xi}_1)\eta=0, \forall {\xi}_1\in R_1\}$, that is to say $\{ \Upsilon (T)=\{ \eta\in T{\mid} \{{\xi}_1,j_1(\eta)\}=0, \forall {\xi}_1\in R_1 \}$ if we choose $j_1(\eta)$ as a lift of $\eta$ in $J_1(T)$. In particular, if 
$\xi\in \Theta$ and thus $j_1(\xi)\in R_1$, we have $\{j_1(\xi),j_1(\eta)\}=[\xi,\eta]$ and 
thus $\Upsilon (T)\subseteq C(\Theta)$. \\
Now, $j_{q-1}([\xi,\eta])=\{j_q(\xi),j_q(\eta)\}$ and thus $C(\Theta)=\{\eta\in T {\mid} \{{\xi}_q,j_q(\eta)\}=0, \forall {\xi}_q\in R_q\}$, providing by projection $\{{\xi}_1,j_1(\eta)\}=0$, that is $C(\Theta)\subseteq \Upsilon (T)$ and thus $C(\Theta)=\Upsilon (T)$.   \\
\hspace*{12cm}       Q.E.D.   \\

It follows that $Z(\Theta)=\Theta\cap C(\Theta)\Rightarrow Z(\Theta)=\{ \eta\in T {\mid} {\cal{D}}\eta=0, L({\xi}_1)\eta=0, \forall {\xi}_1\in R_1\}$ and $Z(\Theta)$ is made by sections of $\Upsilon (T)$ killed by ${\cal{D}}$. The study of $N(\Theta)$ is much more delicate and we first need the next proposition where we notice the importance of involution or at least formal integrability.Ê \\

\noindent
{\bf PROPOSITION  5.24}: The Lie operator ${\cal{D}}:T\longrightarrow F_0$ induces a homomorphism of Lie algebras ${\cal{D}}:\Upsilon (T)\longrightarrow \Upsilon (F_0)$ where the bracket on 
$\Upsilon (T)$ is induced by the ordinary bracket on $T$ and the bracket on $\Upsilon (F_0)$ is induced by the differential bracket on $J_q(T)$.  \\

\noindent
{\it Proof}: We already know that $T$ inherits a structure of Lie algebra on sections from the ordinary bracket of vector fields and the situation is similar for $J_q(T)$ with the differential bracket. Now, if $L({\xi}_1)\eta=0 $ and $L({\xi}_1)\zeta=0$, it follows from Corollary 5.15 that $L({\xi}_1)[\eta,\zeta]=0$ and 
$[\Upsilon(T),\Upsilon(T)]\subset \Upsilon(T)$.  \\
Similarly, as $F_0=J_q(T)/R_q=J^0_q(T)/R^0_q$, if $L({\xi}_q){\eta}^0_q\equiv [{\xi}_q,{\eta}^0_q]\in R^0_q$ and $L({\xi}_q){\zeta}^0_q\equiv [{\xi}_q,{\zeta}^0_q]\in R^0_q$, then $L({\xi}_q)[{\eta}^0_q,{\zeta}^0_q]\in R^0_q$ according to the Jacobi identity for the bracket on $J^0_q(T)$ and the fact that $[R^0_q,R^0_q]\subset R^0_q$. Also, if $L({\xi}_{q+1}){\eta}_q\in R_q$ and $L({\xi}_{q+1}){\zeta}_q\in R_q$, then it is less evident to prove that $L({\xi}_{q+1})[{\eta}_q,{\zeta}_q]\in R_q$. For this, using Corollary 5.15, if we set 
$L({\xi}_{q+1}){\eta}_q={\theta}_q\in R_q$, it is sufficient to notice that $[{\theta}_q,{\zeta}_q]=L({\theta}_{q+1}){\zeta}_q-i(\zeta)D{\theta}_{q+1}\in R_q$ because both terms do belong to $R_q$. Thus, in any case, we obtain $[\Upsilon(F_0),\Upsilon(F_0)]\subset \Upsilon(F_0)$.  \\
Finally, we have $j_q([\eta,\zeta])=[j_q(\eta),j_q(\zeta)]$ and we may take $j_q(\eta)$ as a representative of ${\cal{D}}\eta$ in $J_q(T)$. We shall prove that, if $\eta\in\Upsilon(T)$, that is if $L({\xi}_1)\eta=0, \forall {\xi}_1\in R_1$, then $j_q(\eta)\in J_q(T)$ is such that $L({\xi}_{q+1})j_q(\eta)=0$. Introducing $R_2={\pi}^q_2(R_q)\subseteq {\rho}_1(R_1)$ and choosing any ${\xi}_2\in R_2$ over ${\xi}_1\in R_1$, we have ${\pi}^1_0(L({\xi}_2)j_1(\eta))=L({\xi}_1)\eta=0\Rightarrow L({\xi}_2)j_1(\eta)\in T^*\otimes T$. However, using Proposition 5.11 and the fact that $Dj_q(\eta)=0$, we obtain: \\
\[      i(\zeta)D(L({\xi}_2)j_1(\eta))=L(i(\zeta)D{\xi}_2)\eta=0     \]
because $R_2\subset {\rho}_1(R_1)\Rightarrow DR_2\subset T^*\otimes R_1$ and thus $L({\xi}_2)j_1(\eta)=0$ because there is a monomorphism (even an isomorphism) $0\rightarrow T^*\otimes T\stackrel{\delta}{\rightarrow}T^*\otimes T$. Supposing by induction that $L({\xi}_q)j_{q-1}(\eta)=0$, we should obtain in the same way:  \\
\[      i(\zeta)D(L({\xi}_{q+1})j_q(\eta))=L(i(\zeta)D{\xi}_{q+1})j_{q-1}(\eta)=0   \]
because $R_{q+1}\subseteq {\rho}_1(R_q)\Rightarrow DR_{q+1}\subset T^*\otimes R_q$ and thus $L({\xi}_{q+1})j_q(\eta)=0$ because there is a monomorphism $0\rightarrow S_{q+1}T^*\otimes T\stackrel{\delta}{\rightarrow}T^*\otimes S_qT^*\otimes T$ and the restriction of $D$ to a symbol is $-\delta$.Ê\\
\hspace*{12cm}   Q.E.D.   \\

As the reader will discover in the last computational section, {\it the study of the normalizer is much more delicate}.   \\

\noindent
{\bf DEFINITION  5.25}: The {\it normalizer} ${\tilde{\Gamma}}=N(\Gamma)$ of $\Gamma$ in $aut(X)$ is the biggest Lie pseudogroup in which $\Gamma$ is {\it normal}, that is (roughly) $N(\Gamma)={\tilde{\Gamma}}=\{ {\tilde{f}}\in aut(X){\mid} {\tilde{f}}\circ f\circ {\tilde{f}}^{-1}\in \Gamma, \forall f\in \Gamma \}$ and we write $\Gamma \lhd N(\Gamma)\subset aut(X)$.  \\ 

Of course, $N(\Theta)$ will play the part of a Lie algebra for $N(\Gamma)$ exactly like $\Theta$ did for $\Gamma$. However, we shall see that $N(\Gamma)$ may have many components different from the connected component of the identity, for example two in the case of the algebraic Lie pseudogroup of contact transformations where $N(\Gamma)/\Gamma$ is isomorphic to the permutation group of two objects, a result not evident at fist sight. Passing to the jets, we get $j_q({\tilde{f}}\circ f\circ {\tilde{f}}^{-1})^{-1}(\omega)=j_q({\tilde{f}})\circ j_q(f)^{-1}\circ j_q({\tilde{f}})^{-1}(\omega)=\omega \Leftrightarrow j_q(f)^{-1}(j_q({\tilde{f}})^{-1}(\omega))=j_q({\tilde{f}})^{-1}(\omega)$, that is to say $j_q(f)^{-1}({\bar{\omega}})={\bar{\omega}}$ if we set $j_q({\tilde{f}})^{-1}(\omega)={\bar{\omega}}$ and we find back the equivalence relation of Definition 5.16. It follows that 
${\tilde{\Gamma}}=\{{\tilde{f}}\in aut(X){\mid} j_q({\tilde{f}})^{-1}(\omega)=g(\omega,a),h(c,a)=c\}$ is defined by the system ${\tilde{{\cal{R}}}}_{q+1}=\{ {\tilde{f}}_{q+1}\in {\Pi}_{q+1}{\mid}{\tilde{f}}_{q+1}(R_q)=
R_q\}$ with linearization ${\tilde{R}}_{q+1}=\{ {\tilde{\xi}}_{q+1}{\mid} L({\tilde{\xi}}_{q+1}){\eta}_q\in R_q, \forall {\eta}_q\in R_q\}$, that is to say $\{{\tilde{\xi}}_{q+1},{\eta}_{q+1}\}+i({\tilde{\xi}})D{\eta}_{q+1}\in R_q \Leftrightarrow \{{\tilde{\xi}}_{q+1},{\eta}_{q+1}\}\in R_q$. Accordingly, the system of infinitesimal Lie equations defining ${\tilde{\Theta}}=N(\Theta)$ can be obtained by purely algebraic techniques from the system defining $\Theta$. In particular, we notice that ${\pi}^{q+1}_0:{\tilde{R}}_{q+1}\rightarrow T$ is an epimorphism because ${\pi}^{q+1}_0:R_{q+1}\rightarrow T$ is an epimorphisme by assumption and $R_{q+1}\subseteq {\tilde{R}}_{q+1}$. We obtain on the symbol level $\{{\tilde{g}}_{q+1},{\eta}_{q+1}\}\subset g_q$ and thus $\delta {\tilde{g}}_{q+1}\subset T^*\otimes g_q$ leading to ${\tilde{g}}_{q+1}\subseteq g_{q+1}={\rho}_1(g_q)$ and thus ${\tilde{g}}_{q+1}=g_{q+1}$ because $R_{q+1}\subseteq 
{\tilde{R}}_{q+1}\Rightarrow g_{q+1}\subseteq {\tilde{g}}_{q+1}$. Using arguments from $\delta$-cohomology, it can be proved that ${\tilde{R}}_{q+1}$ is involutive when $R_q $ is involutive ([36], p 351, 390). Another proof will be given in Corollary 5.52.  \\
With more details, using the result of Proposition 5.17, we get the following important local result:\\

\noindent
{\bf PROPOSITION  5.26}: ${\Upsilon}_0=\Upsilon(F_0)=\{ {\Omega}^{\tau}(x)=A^{\alpha}{W}^{\tau}_{\alpha}(\omega(x)) {\mid} A=cst  \}  $ \\

\noindent
{\it Proof}: Recalling that $F_0={\omega}^{-1}(V({\cal{F}}))$, we shall first study the natural bundle ${\cal{F}}_0=V({\cal{F}})$ of order $q$. Adopting local coordinates $(x,u,v)$, any infinitesimal change of source ${\bar{x}}=x+t\xi(x)+... $ can be lifted to ${\cal{F}}_0$ with ${\bar{u}}^{\tau}={u}^{\tau}+t{\xi}^k_{\mu}(x) L^{\tau\mu}_k(u)+ ... , {\bar{v}}^{\tau}=v^{\tau}+t\frac{\partial L^{\tau\mu}_k(u)}{\partial u^{\sigma}}v^{\sigma}+ ... $, according to the definition of a vertical bundle provided by Definition 3.3. The corresponding infinitesimal generators on ${\cal{F}}_0$ will be:  \\
\[    {\xi}^i(x)\frac{\partial}{\partial x^i} + {\xi}^k_{\mu}(x) (L^{\tau\mu}_k\frac{\partial}{\partial u^{\tau}}+\frac{\partial L^{\tau\mu}_k(u)}{\partial u^{\sigma}}v^{\sigma}\frac{\partial}{\partial v^{\tau}} ) \hspace{5mm} 1\leq {\mid}\mu{\mid}\leq q , \forall {\xi}_q\in J_q(T)  \]
It follows that a section $\epsilon:{\cal{F}}\rightarrow {\cal{F}}_0:(x,u)\rightarrow (x,u,v=\epsilon(x,u))$ will be equivariant, that is $v-\epsilon (x,u)=0\Rightarrow \bar{v}-\epsilon (\bar{x},\bar{u})=0$ if and only if $v=\epsilon (u)$ satisfies:   \\
\[      L^{\tau\mu}_k(u)\frac{\partial {\epsilon}^{\sigma}(u)}{\partial u^{\tau}} - \frac{\partial L^{\tau\mu}_k(u)}{\partial u^{\sigma}}{\epsilon}^{\sigma}(u)=0 \hspace{5mm},\hspace{5mm} 1\leq{\mid}\mu{\mid}\leq q \]
Hence, we have $[L^{\mu}_k,\epsilon]=0\Rightarrow {\epsilon}^{\tau}(u)=A^{\alpha}W^{\tau}_{\alpha}(u)$ with $A=cst$.\\
As another approach, working directly with $F_0$, we may consider the invariance of the section $u-\omega (x)=0, v-\Omega(x)=0$ and get:  \\
\[  L^{\tau\mu}_k(\omega(x)){\xi}^k_{\mu}+{\xi}^r{\partial}_r{\omega}^{\tau}(x)=0\hspace{5mm},\hspace{5mm}
- \frac{\partial L^{\tau\mu}_k(\omega(x))}{\partial u^{\sigma}}{\Omega}^{\sigma}(x){\xi}^k_{\mu}+{\xi}^r{\partial}_r{\Omega}^{\tau}(x)=0  \]
The first condition brings at once ${\xi}_q\in R_q$ and we obtain therefore the following central local result for $F_0$:  \\
\[      {\Upsilon}_0=\{ \Omega\in F_0{\mid} -\frac{\partial L^{\tau\mu}_k(\omega(x))}{\partial u^{\sigma}}{\Omega}^{\sigma}{\xi}^k_{\mu}+{\xi}^r{\partial}_r{\Omega}^{\tau}=0\hspace{3mm}, \forall {\xi}_q\in R_q \}  \]
As usual, the study of this system can be cut into two parts. First of all, we have to look for:\\
\[    E_0=\{  \Omega\in F_0 {\mid} \frac{\partial L^{\tau\mu}_k(\omega(x))}{\partial u^{\sigma}}{\Omega}^{\sigma}{\xi}^k_{\mu}=0, \forall {\xi}^0_q\in R^0_q  \}  \subseteq F_0 \]
In a symbolic way with $pri(R^0_q)=\{{\xi}^{\sigma}\}$ and $par(R^0_q)=\{{\xi}^{m+r}\}$, we get:  \\
\[ \frac{\partial L_{\sigma}}{\partial u}\Omega {\xi}^{\sigma}+\frac{\partial ({\cal{E}}^{\sigma}_{m+r}L_{\sigma})}{\partial u}\Omega{\xi}^{m+r}=0\]
 whenever ${\xi}^{\sigma}+{\cal{E}}^{\sigma}_{m+r}{\xi}^{m+r}=0$ and thus $\frac{\partial {\cal{E}}}{\partial u}\Omega=0$ in agrement with Proposition 5.17.\\
Then we have to introduce an $R_q$-connection, that is a section ${\chi}_q\in T^*\otimes R_q$ over $id_T\in T^*\otimes T$ such that:  \\
\[   - L^{\tau\mu}_k(\omega (x)){\chi}^k_{\mu,i}+{\partial}_i{\omega}^{\tau}(x)=0\Rightarrow 
   - \frac{\partial L^{\tau\mu}_k(\omega(x))}{\partial u^{\sigma}}{\chi}^k_{\mu,i}{\Omega}^{\sigma}+{\partial}_i{\Omega}^{\tau}=0  \]
and it just remains to study this last system for ${\Omega}^{\tau}(x)=A^{\alpha}(x)W^{\tau}_{\alpha}(\omega (x))$ as it does not depend on the choice of the connection ${\chi}_q$. Substituting, we obtain successively:  \\
\[ (- \frac{\partial L^{\tau\mu}_k}{\partial u^{\sigma}}{\chi}^k_{\mu,i}W^{\sigma}_{\alpha}+\frac{\partial W^{\tau}_{\alpha}}{\partial u^{\sigma}}{\partial}_i{\omega}^{\sigma})A^{\alpha}+({\partial}_iA^{\alpha})W^{\tau}_{\alpha}=0   \]
\[ ([L^{\mu}_k,W_{\alpha}])^{\tau}{\chi}^k_{\mu,i}A^{\alpha}+({\partial}_iA^{\alpha})W^{\tau}_{\alpha}=0\]
As $[L^{\mu}_k,W_{\alpha}]=0$ and the action is effective, we finally obtain ${\partial}_iA^{\alpha}=0$ that is $A=cst$.  \\
\hspace*{12cm}     Q.E.D.   \\

The following corollary is a direct consequence of the above proposition and explains many  classical results as we shall see in the last computational section. \\

\noindent
{\bf COROLLARY  5.27}: We have the relation:  \\
\[  \Theta=\{\xi\in T {\mid} {\cal{D}}\xi\equiv {\cal{L}}(\xi)\omega=0\}\Rightarrow N(\Theta)=\{ {\cal{D}}\xi\equiv {\cal{L}}(\xi)\omega=AW(\omega)\in {\Upsilon}_0  \}      \]

\noindent
{\bf REMARK  5.28}: If we set ${\tilde{R}}_q={\pi}^{q+1}_q({\tilde{R}}_{q+1})$, then ${\tilde{R}}_{q+1}$ is obtained by a procedure with two steps. First of all, we obtain ${\tilde{R}}_q=\{{\xi}_q\in J_q(T) {\mid} L({\xi}_q)\omega =A(x)W(\omega)\in E_0 \} $ by eliminating the infinitesimal parameters $A$ by means of pure linear algebra. Then, we have to take into account that $A=cst$ in the above corollary and another reason for understanding that $g_{q+1}={\tilde{g}}_{q+1}$ though we have only 
$g_q\subseteq {\tilde{g}}_q$ in general. Such a situation is well known in physics where the Poincar\'{e} group is of codimension one in its normalizer which is the Weyl group obtained by adding dilatations.  \\

According to the definition of the Janet bundles at the end of section 4, using the inclusion ${\wedge}^rT^*\otimes S_qT^*\otimes T\subset {\wedge}^rT^*\otimes J^0_q(T)\subset {\wedge}^rT^*\otimes J_q(T)$, in the case of an involutive system $R_q\subset J_q(T)$ of infinitesimal transitive Lie equations, we have:  \\
\[       F_r= {\wedge}^rT^*\otimes J^0_q(T)/({\wedge}^rT^*\otimes R^0_q+\delta ({\wedge}^{r-1}T^*\otimes S_{q+1}T^*\otimes T))     \]
As $R^0_q$ is associated with $R_q$ and $J^0_q(T)$ is associated with $R_q\subset J_q(T)$, we obtain:  \\

\noindent
{\bf LEMMA  5.29}: The Janet bundles are associated with $R_q$ and we set ${\Upsilon}_r=\Upsilon (F_r)$.  \\

From this lemma we shall deduce the following important but difficult theorem:  \\

\noindent
{\bf THEOREM  5.30}: The first order operators ${\cal{D}}_r:F_{r-1}\rightarrow F_r$ induce maps ${\cal{D}}_r:{\Upsilon}_{r-1}\rightarrow {\Upsilon}_r$ in the {\it deformation sequence}: \\
\[    0\longrightarrow Z(\Theta)\longrightarrow C(\Theta)\stackrel{{\cal{D}}}{\longrightarrow}{\Upsilon}_0\stackrel{{\cal{D}}_1}{\longrightarrow}{\Upsilon}_1\stackrel{{\cal{D}}_2}{\longrightarrow} ... \stackrel{{\cal{D}}_n}{\longrightarrow} {\Upsilon}_n\longrightarrow 0    \]
which is locally described by finite dimensional vector spaces and linear maps.   \\

\noindent
{\it Proof}: Any section of $F_{r-1}$ can be lifted to a section $A^{r-1}_q\in {\wedge} ^{r-1}T^*\otimes J_q(T)$ modulo ${\wedge}^{r-1}T^*\otimes R_q+\delta ({\wedge}^{r-2}T^*\otimes S_{q+1}T^*\otimes T)$ where the second component is in the image of $D$. Then we can lift again this section to a section $A^{r-1}_{q+1}\in {\wedge}^{r-1}T^*\otimes J_{q+1}(T)$, modulo a section of ${\wedge}^{r-1}T^*\otimes R_{q+1}+D({\wedge}^{r-2}T^*\otimes J^q_{q+2}(T))+{\wedge}^{r-1}T^*\otimes S_{q+1}T^*\otimes T$ where $J^q_{q+2}$ is the kernel of the projection ${\pi}^{q+2}_q:J_{q+2}(T)\rightarrow J_q(T)$. The image by ${\cal{D}}_r$ is obtained by applying $D:{\wedge}^{r-1}T^*\otimes J_{q+1}(T)\rightarrow {\wedge}^rT^*\otimes J_q(T)$ and projecting $DA^{r-1}_{q+1}\in {\wedge}^rT^*\otimes J_q(T)$ thus obtained to $F_r$ while taking into account successively the restriction $D:{\wedge}^{r-1}T^*\otimes R_{q+1}\rightarrow {\wedge}^rT^*\otimes R_q$, the fact that $D\circ D=D^2=0$ and the restriction $\delta$ providing an element in $\delta ({\wedge}^{r-1}T^*\otimes S_{q+1}T^*\otimes T)$. Of course, with such a choice we need to use associations with $R_{q+1}$ even though finally only $R_q$ is involved because of the above lemma. The reason is that the restriction of $D$ to $J^0_{q+1}(T)$ has an image in $T^*\otimes J_q(T)$ and {\it not} in $T^*\otimes J^0_q(T)$. Hence, starting with a section of ${\Upsilon}_{r-1}\subset F_{r-1}$, we must have:  \\
\[   L({\xi}_{q+1})A^{r-1}_q\in {\wedge}^{r-1}T^*\otimes R_q+\delta ({\wedge}^{q-2}T^*\otimes S_{q+1}T^*\otimes T)  \hspace{1cm}  \forall {\xi}_{q+1}\in R_{q+1}={\rho}_1(R_q)      \]
Accordingly, looking at the right member in the formula of Proposition 5.13, the section in ${\wedge}^{r-1}T^*\otimes R_q$ is contracted with $r-1$ vectors in order to provide a section of $R_q$ and the skewsymmetrized summation finally produces a section of ${\wedge}^rT^*\otimes R_q$. Similarly, as $\delta ({\wedge}^{r-2}T^*\otimes S_{q+1}T^*\otimes T)\subset {\wedge}^{r-1}T^*\otimes S_qT^*\otimes T$ we should obtain a section of ${\wedge}^rT^*\otimes S_qT^*\otimes T$ but it is not evident at all that such a section is in $\delta ({\wedge}^{r-1}T^*\otimes S_{q+1}T^´\otimes T)$. For this, let us notice that, if the multi-index $I$ has length $r-2$ and ${\mid}\mu{\mid}=q$, we have the local formula:  \\
\[ {\alpha}_{\mu+1_i,I,j}dx^i\wedge (dx^j\wedge dx^I)=- {\alpha}_{\mu+1_i,I,j}dx^j\wedge (dx^i\wedge dx^I)      \]   
explaining why the following diagram is commutative and exact: \\
\[     \begin{array}{cccl}
T^*\otimes {\wedge}^{r-2}T^*\otimes S_{q+1}T^*&\stackrel{(-1)^{r-2}\delta}{\longrightarrow} & {\wedge}^{r-1}T^*\otimes S_{q+1}T^*  & \longrightarrow 0   \\
\downarrow \delta &    & \downarrow \delta   &         \\
T^*\otimes {\wedge}^{r-1}T^*\otimes S_qT^* &   \stackrel{(-1)^{r-1}\delta}{\longrightarrow} & {\wedge}^rT^*\otimes S_qT^*    &   \longrightarrow 0 
\end{array}     \]
where the upper map induced by $\delta : {\wedge}^{r-2}T^*\otimes T^*\rightarrow {\wedge}^{r-1}T^*$ and the lower map induced by $\delta : {\wedge}^{r-1}T^*\otimes T^*\rightarrow {\wedge}^rT^*$ are both epimorphisms while the left vertical map is induced by $\delta : {\wedge}^{r-2}T^*\otimes S_{q+1}T^*\rightarrow {\wedge}^{r-1}T^*\otimes S_qT^*$.  \\
We obtain therefore:  \\
\[   DL({\xi}_{q+2})A^{r-1}_{q+1}-L({\xi}_{q+1})DA^{r-1}_{q+1}\in {\wedge}^rT¬*\otimes R_q + \delta ( {\wedge}^{r-1}T^*\otimes S_{q+1}T^*\otimes T)      \]
Now, from the construction of ${\cal{D}}_r$ already explained, we have:   \\
\[    F_{r-1}={\wedge}^{r-1}T^*\otimes J_{q+1}(T)/({\wedge}^{r-1}T^*\otimes R_{q+1}+D({\wedge}^{r-2}T^*\otimes J^q_{q+2}(T))+{\wedge}^{r-1}T^*\otimes S_{q+1}T^*\otimes T)   \]
It follows that:   \\
\[  L({\xi}_{q+2})A^{r-1}_{q+1}\in {\wedge}^{r-1}T^*\otimes R_{q+1} + D({\wedge}^{r-2}T^*\otimes J^q_{q+2}(T)) + {\wedge}^{r-1}T^*\otimes S_{q+1}T^*\otimes T   \]
and thus:   \\
\[  DL({\xi}_{q+2})A^{r-1}_{q+1}\in {\wedge}^rT^*\otimes R_q + \delta ({\wedge}^{r-1}T^*\otimes S_{q+1}T^*\otimes T)      \]
that is finally:   \\
\[   L({\xi}_{q+1})DA^{r-1}_{q+1}\in {\wedge}^rT^*\otimes R_q + \delta ({\wedge}^{r-1}T^*\otimes S_{q+1}T^*\otimes T)         \]
Accordingly, $DA^{r-1}_{q+1}\in {\wedge}^rT^*\otimes J_q(T)$ is the representative of a section of ${\Upsilon}_r$.   \\
\hspace*{12cm}     Q.E.D.     \\

\noindent
{\bf DEFINITION  5.31}: The deformation sequence is not necessarily exact and only depends on 
$R_q$. We can therefore define as usual {\it coboundaries} $B_r(R_q)$, {\it cocycles} $Z_r(R_q)$ and {\it cohomology groups} $H_r(R_q)=Z_r(R_q)/B_r(R_q)$ with $B_r(R_q)\subseteq Z_r(R_q)\subseteq {\Upsilon}_r$ for $r=0,1,...,n$.  \\

\noindent
{\bf PROPOSITION  5.32}: We have $B_0(R_q)=C(\Theta)/Z(\Theta), Z_0(R_q)=N(\Theta)/\Theta$ and the short exact sequence:   \\
\[    0 \longrightarrow  C(\Theta)/Z(\Theta) \longrightarrow N(\Theta)/\Theta \longrightarrow H_0(R_q) \longrightarrow 0            \]

\noindent
{\it Proof}: The monomorphism on the left is induced by the inclusion $C(\Theta)\subset N(\Theta)$ because $Z(\Theta)=\Theta \cap C(\Theta)$. Then $B_0(R_q)$ is the image of ${\cal{D}}$ in ${\Upsilon}_0$ because of Proposition 5.23 and Proposition 5.24. Finally, $Z_0(R_q)=N(\Theta)/\Theta$ is just a way to rewrite Corollary 5.27 while taking into account the fact that ${\cal{D}}_1\circ{\cal{D}}=0$. The cocycle condition just tells that the label transformations induced by the normalizer do not change the structure constants because the Vessiot structure equations are invariant under {\it any} natural transformation. This result should be compared to Lemma 2.7.  \\
\hspace*{12cm}   Q.E.D.  \\

In order to generalize Proposition 5.26 and to go further on, we need a few more concepts from differential geometry. \\

\noindent
{\bf DEFINITION  5.33}: A chain ${\cal{E}}\stackrel{\Phi}{\longrightarrow} {\cal{E}}' \stackrel{\Psi}{\longrightarrow} {\cal{E}}''$ of fibered manifolds is said to be a {\it sequence with respect to a section} $f''$ {\it of } ${\cal{E}}''$ if $im(\Phi)=ker_{f''}(\Psi)$, that is with local coordinates $(x,y)$ on ${\cal{E}}$, $(x,y')$ on ${\cal{E}}'$, $(x,y'')$ on ${\cal{E}}''$ and $y'=\Phi (x,y), y''=\Psi (x,y')$, we have $\Psi(x,\Phi(x,y))\equiv f''(x), \forall (x,y)\in {\cal{E}}$.  \\

Differentiating this identity with respect to $y$, we obtain:  \\
\[   \frac{\partial \Psi}{\partial y'}(x,\Phi (x,y)). \frac{\partial \Phi}{\partial y}(x,y) \equiv 0  \]
and we have therefore a sequence $V({\cal{E}})\stackrel{V(\Phi)}{\longrightarrow } V({\cal{E}}') \stackrel{V(\Psi)}{\longrightarrow  } V({\cal{E}}'' ) $ of vector bundles pulled back over ${\cal{E}}$ by reciprocal images.  \\

\noindent
{\bf DEFINITION  5.34}: A sequence of fibered manifolds is said to be an {\it exact sequence} if $im(\Phi)=ker_{f''}(\Psi)$ {\it and} the corresponding vertical sequence of vector bundles is exact. \\

\noindent
{\bf PROPOSITION 5.35}: If ${\cal{E}}, {\cal{E}}', {\cal{E}}'' $ are affine bundles over $X$ with corresponding model vector bundles $E, E', E'' $ over $X$, a sequence of such affine bundles is exact if and only if the corresponding sequence of model vector bundles is exact. In that case, there is an exact sequence ${\cal{E}}\stackrel{\Phi}{\longrightarrow} {\cal{E}}' \longrightarrow E''$ which allows to avoid the use of a section of ${\cal{E}}'' $ while replacing it by the zero section of $E'' $.  \\

\noindent
{\it Proof}: We have successively $y'=\Phi (x,y)= A(x)y+B(x), y''= \Psi (x,y')=C(x)y'+D(x)$ and by composition $y''= C(x)A(x)y+C(x)B(x)+D(x)=f''(x)$. Accordingly, we {\it must} have $C(x)A(x)\equiv 0, 
C(x)B(x)+D(x)=f''(x), \forall x\in X$ and we obtain the following commutative diagram:  \\

\[     \begin{array}{rcccc}
  E & \stackrel{V(\Phi)}{\longrightarrow} & E' & \stackrel{V(\Psi)}{\longrightarrow} & E''    \\
  \vdots  &    &   \vdots   & \nearrow & \vdots    \\
  {\cal{E}}  & \stackrel{\Phi}{\longrightarrow} & {\cal{E}}' &  \stackrel{\Psi}{\longrightarrow} & {\cal{E}}''  \\
  \pi \downarrow &   &  {\pi}' \downarrow \hspace{3mm}  &   &{\pi}'' \downarrow \uparrow f''   \\
  X  &  =  & X  &  =  &  X
  \end{array}     \]
  
If $(x,y')\in {\cal{E}}'$ is such that $C(x)y'+D(x)=f''(x)$, we obtain by substraction $C(x)(y'-B(x))=0 \Leftrightarrow C(x)(y'-\Phi(x,y))=0$ with $(x,y'-\Phi(x,y))\in E'$. Supposing the model sequence exact, we may find $(x,v)\in E$ such that $y'-(A(x)y+B(x))=A(x) v\Leftrightarrow y'=A(x)(y+v)+B(x)$ and thus $(x,y')\in im(\Phi)$. The converse is similar and left to the reader.\\
Finally, we just need to notice that $(x,y'')\in {\cal{E}}'' \Leftrightarrow (x,y''-f''(x))\in E'' , \forall x\in X$ with $y''-f''(x)=C(x)(y'-B(x))=C(x)(y'-\Phi(x,y))$ where$(x,y)\in {\cal{E}}$ is {\it any} point over $x\in X$. \\
\hspace*{12cm}      Q.E.D.   \\

Coming back to the construction of the Vessiot structure equations from the knowledge of the generating differential invariants at order $q$, we have already exhibited the system ${\cal{B}}_1\subset J_1({\cal{F}})$ locally defined by affine equations of the form $I(u_1)\equiv A(u)u_x+B(u)=0$. The symbol ${\cal{H}}_1\subset T^*\otimes V({\cal{F}})$ of this system is defined by linear equations of the form $A(u)v_x=0$ and the last proposition provides at once the following commutative and exact diagram of affine bundles and model vector bundles.  \\

\[   \begin{array}{rcccccl}
0\longrightarrow & {\cal{H}}_1 & \longrightarrow & T^*\otimes V({\cal{F}}) & \longrightarrow &{\cal{F}}_1 & \longrightarrow 0  \\
      &    \vdots  &    &   \vdots   &  \nearrow & \parallel &    \\
      0  \longrightarrow & {\cal{B}}_1 &  \longrightarrow & J_1({\cal{F}}) & \longrightarrow & {\cal{F}}_1 & \longrightarrow 0  \\
      & \downarrow &   &  \downarrow  & & \downarrow  &   \\
      &  {\cal{F}}  & =  & {\cal{F}}  & =  &  {\cal{F}}  & 
 \end{array}     \]

More generally, having in mind the diagram at the end of section 4, we may define the {\it nonlinear Janet bundles} ${\cal{F}}_r={\wedge}^rT^*\otimes V({\cal{F}})/\delta ({\wedge}^{r-1}T^*\otimes {\cal{H}}_1)$ as a family of natural vector bundles over ${\cal{F}}$. The following commutative diagram of reciprocal images generalizes the result of Proposition 5.26:  \\

\[  \begin{array}{rcccl}
   & F_r  &  \longrightarrow & {\cal{F}}_r &      \\
   section \hspace{2mm}of\hspace{2mm} {\Upsilon}_r & \downarrow\uparrow & &\downarrow\uparrow & equivariant \hspace{2mm} section  \\
     &  X  &  \stackrel{\omega}{\longrightarrow} & {\cal{F}}&      
     \end{array}    \]

As a byproduct, we may generalize Proposition 5.26 by saying that any infinitesimal change of source $\bar{x}=x+t\xi(x)+ ... $ can be lifted to ${\cal{F}}_r$ with ${\bar{u}}^{\tau}=u^{\tau}+t{\xi}^k_{\mu}(x)L^{\tau\mu}_k(u)+ ... , {\bar{v}}^{\alpha}=v^{\alpha}+t{\xi}^k_{\mu}(x)M^{\alpha}_{\beta}{\mid}^{\mu}_k(u)v^{\beta}+ ... $.\\

\noindent
{\bf REMARK  5.36}: No classical technique could provide this result because all the known methods of computer algebra do construct the Janet sequence "{\it step by step}" and never "{\it as a whole}", that is from the Spencer operator ([36], p 391).  \\

\noindent
{\bf PROPOSITION  5.37}: The affine bundle:   \\
\[     A(R^0_q)=\{ {\chi}_q\in T^*\otimes J_q(T)/(T^*\otimes R^0_q + \delta (S_{q+1}T^*\otimes T)){\mid} {\chi}_0=id_T\in T^*\otimes T \}   \]
is modelled on $F_1$ and associated with $R_q$.  \\

\noindent
{\it Proof}: First of all, $A(R^0_q)$ is modelled on $T^*\otimes J^0_q(T)/(T^*\otimes R^0_q+\delta (S_{q+1}T^*\otimes T))=F_1$ and it just remains to prove that the affine bundle:  \\
\[        A_q=\{  {\chi}_q\in C_1(T){\mid} {\chi}_0=id_T\in T^*\otimes T  \} = id_T^{-1}(C_1(T))   \]
is a natural affine bundle over $X$ modeled on the vector bundle $C^0_1(T)=T^*\otimes J^0_q(T)/\delta (S_{q+1}T^*\otimes T)$ and associated with $J_q(T)$, thus with $R_q$ too because $T^*\otimes R^0_q$ is ass
ociated with $R_q$. \\For this it is sufficient to observe the transition laws of $T^*\otimes J_q(T)$ when $\bar{x}=\varphi(x)$, namely:  \\
\[  \begin{array}{lccl}
T^*\otimes T   &  {\bar{\chi}}^l_{,r}{\partial}_i{\varphi}^r & = & {\chi}^k_{,i}{\partial}_k{\varphi}^l \\
  T^*\otimes J_1(T) & {\bar{\chi}}^l_{s,r}{\partial}_i{\varphi}^r{\partial}_j{\varphi}^s &= & {\chi}^k_{j,i}{\partial}_k{\varphi}^l+{\chi}^k_{,i}{\partial}_{jk}{\varphi}^l    \\
 T^*\otimes J_q(T)  &  {\bar{\chi}}^l_{{r_1}... {r_q},r}{\partial}_{i_1}{\varphi}^{r_1}... {\partial}_{i_q}
 {\varphi}^{r_q}{\partial}_i{\varphi}^r +i ... &=& {\chi}^k_{\mu,i}{\partial}_k{\varphi}^l+ ... +{\chi}^k_{,i}{\partial}_{\mu+1_k}{\varphi}^l   
\end{array}   \]
with $\mu$ replaced by $(i_1, ... ,i_q)$ when ${\mid}\mu{\mid}=q$. Setting ${\chi}^k_{,i}={\delta}^k_i$, we get: \\
\[  id_T^{-1}(T^*\otimes J_1(T)) \hspace{4cm}  {\bar{\chi}}^l_{s,r}{\partial}_i{\varphi}^r{\partial}_j{\varphi}^s={\chi}^k_{j,i}{\partial}_k{\varphi}^l+{\partial}_{ij}{\varphi}^l    \]
\[  id_T^{-1}(T^*\otimes J_q(T)) \hspace{2cm}  {\bar{\chi}}^l_{{r_1}... {r_q},r}{\partial}_{i_1}{\varphi}^{r_1}... {\partial}_{i_q}{\varphi}^{r_q}{\partial}_i{\varphi}^r + ... = {\chi}^k_{\mu,i}{\partial}_k{\varphi}^l+ ... +{\partial}_{\mu+1_i}{\varphi}^l   \]
and thus:  \\
\[  \begin{array}{rcl}
i(\eta)L({\xi}_{q+1}){\chi}_q & = & L({\xi}_{q+1}){\chi}_q(\eta)-{\chi}_q(L({\xi}_1)\eta)   \\
                                                 & = & [{\xi}_q,{\chi}_q(\eta)]+ i(\eta)D{\xi}_{q+1}-{\chi}_q(L({\xi}_1)\eta)
                                                 \end{array}    \]                                                  
because  ${\chi}_0(\eta)=\eta$. If we choose ${\xi}_{q+1},{\bar{\xi}}_{q+1}\in J_{q+1}(T)$ over ${\xi}_q\in J_q(T)$, the difference will be $i(\eta)D({\bar{\xi}}_{q+1}-{\xi}_{q+1})=-i(\eta)\delta ({\bar{\xi}}_{q+1}-{\xi}_{q+1})=0$ by residue because we are in $C_1(T)$. We notice that no formal integrability assumption is needed for $R_q$. \\
\hspace*{12cm}    Q.E.D.    \\                                               
                                               
Similarly, we obtain:  \\

\noindent
{\bf PROPOSITION  5.38}: $id_q^{-1}(J_1({\Pi}_q))$ is an affine natural bundle of order $q+1$ over $X$, modelled on $T^*\otimes J_q(T)$.  \\

\noindent
{\it Proof}: We only prove the proposition when $q=1$ as the remaining of he proof is similar to the previous one. Indeed, if $\bar{x}=\varphi (x), \bar{y}=\psi (y) $ are the changes of coordinates on $X\times Y$, we get:  \\
\[   \begin{array}{lrcl}
{\Pi}_1  &  {\bar{y}}^u_r{\partial}_i{\varphi}^r(x) & = & \frac{\partial {\psi}^u}{\partial y^k}(y)y^k_i   \\
J_1({\Pi}_1)  &  {\bar{y}}^u_{r,s} {\partial}_i{\varphi}^r(x){\partial}_j{\varphi}^s(x)+{\bar{y}}^u_r{\partial}_{ij}{\varphi}^r(x) & = &  \frac{\partial {\psi}^u}{\partial y^k}(y)y^k_{i,j}+\frac{{\partial}^2{\psi}^u}{\partial y^k\partial y^l}(y)y^k_iy^l_{,j}  \\
id_1^{-1}(J_1({\Pi}_1))  &  {\bar{y}}^u_{r,s}{\partial}_i{\varphi}^r{\partial}_j{\varphi}^s +{\partial}_{ij}{\varphi}^u & = & {\partial}_k{\varphi}^uy^k_{i,j}+{\partial}_{ik}{\varphi}^uy^k_{,j}   
\end{array}   \]                                            
Then, we just need to  set $y=x, y^k_j={\delta}^k_j$ and compare to $T^*\otimes J_1(T)$ above.\\
Finally, as $J_1({\Pi}_q)$ is an affine bundle over ${\Pi}_q$ modelled on $T^*\otimes V({\Pi}_q)$, then $id_q^{-1}(J_1({\Pi}_q))$ is an affine bundle over $X$, modelled on $id_q^{-1}(T^*\otimes V({\Pi}_q))=T^*\otimes id_q^{-1}(V({\Pi}_q))=T^*\otimes J_q(T)$.  \\
\hspace*{12cm}   Q.E.D.    \\

\noindent
{\bf PROPOSITION  5.39}: We have the following exact sequences of affine bundles over $X$ and model vector bundles:  \\                            
\[   \begin{array}{rcccccl}
0 \longrightarrow & T^*\otimes R_q & \longrightarrow & T^*\otimes J_q(T) & \stackrel{V(\Phi)}{\longrightarrow} & T^*\otimes F_0 & \longrightarrow 0  \\
   &  \vdots  &   &  \vdots   &  &  \vdots  &          \\
   0 \longrightarrow  & id_q^{-1}(J_1({\cal{R}}_q)) &  \longrightarrow & id_q^{-1}(J_1({\Pi}_q)) & \longrightarrow & {\omega}^{-1}(J_1({\cal{F}}))  & \longrightarrow 0  \\
      &  \downarrow  &   &  \downarrow  &  &  \downarrow   &    \\
        &  X   &  =  & X &   =  &  X  & 
\end{array}      \]

\noindent
{\it Proof}: We just need to use $1\leq \mid\mu\mid \leq q$ and set $y_q=id_q(x)$ in order to obtain: \\
\[  u^{\tau}={\Phi}^{\tau}(\omega(y), y_{\mu}) \stackrel{id_q^{-1}}{\longrightarrow}{\omega}^{\tau}(x)  \]
\[ u^{\tau}_i=\frac{\partial {\Phi}^{\tau}}{\partial u^{\sigma}}(y_q)\frac{\partial {\omega}^{\sigma}}{\partial y^k}(y)y^k_{,i}+\frac{\partial {\Phi}^{\tau}}{\partial y^k_{\mu}}(y_q)y^k_{\mu,i} \stackrel{id_q^{-1}}{\longrightarrow } -L^{\tau\mu}_k({\omega}(x))y^k_{\mu,i}+y^k_{,i}{\partial}_k{\omega}^{\tau}(x)  \]
in a coherent way with the upper model sequence and the various distinguished sections.  \\
\hspace*{12cm}    Q.E.D.    \\

\noindent
{\bf PROPOSITION  5.40}: We have the following exact sequences of affine bundles over $X$ and model vector bundles:  \\                            
\[ \begin{array}{rcccccl}
0 \longrightarrow & T^*\otimes R^0_q  & \longrightarrow & T^*\otimes J^0_q(T) & \stackrel{V(\Phi)}{\longrightarrow} & T^*\otimes F_0   &  \longrightarrow 0  \\
   &  \vdots  &   & \vdots  &  & \vdots  &      \\
0 \longrightarrow  & id_T^{-1}(T^*\otimes R_q) & \longrightarrow & id_T^{-1}(T^*\otimes J_q(T)) & \longrightarrow &   {\omega}^{-1}(J_1({\cal{F}})) &  \longrightarrow 0  \\
   &  \downarrow &  & \downarrow &  & \downarrow  &       \\
   &  X  & =   &  X  &  =  &  X  &  
   \end{array}   \]

\noindent
{\it Proof}: As $id_q$ is a section of ${\Pi}_q$, then $j_1(id_q)$ is a section of $J_1({\Pi}_q)$ over $id_q$ and thus a section of $id_q^{-1}(J_1({\Pi}_q))$ denoted by $id_{q,1}$. If ${\chi}_q=({\delta}^k_i,{\chi}^k_{\mu,i})$ with $1\leq \mid\mu\mid \leq q$ is a section of $id_T^{-1}(T^*\otimes J_q(T))$, that is a $J_q(T)$-connection, it is of course a section of $T^*\otimes J_q(T)$ and we may consider the section 
$y_{q,1}=id_{q,1}+{\chi}_q$ of $id_q^{-1}(J_1({\Pi}_q))$. The image in ${\omega}^{-1}(J_1({\cal{F}}))$ is 
$u^{\tau}_i=(-L^{\tau\mu}_k(\omega (x)){\chi}^k_{\mu,i}+{\partial}_i{\omega}^{\tau}(x))+{\partial}_i{\omega}^{\tau}(x)$ and we obtain: $u^{\tau}_i={\partial}_i{\omega}^{\tau}(x) \Leftrightarrow u_1=j_1(\omega)(x)\Leftrightarrow {\chi}_q\in id_T^{-1}(T^*\otimes R_q)$.  \\
\hspace*{12cm}     Q.E.D.  \\

\noindent
{\bf THEOREM  5.41}:  \hspace{5mm} $A(R^0_q)\simeq {\omega}^{-1}(J_1({\cal{F}})/{\cal{B}}_1) $   \\

\noindent
{\it Proof}: The affine sub-bundle ${\cal{B}}_1\subset J_1({\cal{F}})$ over ${\cal{F}}$ is the imager of the first prolongation ${\rho}_1(\Phi): {\Pi}_{q+1} \longrightarrow J_1({\cal{F}})$ where ${\Pi}_{q+1}$ is an affine bundle over ${\Pi}_q$ modelled on $S_{q+1}T^*\otimes V(X\times X)$. We may proceed similarly by introducing the affine bundle $id_q^{-1}({\Pi}_{q+1})$ which is modelled on $id_q^{-1}(S_{q+1}T^*\otimes V(X\times X))=S_{q+1}T^*\otimes id^{-1}(V(X\times X))=S_{q+1}T^*\otimes T$ and obtain the commutative and exact diagram of affine bundles and model vector bundles where the model sequence is the symbol sequence that has been used in order to introduce $F_1$:  \\

\[    \begin{array}{rcccccl}
 & S_{q+1}T^*\otimes T  &\longrightarrow & T^*\otimes F_0 & \longrightarrow& F_1 & \longrightarrow  0  \\
&   \vdots   &   &  \vdots  &   &  \vdots   &   \\
 &  id_q^{-1}({\Pi}_{q+1}) & \longrightarrow & {\omega}^{-1}(J_1({\cal{F}})) & \longrightarrow &  {\omega}^{-1}(J_1({\cal{F}})/{\cal{B}}_1) & \longrightarrow  0  \\
  &  \downarrow &   & \hspace{11mm}\downarrow \uparrow j_1(\omega) &   &  \downarrow  \\
   &  X  &  =  &  X  & =   &  X  &    
\end{array}    \]

With $1\leq \mid\mu\mid \leq q$, the first morphism of affine bundles is now described by:  \\
\[  u^{\tau}={\omega}^{\tau}(x),u^{\tau}_i=\frac{\partial {\Phi}^{\tau}}{\partial u^{\sigma}}(y_q)\frac{\partial {\omega}^{\sigma}}{\partial y^k}(y)y^k_i+\frac{\partial {\Phi}^{\tau}}{\partial y^k_{\mu}}(y_q)y^k_{\mu+1_i }
\stackrel{id_q^{-1}}{\longrightarrow } -\sum_{\mid\mu\mid=q}L^{\tau\mu}_k(\omega (x))y^k_{\mu+1_i}+{\partial}_i{\omega}^{\tau}(x)  \]
because $x$ does not appear in ${\Phi}^{\tau}(y_q)$, in such a way that $ u^{\tau}_i-{\partial}_i{\omega}^{\tau}(x)=- \sum_{\mid\mu\mid=q}L^{\tau\mu}_k(\omega (x))y^k_{\mu+1_i}\in T^*\otimes F_0$. It just remains to notice that $\delta : S_{q+1}T^*\otimes T{\longrightarrow} T^*\otimes S_qT^*\otimes T$ is a monomorphism.   \\
\hspace*{12cm}   Q.E.D.     \\

If ${\chi}_q\in T^*\otimes J_q(T)$, with projections ${\chi}_0=id_T\in T^*\otimes T$ and ${\chi}_1\in T^*\otimes J_1(T)$, is a representative of an element $c\in A(R^0_q)$, we may choose a lift ${\chi}_{q+1}\in T^*\otimes J_{q+1}(T)$ and define $\frac{1}{2}\{{\chi}_q,{\chi}_q\}\in {\wedge}^2T^*\otimes J^0_q(T)$ by the formula:  \\
\[  \frac{1}{2} \{ {\chi}_q,{\chi}_q \}(\eta,\zeta)=\{{\chi}_{q+1}(\eta),{\chi}_{q+1}(\zeta)\}-{\chi}_q(\{{\chi}_1(\eta),{\chi}_1(\zeta)\})   \]

\noindent
{\bf PROPOSITION  5.42}: The map ${\chi}_q \longrightarrow \frac{1}{2}\{{\chi}_q,{\chi}_q\}$ provides a well defined map $c \longrightarrow \frac{1}{2}\{c,c\}$ from elements of $A(R^0_q)$ invariant by $R^0_q$ to $F_2={\wedge}^2T^*\otimes J^0_q(T)/({\wedge}^2T^*\otimes R^0_q+\delta(T^*\otimes S_{q+1}T^*\otimes T))$. \\

\noindent
{\it Proof}: First of all, if we change the lift ${\chi}_{q+1}$ to ${\bar{\chi}}_{q+1}$, then ${\bar{\chi}}_{q+1}-{\chi}_{q+1}\in T^*\otimes S_{q+1}T^*\otimes T$ and the difference will therefore be in $\delta (T^*\otimes S_{q+1}T^*\otimes T)$ with a zero projection in $F_2$.  \\
Then, if we modify ${\chi}_q$ by an element $M_q\in \delta(S_{q+1}T^*\otimes T)\subset T^*\otimes  S_qT^*\otimes T\subset T^*\otimes J_q(T)$, the difference in $J_{q-1}(T)$ will be only produced by 
$\{{\chi}_q(\eta),{\chi}_q(\zeta)\}$ and will be $\{M_q(\eta),{\chi}_q(\zeta)\}+\{{\chi}_q(\eta),M_q(\zeta)\}=\delta M_q(\eta,\zeta)=0$ because $\delta\circ\delta=0$. However, when lifting $M_q$ at order $q+1$, new terms will modify the quadratic application and, in particular, even one more if $q=1$, namely ${\chi}_1(\{M_1(\xi),{\chi}_1(\eta)\}+\{{\chi}_1(\xi),M_1(\eta)\})={\chi}_1(\delta M_1(\xi,\eta))={\chi}_1(0)=0$ and such a situation will not differ from the general one we now consider. According to the above results, the only pertubating term to study is:  \\
\[\{{\sigma}^{q-1}_{q+1}(\xi),{\chi}_{q+1}(\eta)\}+\{{\chi}_{q+1}(\xi),{\sigma}^{q-1}_{q+1}(\eta)\}-M_q(\{{\chi}_1(\xi),{\chi}_1(\eta)\})\in S_qT^*\otimes T   \]
and we just need to prove that such a term comes from an element in $\delta (T^*\otimes S_{q+1}T^*\otimes T)$. With more details and $\mid\mu\mid=q-1$, we get for the $6=3+3$ factors of ${\xi}^i,{\eta}^j$: \\
\[ ({\chi}^r_{t,i}M^k_{\mu+1_r+1_j}-{\chi}^r_{t,j}M^k_{\mu+1_r+1_i})+(M^r_{\mu+1_t+1_i}{\chi}^k_{r,j}-M^r_{\mu+1_t+1_j}{\chi}^k_{r,i})-M^k_{\mu+1_r+1_t}({\chi}^r_{i,j}-{\chi}^r_{j,i})  \]
Setting:  \\
\[ N_{q+1}=(N^k_{\mu+1_t+1_j,i}=({\chi}^r_{t,i}M^k_{\mu+1_r+1_j}+{\chi}^r_{j,i}M^k_{\mu+1_r+1_t})-{\chi}^k_{r,i}M^r_{\mu+1_t+1_j}) \in T^*\otimes S_{q+1}T^*\otimes T  \]
where the first sum in the parenthesis insures the symmetry in $t/j$, a tedious but straightforward though unexpected calculation proves that all the above terms are just described by $\delta N_{q+1}$. \\
As another proof, we may also consider the $3$-cyclic sum obtained by preserving $\mu$ but replacing 
$(i,j,t)$ successively by $(j,t,i)$ and $(t,i,j)$ in order to discover that {\it all the terms disappear two by two}. It is then sufficient to use the exactness of the $\delta$-sequence:   \\
\[  0 \longrightarrow S_{q+2}T^*\otimes T\stackrel{\delta}{\longrightarrow} T^*\otimes S_{q+1}T^*\otimes T \stackrel{\delta}{\longrightarrow} {\wedge}^2T^*\otimes S_qT^´\otimes T\stackrel{\delta}{\longrightarrow} {\wedge}^3T^*\otimes S_{q-1}T^*\otimes T   \]
Hence it just remains to modify ${\chi}_q$ by ${\sigma}^0_q\in T^*\otimes R^0_q$ with lift ${\sigma}^0_{q+1}\in T^*\otimes J^0_{q+1}(T)$ and the difference will be:  \\
\[ \{{\chi}_{q+1}(\eta),{\sigma}^0_{q+1}(\zeta)\}+\{ {\sigma}^0_{q+1}(\eta),{\chi}_{q+1}(\zeta)\}+[{\sigma}^0_q(\eta),{\sigma}^0_q(\zeta)]\]
\[-{\chi}_q(\{{\sigma}^0_1(\eta),{\chi}_1(\zeta)\}+\{{\chi}_1(\eta),{\sigma}^0_1(\zeta)\})-{\sigma}^0_q(\{{\chi}_1(\eta),{\chi}_1(\zeta)\})   \]
As the third and the last terms already belong to $R^0_q$, we need just consider:  \\
\[  \{ {\chi}_{q+1}(\eta),{\sigma}^0_{q+1}(\zeta)\}+\{ {\sigma}^0_{q+1}(\eta),{\chi}_{q+1}(\zeta)\}-{\chi}_q(\delta {\sigma}^0_1(\eta,\zeta))      \]
Now we have:  \\
\[  i(\zeta)L({\sigma}^0_{q+1}(\eta)){\chi}_q=L({\sigma}^0_{q+1}(\eta)){\chi}_q(\zeta)-{\chi}_q(L({\sigma}^0_1(\eta))\zeta)=\{ {\sigma}^0_{q+1}(\eta),{\chi}_{q+1}(\zeta)\}-{\chi}_q(L({\sigma}^0_1(\eta))\zeta)  \]
and the above remainder becomes:\\
\[  i(\zeta)L({\sigma}^0_{q+1}(\eta)){\chi}_q-i(\eta)L({\sigma}^0_{q+1}(\zeta)  \in J^0_q(T)  \]
because:  \
\[      L({\sigma}^0_1(\eta))\zeta-L({\sigma}^0_1(\zeta))\eta-\delta{\sigma}^0_1(\eta,\zeta)=0  \]
Finally, $L({\xi}^0_q)c=0, \forall {\xi}^0_q\in R^0_q$ in $F_1$ $\Leftrightarrow$ $L({\xi}^0_{q+1}){\chi}_q\in T^*\otimes R^0_q+\delta (S_{q+1}T^*\otimes T)$. Hence, modifying ${\sigma}^0_{q+1}$ if necessary by an element in $T^*\otimes S_{q+1}T^*\otimes T$, the above remainder belongs to $R^0_q$. \\
It follows that $\frac{1}{2}\{{\chi}_q,{\chi}_q\}$ will be modified by an element in ${\wedge}^2T^*\otimes R^0_q+\delta (T^*\otimes S_{q+1}T^*\otimes T)$ that will not change the projection in $F_2$.  \\
\hspace*{12cm}    Q.E.D.   \\

If we start now with any system $R_q\subset J_q(T)$ of infinitesimal Lie equations satisfying $[R_q,R_q]\subset R_q$, we may choose {\it any} $R_q$-connection ${\chi}_q\in id_T^{-1}(T^*\otimes R_q)\subset id_T^{-1}(T^*\otimes J_q(T))$ and obtain by projection an element $c\in A(R^0_q)$ that we use in the following proposition.   \\

\noindent
{\bf PROPOSITION  5.43}: ${\pi}^{q+1}_q:R_{q+1}\longrightarrow R_q$ is surjective if and only if 
$L({\xi}_q)c=0, \forall {\xi}_q\in R_q$. In this case, we have $\frac{1}{2}\{c,c\}=0$ in $F_2$.  \\

\noindent
{\it Proof}: We have ${\chi}_0(\eta)=\eta$ and we may find ${\xi}_{q+1}\in J_{q+1}(T)$ such that:  \\
\[  \begin{array}{rcl}
 i(\eta)L({\xi}_{q+1}){\chi}_q &= &L({\xi}_{q+1}){\chi}_q(\eta)-{\chi}_q(L({\xi}_1)\eta)  \\
          & = &  [{\xi}_q,{\chi}_q(\eta)]+i(\eta)D{\xi}_{q+1}-{\chi}_q(L({\xi}_1)\eta)\in J^0_q(T)
          \end{array}   \]
As $[R_q,R_q]\subset R_q$ and ${\chi}_q\in T^*\otimes R_q$, modifying ${\xi}_{q+1}$ by an element in $S_{q+1}T^*\otimes T$ if necessary, we get at once:  \\
\[   L({\xi}_{q+1}){\chi}_q\in T^*\otimes R^0_q \Leftrightarrow L({\xi}_q)c=0  \Leftrightarrow \exists {\xi}_{q+1}\in J_{q+1}(T), {\pi}^{q+1}_q({\xi}_{q+1})={\xi}_q\in R_q, D{\xi}_{q+1}\in T^*\otimes R_q   \]
Proceeding backwards in the definition of $D$ given at the beginning of section 4, we obtain therefore:\\
\[    L({\xi}_q)c=0  \Leftrightarrow \exists {\xi}_{q+1}\in R_{q+1}={\rho}_1(R_q), {\pi}^{q+1}_q({\xi}_{q+1})={\xi}_q\in R_q   \]
It follows that ${\pi}^{q+1}_q: R_{q+1}\longrightarrow R_q$ is surjective and we can choose a lift ${\chi}_{q+1}\in T^*\otimes R_{q+1}$ over ${\chi}_q\in T^*\otimes R_q$ in such a way that $\frac{1}{2}\{{\chi}_q,{\chi}_q\}\in {\wedge}^2T^*\otimes R^0_q\Rightarrow \frac{1}{2}\{c,c\}=0$ in $F_2$.  \\
\hspace*{12cm}   Q.E.D.   \\

\noindent
{\bf REMARK  5.44}: As already noticed, the IC $L({\xi}_q)c=0$ may be obtained in two steps. First, in order to obtain the surjectivity ${\pi}^{q+1}_q:R^0_{q+1}\rightarrow R^0_q$, we must have $L({\xi}^0_q)c=0, \forall {\xi}^0_q\in R^0_q$, that is $c$ {\it must be} $R^0_q$-{\it invariant}. Second, in order to obtain the surjectivity ${\pi}^{q+1}_0: R_{q+1}\rightarrow T$, we must have $L({\chi}_q(\xi))c=0, \forall \xi\in T$ for any $R_q$-connection ${\chi}_q$.  \\

The following definition only depends on $R^0_q$:  \\

\noindent
{\bf DEFINITION  5.45}: A {\it truncated Lie algebra} is a couple $(R^0_q,c)$ where $R^0_q\subset J^0_q(T)$ is such that $[R^0_q,R^0_q]\subset R^0_q$ and $c\in A(R^0_q)$ is $R^0_q$-invariant with $\frac{1}{2}\{c,c\}=0$ in $F_2$.   \\

\noindent
{\bf THEOREM   5.46}:  The infinitesimal deformation of $Cart (R_q)=\{c\in A(R^0_q){\mid} c\hspace{2mm} is \hspace{2mm} R_q-invariant,\\ \frac{1}{2}\{c,c\}=0\}$ is just $Z_1(R_q)$.   \\

\noindent
{\it Proof}: If $c_t=c+tC+ ...$ is a deformation of $c$ and ${\chi}_q(t)={\chi}_q+t{X}^0_q+ ... $is a representative of $c_t$, we must have $C\in F_1$ and ${X}^0_q\in T^*\otimes J^0_q(T)$. Denoting simply  by $\{c,C\}=0$ the linearization of $\frac{1}{2}\{c,c\}=0$, we must have $ \{{\chi}_q,{X}^0_q\}\in {\wedge}^2T^*\otimes R^0_q+\delta (T^*\otimes S_{q+1}T^*\otimes T)$ and get successively with different lifts: \\
\[\{{\chi}_{q+1}(\xi),X^0_{q+1}(\eta)\}+\{X^0_{q+1}(\xi),{\chi}_{q+1}(\eta)\}-X^0_q(\{{\chi}_1(\xi),{\chi}_1(\eta)\})-{\chi}_q(\{X^0_1(\xi),{\chi}_1(\eta)\}+\{{\chi}_1(\xi),X^0_1(\eta)\})   \]
   \[  \begin{array}{rcl}
   & = & L({\chi}_q(\xi))X^0_q(\eta)-L({\chi}_q(\eta))X^0_q(\xi)   \\

    &   & -i(\xi)DX^0_{q+1}(\eta)+i(\eta)DX^0_{q+1}(\xi) -X^0_q(\{  \})-{\chi}_q(   ) \\
    & = & i(\eta)L({\chi}_q(\xi))X^0_q-i(\xi)L({\chi}_q(\eta))X^0_q-DX^0_{q+1}(\xi,\eta)-{\chi}_q(  )
    \end{array}      \]
where we have used the relation: \\
\[ i(\eta)L({\chi}_q(\xi))X^0_q=L({\chi}_q(\xi))X^0_q(\eta)-X^0_q(L({\chi}_1(\xi)\eta)  \]
with $L({\chi}_1(\xi))\eta=[\xi,\eta]+i(\eta)D{\chi}_1(\xi)$ and the relation:  \\
\[   i(\eta)DX^0_{q+1}(\xi)-i(\xi)DX^0_{q+1}(\eta)+X^0_q([\xi,\eta])=DX^0_{q+1}(\xi,\eta)  \]
that is finally: \\
\[  \{{\chi}_q,X^0_q\}(\xi,\eta)+DX^0_{q+1}(\xi,\eta)=i(\eta){\nabla}_{\xi}X^0_q-i(\xi){\nabla}_{\eta}X^0_q-{\chi}_q(  )   \]
Now $L({\xi}_q)C=0, \forall {\xi}_q\in R_q=R^0_q\oplus {\chi}_q(T) \Rightarrow {\nabla}_{\xi}C=0$ in $F_1, \forall \xi\in T\Rightarrow {\nabla}_{\xi}X^0_q\in T^*\otimes R^0_q+\delta (S_{q+1}T^*\otimes T)$ and thus: $\{ {\chi}_q, X^0_q\}+DX^0_{q+1}\in {\wedge}^2T^*\otimes R_q$. Modifying if necessary the lift $X^0_{q+1}$, we need only $DX^0_{q+1}\in {\wedge}^2T^*\otimes R_q$, a result leading to $L({\xi}_q)C=0, \forall {\xi}_q\in R_q$, that is $C\in {\Upsilon}_1$ and ${\cal{D}}_2C=0$ in $F_2$. \\
\hspace*{12cm}   Q.E.D.   \\

As the $F_r$ only depend on $R^0_q$ both with the sub-bundles $E_r\subset F_r$, we finally obtain ([37], p 721): \\

\noindent
{\bf COROLLARY  5.47}: The deformation cohomology $H_r(R_q)$ only depends on the truncated Lie algebra $(R^0_q,c)$.  \\

It finally remains to study a last but delicate problem, namely to compare the deformation cohomology of $R_q$ to that of $R_{q+1}$. First of all, as $R_q$ is supposed to be involutive in order to construct the corresponding Janet sequence with Janet bundles $F_r$, then $R_{q+1}$ is of course involutive too and we may construct the corresponding Janet sequence with Janet bundles $F'_r$. The link between the two Janet sequences is described by the next theorem. \\

\noindent
{\bf THEOREM  5.48}: There is the following commutative diagram with exact columns and locally exact top row:  \\

\[      \begin{array}{rcccccccccl}
 & & & & & 0 & & 0 & & 0 &  \\
 & & & & & \downarrow & & \downarrow & & \downarrow &   \\
 & & & 0 & \rightarrow & K_0 & \rightarrow & K_1 & \rightarrow \hdots \rightarrow & K_n &  \rightarrow 0 \\
 & & & \downarrow & & \downarrow & & \downarrow &  &\downarrow &  \\
0\rightarrow & \Theta & \rightarrow & T & \stackrel{{\cal{D}}'}{\rightarrow} & F'_0 & \stackrel{{\cal{D}}'_1}{\rightarrow}& F'_1 & \stackrel{{\cal{D}}'_2}{\rightarrow} \hdots \stackrel{{\cal{D}}'_n}{\rightarrow} & F'_n & \rightarrow 0  \\
 & \parallel & & \parallel & &\hspace{5mm} \downarrow {\Psi}_0 & & \hspace{5mm} \downarrow {\Psi}_1 &  & \hspace{5mm} \downarrow {\Psi}_n & \\
0 \rightarrow & \Theta & \rightarrow & T & \stackrel{{\cal{D}}}{\rightarrow} & F_0 & \stackrel{{\cal{D}}_1}{\rightarrow} & F_1 & \stackrel{{\cal{D}}_2}{\rightarrow} \hdots \stackrel{{\cal{D}}_n}{\rightarrow} &F_n & \rightarrow 0  \\
& & & \downarrow & & \downarrow & & \downarrow & & \downarrow &    \\
& & & 0 & & 0 & & 0 & & 0 & 
\end{array}    \]

\noindent
{\it Proof}: We may use the following commutative and exact diagram in order to construct successively the epimorphisms ${\Psi}_r:F'_r \rightarrow F_r$ with kernels $K_r$ by using an induction on $r$ starting with $r=0$. \\

\[  \begin{array}{rcccccccl}
 & 0 & & 0 & & 0 & & 0 &  \\
 & \downarrow & & \downarrow & & \downarrow & & \downarrow &  \\
 0 \rightarrow & g_{q+r+1} & \rightarrow & S_{q+r+1}T^*\otimes T & \rightarrow & J_r(K_0) & \rightarrow \hdots \rightarrow & K_r &  \rightarrow 0  \\
 & \downarrow & & \downarrow & & \downarrow & & \downarrow &  \\
0 \rightarrow & R_{q+r+1} & \rightarrow & J_{q+r+1}(T) & \rightarrow & J_r(F'_0) & \rightarrow \hdots \rightarrow & F'_r &  \rightarrow 0  \\
  & \downarrow &  & \hspace{12mm} \downarrow {\pi}^{q+r+1}_{q+r} & & \hspace{11mm} \downarrow  J_r({\Psi}_0) & & \hspace{5mm} \downarrow {\Psi}_r &   \\
  0 \rightarrow & R_{q+r} & \rightarrow & J_{q+r}(T) & \rightarrow & J_r(F_0) & \rightarrow \hdots \rightarrow & F_r & \rightarrow 0  \\
  & \downarrow & & \downarrow & & \downarrow & & \downarrow &   \\
  & 0 & & 0 & & 0 & & 0 & 
\end{array}       \]

It follows that {\it the Janet sequence for} ${\cal{D}}'$ {\it projects onto the Janet sequence for} ${\cal{D}}$ and we just need to prove that the kernel sequence, made up by first order operators only, is exact. For this, introducing the short exact sequences : \\
\[      0 \rightarrow g_{q+r+1} \rightarrow S_{q+r+1}T^*\otimes T \rightarrow h_{r+1} \rightarrow 0   \]
with $h_{r+1}\subset S_{r+1}T^*\otimes F_0$ the $r$-prolongation of the symbol $h_1\subset T^*\otimes F_0$ of the system $B_1\subset J_1(F_0)$, image of the first prolongation $J_{q+1}(T) \rightarrow J_1(F_0)$ of the epimorphism $J_q(T) \rightarrow F_0$ used in order to define ${\cal{D}}$. Here $h_1$ is identified with the reciprocal image of the symbol of ${\cal{B}}_1\subset J_1({\cal{F}})$ by $\omega$. Using the Spencer operator $D$ and its various extensions, we obtain the following commutative diagram where it is known that the vertical $D$-sequences are locally exact ([36],[49]). \\

\[  \begin{array}{rcccccl}
& & & 0 & & 0 &  \\
& & & \downarrow & & \downarrow &  \\
 & 0 & \rightarrow & K_0 & \rightarrow \hdots \rightarrow & K_n & \rightarrow 0 \\
 & \downarrow & & \hspace{7mm} \downarrow j_{r+n} & & \hspace{3mm} \downarrow j_r &   \\
0 \rightarrow & h_{r+n+1} & \rightarrow & J_{r+n}(K_0) & \rightarrow \hdots \rightarrow  & J_r(K_n) & \rightarrow 0  \\
 & \hspace{5mm} \downarrow -\delta &  & \hspace{4mm} \downarrow D &  & \hspace{4mm} \downarrow D &   \\
 0 \rightarrow & T^*\otimes h_{r+n} & \rightarrow & T^*\otimes J_{r+n-1}(K_0) & \rightarrow \hdots \rightarrow & T^*\otimes J_{r-1}(K_n) & \rightarrow 0 \\
 & \hspace{4mm} \downarrow -\delta & & \hspace{4mm} \downarrow D &  & \hspace{4mm} \downarrow D &   \\
 & \vdots & & \vdots & & \vdots &   \\
 & \hspace{4mm} \downarrow -\delta & & \hspace{4mm} \downarrow D & & & \\
0 \rightarrow  & {\wedge}^nT^*\otimes h_{r+1} & \rightarrow & {\wedge}^nT^*\otimes J_r(K_0) & \rightarrow \hdots \hspace{4mm}&  &   \\
 & \downarrow & & \downarrow &  &  &  \\
 & 0 & & 0 & & &  
 \end{array}      \]

Now the left vertical sequence is exact because, applying the $\delta$-sequence to the previous short exact sequence, then $h_1$ is involutive whenever $g_q$ is involutive. The exactness of the top row finally follows from a diagonal chase in the last diagram because the other rows are exact by induction according to the previous diagram. \\

\hspace*{12cm}   Q.E.D.   \\

In order to use the last theorem in a natural way, we shall, for simplicity, restrict to the use of the classical Lie derivative by using the formulas:  \\
\[  {\cal{L}}(\xi){\eta}_q=[j_q(\xi),{\eta}_q]       \hspace{2cm}      \forall \xi\in T, \forall {\eta}_q\in J_q(T)   \]
\[   {\cal{L}}(\xi)\{{\eta}_q,{\zeta}_q\}=\{{\cal{L}}(\xi){\eta}_q,{\zeta}_q\}+\{{\eta}_q,{\cal{L}}(\xi){\zeta}_q\}  \]
\[  {\cal{L}}(\xi)[{\eta}_q,{\zeta}_q]=[{\cal{L}}(\xi){\eta}_q,{\zeta}_q]+[{\eta}_q,{\cal{L}}(\xi){\zeta}_q]    \]
\[      D {\cal{L}}(\xi)={\cal{L}}(\xi)D    \]
which are direct consequences of the fact that the algebraic bracket, the differential bracket and $D$ only contain natural operations. Accordingly, when restricting to an involutive system $R_q\subset J_q(T)$ of infinitesimal Lie equations, we have to preserve the section $\omega$ of the natural bundle ${\cal{F}}$ in order to construct the Janet sequence and we obtain at once:  \\

\noindent
{\bf PROPOSITION  5.49}: We have the formulas:  \\
${\cal{L}}(\xi){\cal{D}}={\cal{D}}{\cal{L}}(\xi) , \hspace{1cm}   {\cal{L}}(\xi){\cal{D}}_r={\cal{D}}_r{\cal{L}}(\xi), \hspace{1cm}  \forall \xi \in \Theta, \forall r=0,..., n  $. \\

\noindent
{\bf EXAMPLE  5.50}: Whenever ${\cal{D}}$ is a Lie operator already defined by the Lie derivative of $\omega$ with respect to a vector field, we get:  \\
\[ ({\cal{L}}(\xi){\cal{D}}-{\cal{D}}{\cal{L}}(\xi))\eta={\cal{L}}(\xi){\cal{L}}(\eta)\omega-{\cal{D}}[\xi,\eta]=({\cal{L}}(\xi){\cal{L}}(\eta)-{\cal{L}}([\xi,\eta])) \omega={\cal{L}}(\eta){\cal{L}}(\xi)\omega={\cal{L}}(\eta){\cal{D}} \xi =0   \]
Such a result may be applied at once to all the known structures. The case of a symplectic structure with $n=2p, \omega \in {\wedge}^2T^*, det(\omega)\neq 0, d\omega =0 $ is particularly simple because a part of the corresponding Janet sequence is made by a part of the Poincar\'{e} sequence and it is well known that the Lie derivative ${\cal{L}}(\xi)$ commutes with the exterior derivative $d$ (See [37], p 682 for more details). \\

The two following propositions will be obtained from the diagram of the last theorem by means of unusual chases in the tridimensional diagram obtained by applying ${\cal{L}}(\xi)$ for $\xi \in \Theta$ to the diagram of the last theorem. However, before providing them, it is essential to notice that a short exact sequence $0 \rightarrow F \rightarrow F' \rightarrow F'' \rightarrow 0 $ of bundles associated with $R_q$ may only provide in general the exact sequence $0 \rightarrow \Upsilon \rightarrow {\Upsilon}' \rightarrow {\Upsilon}'' $ by letting ${\cal{L}}(\xi)$ acting on each bundle whenever $\xi\in \Theta$, unless one can split this sequence by using a natural map from $F''$ to $F'$. For example, in the case of the short exact sequence $0 \rightarrow J^0_q(T) \rightarrow J_q(T) \rightarrow T \rightarrow 0$ one can use $j_q:T\rightarrow J_q(T)$. We also recall that an isomorphism $Z'/B'=H'\simeq H=Z/B$ induced by maps $Z' \rightarrow Z$ and $B' \rightarrow B$ does not necessarily imply any property of these maps which may be neither monomorphisms nor epimorphisms.   \\

\noindent
{\bf PROPOSITION  5.51}: There is an isomorphism $H'_0=H_0(R_{q+1}) \simeq H_0(R_q)=H_0$.  \\

\noindent
{\it Proof}: We shall provide two different proofs:  \\
1)  First of all, this result is a direct consequence of the short exact sequence provided by Proposition 5.32.  \\
2)  Let us prove that there is a monomorphism $0 \rightarrow H'_0 \rightarrow H_0$. First, if $a' \in {\Upsilon}'_0\subset F'_0$ is killed by ${\cal{D}}'_1$ and such that ${\Psi}_0(a')=a={\cal{D}}\eta$ for a certain $\eta\in \Upsilon (T)\subset T$, then $a={\Psi}_0(a')={\Psi}_0({\cal{D}}' \eta)$ and thus ${\Psi}_0(a'-{\cal{D}}' \eta)=0$, that is $a'-{\cal{D}}' \eta=b' \in K_0$ with ${\cal{D}}'_1 b'={\cal{D}}'_1 a'-{\cal{D}}'_1 {\cal{D}}' \eta=0$, that is $b'=0 \Rightarrow  a' = {\cal{D}}' \eta$.  \\
Let us then prove that there is an epimorphism $H'_0 \rightarrow H_0 \rightarrow 0 $. For this, if $a\in {\Upsilon}_0\subset F_0$ is killed by ${\cal{D}}_1$, we may find $a' \in F'_0$ such that $a={\Psi}_0(a' )$ and ${\cal{L}}(\xi) a' \in K_0$. It follows that $c'={\cal{D}}'_1 a' \in K_1$ with ${\cal{D}}'_2 c'=0$, that is $c'={\cal{D}}'_1b' $ for a certain $b' \in K_0$. Accordingly, modifying $a' $ if necessary by replacing it by 
$a' -b' $, we may suppose that ${\cal{D}}'_1a'=0$ with ${\Psi}_0(a' )=a$ and we just need to prove that ${\cal{L}}(\xi)a'=0$. Indeed, we have ${\Psi}_0({\cal{L}}(\xi)a')={\cal{L}}(\xi)a=0 \Rightarrow {\cal{L}}(\xi)a' \in K_0$ and thus ${\cal{D}}'_1{\cal{L}}(\xi)a'={\cal{L}}(\xi){\cal{D}}'_1a'=0 \Rightarrow {\cal{L}}(\xi)a'=0 $ because the restriction of ${\cal{D}}'_1$ to $K_0$ is a monomorphism.  \\

\hspace*{12cm}   Q.E.D.   \\

\noindent
{\bf PROPOSOTION  5.52}: If $R_q\subset J_q(T)$ is involutive, the normalizer $N(\Theta)$ of $\Theta$ in $T$ is defined by the involutive system ${\tilde{R}}_{q+1}\subset J_{q+1}(T)$ with involutive symbol ${\tilde{g}}_{q+1}=g_{q+1}\subset S_{q+1}T^*\otimes T$, defined by the purely algebraic condition $\{R_{q+1},{\tilde{R}}_{q+1}\}\subset R_q$. Moreover, $Z_0(R_{q+1})\simeq Z_0(R_q)\simeq {\tilde{R}}_{q+1}/R_{q+1}$. \\
 
\noindent
{\it Proof}: With $Z_0=Z_0(R_q)$ and $Z'_0=Z_0(R_{q+1})$ as before, let us prove that there is an isomorphism $0 \rightarrow Z'_0 \rightarrow Z_0 \rightarrow 0 $. For this, if $a'\in Z'_0$ is such that ${\Psi}_0(a')=0$, then $a'\in K_0$ with ${\cal{D}}'_1a'=0 \Rightarrow a'=0$ because the restriction of ${\cal{D}}'_1$ to $K_0$ is injective and we get a monomorphism $0 \rightarrow Z'_0 \rightarrow Z_0 \rightarrow 0$.  \\
Similarly, if $a\in Z_0$, using the previous proposition, we can find $a'\in F'_0$ with ${\cal{D}}'_1a'=0$ and ${\Psi}_0(a')=a$. It follows that ${\cal{L}}(\xi)a'=d'\in K_0$ with ${\cal{D}}'_1d'={\cal{D}}'_1{\cal{L}}(\xi)a'={\cal{L}}(\xi){\cal{D}}'_1a'=0 \Rightarrow d'=0 \Rightarrow {\cal{L}}(\xi)a'=0 \Rightarrow a'\in Z'_0$ and we get an epimorphism $Z'_0 \rightarrow Z_0 \rightarrow 0$, that is $Z_ˆ(R_{q+1})\simeq Z_0(R_q)$.  \\
Now we have ${\cal{L}}({\xi}_{q+1}){\eta}_q=\{{\xi}_{q+1},{\eta}_{q+1}\}+i(\xi)D{\eta}_{q+1}\in R_q$ for {\it any lift} ${\eta}_{q+1}\in J_{q+1}(T)$ of ${\eta}_q\in J_q(T)$. However, cocycles in $Z_0$ are also killed by ${\cal{D}}_1$ and the representative ${\eta}_q\in J_q(T)$ must be such that $D{\eta}_{q+1}\in T^*\otimes R_q+\delta (S_{q+1}T^*\otimes T)$. Therefore , modifying ${\eta}_{q+1}$ if necessary without changing ${\eta}_q$, we may suppose that $D{\eta}_{q+1}\in T^*\otimes R_q$, a result leading to the condition $\{{\xi}_{q+1},{\eta}_{q+1}\}\in R_q \Rightarrow {\eta}_{q+1}\in {\tilde{R}}_{q+1}$. As we already know that ${\tilde{g}}_{q+1}=g_{q+1}$, introducing the projection ${\tilde{R}}_q\subset J_q(T)$ of ${\tilde{R}}_{q+1}\subset J_{q+1}(T)$, we obtain $Z_0(R_q)\simeq {\tilde{R}}_q/R_q\simeq {\tilde{R}}_{q+1}/R_{q+1}$ and similarly $Z_0(R_{q+1})\simeq {\tilde{R}}_{q+2}/{\tilde{R}}_{q+2}$ in the following commutative and exact diagram:  \\
\[   \begin{array}{rcccccl}
  & 0 & & 0 & &  &   \\
  & \downarrow &  & \downarrow &  &  &    \\
0 \rightarrow & g_{q+2}  & \rightarrow & {\tilde{g}}_{q+2} & \rightarrow & 0 &   \\
  & \downarrow  &  & \downarrow & & \downarrow &  \\
0 \rightarrow & R_{q+2} & \rightarrow & {\tilde{R}}_{q+2} & \rightarrow & Z_0(R_{q+1})  &  \rightarrow 0  \\
  &    \downarrow &  &  \downarrow  &  & \downarrow  &   \\
 0 \rightarrow  & R_{q+1} & \rightarrow  & {\tilde{R}}_{q+1} & \rightarrow  & Z_0(R_q) & \rightarrow  0  \\
   &  \downarrow  &  & & &  \downarrow &    \\
   &  0  &  &  &  &  0  &
      \end{array}   \]
Chasing in this diagram or counting the dimensions in order to get:  \\
\[  dim ({\tilde{R}}_{q+2})- dim ({\tilde{R}}_{q+1}) = dim (R_{q+2}) - dim (R_{q+1}) = dim (g_{q+2})= dim ({\tilde{g}}_{q+2})   \]
it follows that ${\pi}^{q+2}_{q+1}:{\tilde{R}}_{q+2} \rightarrow {\tilde{R}}_{q+1}$is an epimorphism.  \\
Next, applying the formula of Proposition 5.6, we get:   \\
\[ \{{\xi}_{q+1},i(\zeta)D{\eta}_{q+2}\}=i(\zeta)D\{{\xi}_{q+2},{\eta}_{q+2}\}-\{i(\zeta)D{\xi}_{q+2},{\eta}_{q+1}\}\in R_q  \Rightarrow D{\eta}_{q+2}\in T^*\otimes {\tilde{R}}_{q+1}  \]
It follows that ${\tilde{R}}_{q+2}={\rho}_1({\tilde{R}}_{q+1})$ because both systems project onto ${\tilde{R}}_{q+1}$ and have the same symbol ${\tilde{g}}_{q+2}=g_{q+2}={\rho}_1(g_{q+1})={\rho}_1({\tilde{g}}_{q+1})$. \\
Finally, as ${\tilde{g}}_{q+1}=g_{q+1}$ is an involutive symbol when $g_q$ is involutive and ${\pi}^{q+2}_{q+1}:{\tilde{R}}_{q+2} \rightarrow {\tilde{R}}_{q+1}$ is surjective, then ${\tilde{R}}_{q+1}$ is involutive because of the Janet/Goldschmidt/Spencer criterion of formal integrability ([23],[36],[39],[49]). \\

\hspace*{12cm}     Q.E.D.   \\

\noindent
{\bf PROPOSITION  5.53}: There is a monomorphism $0 \rightarrow H'_1\rightarrow H_1$. Accordingly, $H_1(R_{q+1})=0$ whenever $H_1(R_q)=0$ and thus $R_{q+1}$ is rigid whenever $R_q$ is rigid and we may say that $\Theta$ is rigid.  \\

\noindent
{\it Proof}: If $a' \in {\Upsilon}'_1 \subset F'_1$ is killed by ${\cal{D}}'_2$ and such that ${\Psi}_1(a')=a={\cal{D}}_1c$ for a certain $c\in {\Upsilon}_0\subset F_0$, we may find $c'\in F'_0$ such that ${\Psi}_0(c')=c$ and thus $a={\Psi}_1(a' )={\cal{D}}_1{\Psi}_0c'={\Psi}_1{\cal{D}}'_1c' \Rightarrow {\Psi}_1 (a' - {\cal{D}}'_1 c' )=0 \Rightarrow a' - {\cal{D}}'_1 c' = b' \in K_1$. As before, we get ${\cal{D}}'_2 b' = {\cal{D}}'_2 a' - {\cal{D}}'_2 {\cal{D}}'_1 c'=0 \Rightarrow \exists d' \in K_0$ with ${\cal{D}}'_1 d' = b' $ and thus  $a' = {\cal{DD}}'_1 (c' + d' )$. Hence, modifying $c'$ if necessary, we may find $c' \in F'_0$ such that ${\Psi}_0 (c' )=c$ with $a' ={\cal{D}}'_1 c'$ and it only remains to prove that $c' \in {\Upsilon}'_0$. However, we have ${\cal{L}}(\xi)c' \in K_0$ because ${\cal{L}}(\xi)c=0$ and ${\cal{D}}'_1{\cal{L}}(\xi)c' = {\cal{L}}(\xi){\cal{D}}'_1c' = {\cal{L}}(\xi)a' =0$. As the restriction of ${\cal{D}}'_1$ to $K_0$ is injective, we finally obtain ${\cal{L}}(\xi)c' =0$ as we wanted.  \\

\hspace*{12cm}  Q.E.D.   \\

Despite many attempts we have not been able to find additional results, in particular to compare $H_2(R_q)$ and $H_2(R_{q+1})$. However, integrating the Vessiot structure equations with structure constants $c_t$ instead of $c$, we may exhibit a deformation ${\omega}_t$ of $\omega$ as a new section of ${\cal{F}}$. As ${\cal{D}}$ depends on $j_1(\omega)$, the ${\cal{D}}_r$ will also only depend on $\omega$ and various jets. Replacing $\omega$ by ${\omega}_t$, we shall obtain operators ${\cal{D}}(t)$ and ${\cal{D}}_r(t)$ for $0\leq r \leq n$. Therefore we may apply all the techniques and results of section 
2 with only slight changes.  \\
\vspace*{5mm}  \\
\noindent
{\bf 6  EXPLICIT COMPUTATIONS}:  \\

\noindent
{\bf EXAMPLE  6.1}: (1.1 {\it revisited}) This is a good example for understanding that the Janet bundles are only defined up to an isomorphism, contrary to the natural bundles. In this example, we have of course the basic natural bundle ${\cal{F}}=T^*\otimes V=T^*{\times}_X ... {\times}_XT^*$ ($n$-times) with $V={\mathbb{R}}^n$. This is a vector bundle (even a tensor bundle) and we may therefore identity ${\cal{F}}$ with ${\cal{F}}_0=V({\cal{F}})$ and $F_0$. However, in this case, $g_1=R^0_1=0$ and we get $F_0=J^0_1(T)/R^0_1=T^*\otimes T$. Similarly, we should get $F_r={\wedge}^rT^*\otimes T^*\otimes T/\delta({\wedge}^{r-1}T^*\otimes S_2T^*\otimes T)={\wedge}^{r+1}T^*\otimes T$ because we have the exact sequence :\\
\[  {\wedge}^{r-1}T^*\otimes S_2T^*\stackrel{\delta}{\longrightarrow}{\wedge}^rT^*\otimes T^*\stackrel{\delta}{\longrightarrow}{\wedge}^{r+1}T^* \longrightarrow 0  \]
Using geometric objects, we should obtain: \\
\[ ({\cal{D}}\xi)^{\tau}_i\equiv {\omega}^{\tau}_r(x){\partial}_i{\xi}^r+{\xi}^r{\partial}_r{\omega}^{\tau}_i(x)={\Omega}^{\tau}_i \Rightarrow ({\cal{D}}_1\Omega)^{\tau}_{ij}\equiv {\partial}_i{\Omega}^{\tau}_j-{\partial}_j{\Omega}^{\tau}_i-c^{\tau}_{\rho\sigma}({\omega}^{\rho}_i{\Omega}^{\sigma}_j+{\omega}^{\sigma}_j{\Omega}^{\rho}_i)=0 \]
 and thus $F_1={\wedge}^2T^*\otimes V$. On the contrary, using the solved form $\xi\rightarrow (\nabla\xi)^k_i\equiv {\partial}_i{\xi}^k+{\xi}^r{\alpha}^k_{\tau}(x){\partial}_r{\omega}^{\tau}_i(x)$ providing at once a unique $R_1$-connexion, a covariant derivative $\nabla$ and thus a well defined $\nabla$-sequence leading to $F_r={\wedge}^{r+1}T^*\otimes T $ and thus $F_1={\wedge}^2T^*\otimes T$. In fact, {\it the explanation is a confusion between the Janet sequence and the Spencer sequence} because we have indeed $C_0=R_1=T\Rightarrow C_{r+1}\simeq F_r\simeq {\wedge}^{r+1}T^*\otimes T$.  \\
The sections of ${\wedge}^rT^*\otimes V$ invariant by $R_1$ are of the form $A^{\tau}_{{\sigma}_1... {\sigma}_r}{\omega}^{{\sigma}_1}\wedge ...\wedge{\omega}^{{\sigma}_r}$ while the sections of ${\wedge}^rT^*\otimes T$ invariant by $R_1$ are of the form ${\alpha}^k_{\tau}A^{\tau}_{{\sigma}_1... {\sigma}_r}{\omega}^{{\sigma}_1}\wedge ...\wedge{\omega}^{{\sigma}_r}$ and are thus different though depending on the same {\it constants} $A^{\tau}_{{\sigma}_1... {\sigma}_r}$, providing therefore a unique cohomology. From the Maurer-Cartan equations $d{\omega}^{\tau}=-c^{\tau}_{\rho\sigma}{\omega}^{\rho}\wedge{\omega}^{\sigma}$ where we change the sign of the structure constants fo convenience, we obtain terms of the form $Ac$ which are {\it exactly} the first terms in the definition of the Chevalley-Eilenberg operator. Now, looking at the terms containing $\alpha$, we obtain (care to the minus sign):   \\
\[  ({\partial}_i{\alpha}^k_{\tau}+{\alpha}^r_{\tau}{\alpha}^k_{\sigma}{\partial}_r{\omega}^{\sigma}_i)A^{\tau}dx^i={\alpha}^r_{\tau}{\alpha}^k_{\sigma}({\partial}_r{\omega}^{\sigma}_i-{\partial}_i{\omega}^{\sigma}_r)A^{\tau}dx^i=-{\alpha}^k_{\tau}c^{\tau}_{\rho\sigma}A^{\rho}{\omega}^{\sigma} \]
and recover the second terms of this operator.\\
{\it It is however essential to notice that the two approaches are totally different}. In particular, elements in $A(R^0_q)$ are $1$-forms with value in some vector bundles and {\it it is therefore a pure chance that they could become $2$-forms in this example} (See electromagnetism and gravitation in [43]).  \\
We may also compare the concept of "change of basis" of the Lie algebra with the "label transformations". Indeed, using the Lie form, we get (care again to the minus sign):  \\
\[  \begin{array}{rcl}
{\bar{\omega}}\sim \omega & \Leftrightarrow & {\bar{\alpha}}^k_{\tau}{\partial}_r{\bar{\omega}}^{\tau}_i={\alpha}^k_{\rho}{\partial}_r{\omega}^{\rho}_i=-{\omega}^{\rho}_i{\partial}_r{\alpha}^k_{\rho} \\
  &  \Leftrightarrow  & {\alpha}^i_{\rho}{\partial}_r{\bar{\omega}}^{\tau}_i+{\bar{\omega}}^{\tau}_k{\partial}_r{\alpha}^k_{\rho} =0 \\
    & \Leftrightarrow  &  {\partial}_r({\alpha}^i_{\rho}{\bar{\omega}}^{\tau}_i)=0  \\
     & \Leftrightarrow & {\bar{\omega}}^{\tau}_i=a^{\tau}_{\sigma}{\omega}^{\sigma}_i, a=cst  
     \end{array}     \]
Accordingly, the sections of ${\Upsilon}_0\subset T^*\otimes V$ are of the form $A^{\tau}_{\sigma}{\omega}^{\sigma}_i$ with $A=cst$ and it is rather extraordinary that two such different approaches can provide the same result.\\
Finally, we study $\Theta, Z(\Theta), C(\Theta)$ and $N(\Theta)$ in this framework. First of all, $\Theta$ is described by a linear combination with constant coefficients of the infinitesimal generators $\{{\theta}_{\tau}={\theta}^i_{\tau}{\partial}_i\}$. Then we must have $-{\xi}^k_i{\eta}^i+{\xi}^r{\partial}_r{\eta}^k=0$, whenever $ {\xi}^k_i+{\xi}^r{\alpha}^k_{\tau}{\partial}_r{\omega}^{\tau}_i=0$, that is ${\alpha}^k_{\tau}{\partial}_r{\omega}^{\tau}_i{\eta}^i+{\partial}_r{\eta}^k=0$ or ${\partial}_r({\omega}^{\tau}_i{\eta}^i)=0$ and thus ${\eta}^i={\lambda}^{\tau}{\alpha}^i_{\tau}$ with $\lambda=cst$ ({\it reciprocal distribution}). We already know that $[{\alpha}_{\rho},{\alpha}_{\sigma}]= c^{\tau}_{\rho\sigma}{\alpha}_{\sigma}$ with our choice of sign. The isomorphism between $\Theta$ and $C(\Theta)$ is therefore a pure chance. As for $N(\Theta)$, it follows from Corollary 5.27 that the defining equations are $({\cal{L}}(\xi)\omega)^{\tau}_i\equiv {\omega}^{\tau}_r(x){\partial}_i{\xi}^r+{\xi}^r{\partial}_r{\omega}^{\tau}_i(x)=A^{\tau}_{\sigma}{\omega}^{\sigma}_i(x)$ where $A=cst$ is a derivation of the Lie algebra ${\cal{G}}$ with the same structure constants $c$. \\

\noindent
{\bf EXAMPLE  6.2}: (1.2 {\it revisited}) Let us study the equivalence relation ${\bar{\omega}}\sim \omega$ by asking first that ${\bar{R}}^0_1=R^0_1$. We must have ${\bar{\omega}}_{rj}{\xi}^r_i+{\bar{\omega}}_{ir}{\xi}^r_j=0\Leftrightarrow {\omega}_{rj}{\xi}^r_i+{\omega}_{ir}{\xi}^r_j=0$. Setting ${\xi}_{i,j}={\omega}_{rj}{\xi}^r_i$, it amounts to check that ${\bar{\omega}}_{rj}{\omega}^{rt}{\xi}_{i,t}+{\bar{\omega}}_{ir}{\omega}^{rt}{\xi}_{j,t}=0 \Leftrightarrow {\xi}_{i,j}+{\xi}_{j,i}=0$ and we must therefore have: \\
\[  {\bar{\omega}}_{rj}{\omega}^{rt}{\delta}^s_i+{\bar{\omega}}_{ir}{\omega} ^{rt}{\delta}^s_j-{\bar{\omega}}_{rj}{\omega}^{rs}{\delta}^t_i-{\bar{\omega}}_{ir}{\omega}^{rs}{\delta}^t_j=0  \]
Contracting in $s$ and $j$, we get $n{\bar{\omega}}_{rj}{\omega}^{rt}-{\bar{\omega}}_{rs}{\omega}^{rs}{\delta}^t_j=0$, that is ${\bar{\omega}}_{rj}{\omega}^{rt}=a(x){\delta}^t_j\Rightarrow {\bar{\omega}}_{ij}=a(x){\omega}_{ij}$. \\
Substituting in ${\bar{R}}_1=R_1$, we finally get $\bar{\omega}=a\omega, a=cst$ and the group of label transformations is just the multiplicative group.It then follows that ${\bar{\gamma}}=\gamma$ but such a property can be checked directly from the solved form of the second order equations ${\cal{L}}(\xi)\gamma=0$. It follows that ${\bar{\rho}}=\rho$ and thus $\bar{c}=(1/a)c$ and, in any case, $N(\Theta)$ is defined by the infinitesimal Lie equations ${\cal{L}}(\xi)\omega=A\omega$ with $cA=0$ as there is no Jacobi condition on the single structure constant $c$. Accordingly, $N(\Theta)=\Theta$ if $c\neq 0$ and $dim(N(\Theta)/\Theta)=1$ if $c=0$. This is the reason for which the Weyl group is the normalizer of the Poincar\'{e} group, obtained by adding the generator $x^i{\partial}_i$ of the dilatation on space-time. It is important to notice that the Galil\'{e}e group is of codimension $2$ in its normalizer obtained by dilatating {\it separately} space and time, another reason for which the Poincar\'{e} group cannot be obtained from the Galil\'{e}e group by a continuous deformation. As for $C(\Theta)$, using again $R^0_1$, we have to solve ${\eta}^i{\omega}^{kj}{\xi}_{i,j}=0\Leftrightarrow {\eta}^i{\xi}_{i,j}=0, \forall {\xi}_{i,j}+{\xi}_{j,i}=0$ and thus $C(\Theta)=0\Rightarrow Z(\Theta)=0$.  \\

\noindent
{\bf EXAMPLE  6.3}: (1.3 {\it revisited})  This example is by far the most difficult to treat in dimension $n=2p+1$ ([37], p 684). In the present case when $n=3$, this is one of the best examples where the Lie equations obtained by eliminating $\rho$, namely ${\xi}^1_3-x^3{\xi}^2_3=0, {\xi}^1_2-x^3{\xi}^2_2+x^3
{\xi}^1_1-(x^3)^2{\xi}^2_1-{\xi}^3=0$ could be used without even knowing about the underlying geometric object. However, {\it in this case} $R_1$ {\it is not involutive} and we must start afresh with the involutive system $R^{(1)}_1\subset R_1$ for constructing the canonical Janet sequence. We first notice that there is only one CC for the new involutive system $R_1\subset J_1(T)$ of infinitesimal Lie equations and thus surely no Jacobi condition for the only structure constant $c$. It thus remains to study the inclusions $Z(\Theta)\subseteq \Theta\subseteq N(\Theta)$. Using $R^0_1\subset J^0_1(T)$ for $C(\Theta)$, we have to solve ${\eta}^s{\xi}^k_s=0$ whenever ${\omega}_r{\xi}^r_i-\frac{1}{2}{\omega}_i{\xi}^r_r=0$ and thus ${\eta}^s{\omega}_s{\xi}^r_r=0\Rightarrow {\eta}^s{\omega}_s=0$. As $\omega\neq 0$, changing coordinates if necessary, we may suppose that ${\omega}_1\neq 0$. It follows from the involutive assumption that at least one jet coordinate of each class is parametric and in fact, changing coordinates if necessary, we have in all the examples presented $pri=\{{\xi}^3_3, {\xi}^1_3, {\xi}^1_2\}\Rightarrow par=\{{\xi}^1_1,{\xi}^2_1, {\xi}^3_1, {\xi}^2_2, {\xi}^3_2, {\xi}^2_3\}$. Hence, choosing successively ${\xi}^2_1, {\xi}^2_2, {\xi}^2_3$ as unique non-zero parametric jet, we obtain ${\eta}^1=0, {\eta}^2=0, {\eta}^3=0\Rightarrow \eta=0$ and thus $C(\Theta)=0\Rightarrow Z(\Theta)=0$. As for $N(\Theta)$, we have to study the equivalence 
${\bar{\omega}}\sim \omega$, that is to study when we have ${\bar{\omega}}_r{\xi}^r_i-\frac{1}{2}{\bar{\omega}}_i{\xi}^r_r=0\Leftrightarrow {\omega}_r{\xi}^r_i-\frac{1}{2}{\omega}_i{\xi}^r_r$. Though it looks like to be a simple algebraic problem, one needs an explicit computation or computer algebra and we prefer to use another more powerful technique ([37], p 688). Introducing the completely skewsymmetrical symbol $\epsilon=({\epsilon}^{i_1i_2i_3})$ where ${\epsilon}^{i_1i_2i_3}=1$ if $(i_1i_2i_3)$ is an even permutation of $(123)$ or $-1$ if it is an odd permutation and $0$ otherwise, let us introduce the skewsymmetrical $2$-contravariant density ${\omega}^{ij}={\epsilon}^{ijk}{\omega}_k$. Then one can rewrite the lie equations as :  \\
\[   R_1   \hspace{5cm} -{\omega}^{rj}(x){\xi}^i_r-{\omega}^{ir}(x){\xi}^j_r-\frac{1}{2}{\omega}^{ij}(x){\xi}^r_r+{\xi}^r{\partial}_r{\omega}^{ij}(x)=0   \]
and we may exhibit a section ${\xi}^i_r={\omega}^{is}A_{rs}$ with $A_{rs}=A_{sr}$ and thus ${\xi}^r_r=0$. It is important to notice that $det(\omega)=0$ when $n=2p+1$, contrary to the Riemann or symplectic case and $\omega$ cannot therefore be used in order to raise or lower indices. As we must have ${\bar{R}}^0_1=R^0_1$, the same section must satisfy $({\bar{\omega}}^{rj}{\omega}^{is}+{\bar{\omega}}^{ir}{\omega}^{js})A_{rs}=0, \forall A_{rs}=A_{sr}$, and we must have $({\bar{\omega}}^{rj}{\omega}^{is}+{\bar{\omega}}^{ir}{\omega}^{js})+({\bar{\omega}}^{sj}{\omega}^{ir}+{\bar{\omega}}^{is}{\omega}^{jr})=0$. Setting $s=j$, we get ${\bar{\omega}}^{rj}{\omega}^{ij}={\bar{\omega}}^{ij}{\omega}^{rj}\Rightarrow {\bar{\omega}}^{ij}(x)=a(x){\omega}^{ij}(x)$. Substituting and substracting, we get ${\omega}^{ij}(x){\xi}^r{\partial}_ra(x)=0\Rightarrow a(x)=a=cst$ or similarly ${\bar{\omega}}_i(x)=a(x){\omega}_i(x)\Rightarrow {\omega}_i(x){\xi}^r{\partial}_ra(x)=0\Rightarrow a(x)=a=cst$ because $\omega\neq 0$ and one of the components at least must be nonzero. It follows at once that one has $N(\Gamma)=\{ f\in aut(X)\mid j_1(f)^{-1}(\omega)=a\omega, a^2c=c\}$. Accordingly, $N(\Theta)=\{ \xi\in T\mid {\cal{L}}(\xi)\omega=A\omega, Ac=0\}$ and $\Theta$ is of codimension $1$ in its normalizer if $c=0$ or $N(\Theta)=\Theta$ if $c\neq 0$. For example, in the case of a contact structure with $c=1$, we have $N(\Theta)=\Theta$ but, when $\omega=(1,0,0)\Rightarrow c=0$, we have to eliminate the constant $A$ among the equations ${\partial}_3{\xi}^3+{\partial}_2{xi}^2-{\partial}_1{\xi}^1=-2A, {\partial}_3{\xi}^1=0, {\partial}_2{\xi}^1=0$ and we may add the infinitesimal generator $x^i{\partial}_i$ of a dilatation providing $A=-\frac{1}{2}$. \\

\noindent
{\bf EXAMPLE  6.4}: (1.4 {\it revisited}) This is by far the most interesting example. First of all, contrary to all the previous examples, we have at once ${\cal{F}}={\cal{F}}_0=T^*{\times}_X{\wedge}^2T^*, {\cal{F}}_1={\wedge}^2T^*{\times}_X{\wedge}^3T^*, {\cal{F}}_2={\wedge}^3T^*$ and we may identify $F_r$ with ${\cal{F}}_r$ for $r=0,1,2,3$ in order to obtain the Janet sequence: \\
\[   0 \longrightarrow \Theta \longrightarrow T \stackrel{{\cal{D}}}{\longrightarrow} F_0\stackrel{{\cal{D}}_1}{\longrightarrow}F_1\stackrel{{\cal{D}}_2}{\longrightarrow} F_2\longrightarrow 0 \]
with $dim(T)=3, dim(F_0)=6, dim(F_1)=4, dim(F_2)=1$ and $3-6+4-1=0$. Using local coordinates $(u^1,u^2,u^3)$ for $T^*$ and $(u^4,u^5,u^6)$ for ${\wedge}^2T^*$ in $F_0$, we obtain $W=\{W_1=u^1\frac{\partial}{\partial u^1}+u^2\frac{\partial}{\partial u^2}+u^3\frac{\partial}{\partial u^3}, W_2=u^4\frac{\partial}{\partial u^4}+u^5\frac{\partial}{\partial u^5}+u^6\frac{\partial}{\partial u^6}\}$. With $\omega=(\alpha,\beta)$, we have ${\cal{D}}\xi=({\cal{L}}(\xi)\alpha,{\cal{L}}(\xi)\beta)$ and recall the IC made by the Vessiot structure equations $d\alpha=c' \beta, d\beta=c'' \alpha\wedge\beta$ with $\gamma=\alpha\wedge\beta\neq 0$. As before, we may easily find the label transformations $(\bar{\alpha}=a \alpha, \bar{\beta}=b \beta)\Rightarrow ({\bar{c}}'=\frac{a}{b}c', {\bar{c}}''=\frac{1}{a}c'')$ and obtain ${\Upsilon}_0=(A\alpha,B\beta)\subset F_0$ with $A,B=cst$. Similarly, we get ${\Upsilon}_1=(C'\beta,C''\gamma)\subset F_1$ with $C',C''=cst$ and the induced map ${\cal{D}}_1:{\Upsilon}_0\rightarrow {\Upsilon}_1: (A,B)\rightarrow (c'(A-B)=C',-c''A=C'')$ in a coherent way with the linearization of the Vessiot structure equations:  \\
\[  Ad\alpha-c'B\beta=c'(A-B)\beta, Bd\beta-c''A\alpha\wedge\beta-c'' B\alpha\wedge\beta=-c''A\gamma \]
We obtain therefore $N(\Theta)=\{\xi\in T\mid {\cal{L}}(\xi)\alpha=A\alpha, {\cal{L}}(\xi)\beta=B\beta, c'(A-B)=0,c''A=0 \}$ and thus:  $c=(0,0)\Rightarrow dim(N(\Theta)/\Theta)=2, c=(1,0)\Rightarrow A=B\Rightarrow dim(N(\Theta)/\Theta)=1, c=(0,1)\Rightarrow A=0\Rightarrow dim(N(\Theta)/\Theta)=1$. The only Jacobi condition $c'c''=0$ provides $c''C'+c'C''=0$ in any case. For example, in the unimodular contact case, we must add $x^1{\partial}_1+x^2{\partial}_2$ to $\Theta$ in order to obtain $N(\Theta)$. \\
As for the centralizer $C(\Theta)$, we must look first for $\eta\in T$ such that ${\eta}^r{\xi}^k_r=0, \forall {\xi}^0_1\in R^0_1$. For this, multiplying $ {\beta}_{rj}{\xi}^r_i+{\beta}_{ir}{\xi}^r_j=0$ by ${\eta}_i$ and contracting on $i$, we get $({\eta}^i{\beta}_{ir}){\xi}^r_j=0$ whenever ${\alpha}_r{\xi}^r_i=0$. Accordingly, we obtain ${\eta}^i{\beta}_{ir}=L(x){\alpha}_r$ as $1$-forms. Now, as $n=3$, we may introduce the pseudo-vector $({\beta}^1={\beta}_{23}, {\beta}^2={\beta}_{31}, {\beta}^3={\beta}_{12})$ transforming like a vector up to a division by the Jacobian $\Delta$ and thus, using the volume form $\gamma=\alpha\wedge\beta$, it follows that ${\tilde{\beta}}=({\beta}^1/\gamma, {\beta}^2/\gamma,{\beta}^3/\gamma)$ is a true vector field. As ${\beta}_{ir}{\tilde{\beta}}^r=0$ because $n=3$, we obtain by contraction $L(x){\alpha}_r{\tilde{\beta}}^r=L(x)(\alpha\wedge\beta/\gamma)=L(x)=0\Rightarrow {\eta}^i{\beta}_{ir}=0\rightarrow {\eta}^k=K(x){\tilde{\beta}}^k$. But $\alpha, \beta, \gamma$ and thus ${\tilde{\beta}}$ are invariant by any $\xi \in \Theta$ and thus $K(x)=K=cst$, that is $C(\Theta)=\{\theta\in T{\mid}\eta=K{\tilde{\beta}}, K=cst\}$. In all the three special sections considered, we have simply $C(\Theta)=\{(K, 0, 0)\in T{\mid}K=cst\}$.\\

We obtain finally the folowing recapitulating picture:  \\

\[    \begin{array}{rcccccccl}
& (K) & \longrightarrow & (A, B) & \longrightarrow & (C' , C'' ) & \longrightarrow & (D )& \longrightarrow 0  \\
   &(K{\tilde{\beta}}) &\longrightarrow &  (A\alpha , B\beta)& \longrightarrow & (C' \beta , C'' \gamma) & \longrightarrow & (D\gamma ) &   \longrightarrow 0  \\
0 \longrightarrow \Theta \longrightarrow & T & \stackrel{{\cal{D}}}{\longrightarrow } & T^*\times_X {\wedge}^2T^* & \stackrel{{\cal{D}}_1}{\longrightarrow } & {\wedge}^2T^*\times_X{\wedge} ^3T^* & \stackrel{{\cal{D}}_2}{\longrightarrow} & {\wedge}^3T^*&  \longrightarrow  0 \\
  & \xi & \longrightarrow & (\Phi,\Psi) & \longrightarrow & (U, V) &\longrightarrow & (W) & \longrightarrow 0
\end {array}   \]
 
The {\it purely differential lower part} is describing the operators involved in the Janet sequence, namely and successively:  \\
\[  \begin{array}{lc}
{\cal{D}}     \hspace{2cm} &\xi \longrightarrow ({\cal{L}}(\xi)\alpha = \Phi, {\cal{L}}(\xi)\beta =\Psi ) \\
{\cal{D}}_1 \hspace{2cm}  & (\Phi , \Psi ) \longrightarrow (d\Phi - c' \Psi = U, d\Psi - c'' (\beta \wedge \Phi + \alpha \wedge \Psi) =V) \\
{\cal{D}}_2 \hspace{2cm} &   (U , V ) \longrightarrow  (dU+c' V = W )    
\end{array}   \]
                                                                     
The {\it purely algebraic upper part} is induced by these operators acting on the second line which is describing the invariant sections of the respective Janet bundles and we obtain the linear maps:  \\
\[   \begin{array}{lc}
{\cal{D}} \hspace{4cm}  & (K) \longrightarrow (0, c'' K)  \\
{\cal{D}}_1   \hspace{4cm}  &  (A, B) \longrightarrow  (c' (A-B)=C' , -c'' A = C'' )  \\
{\cal{D}}_2  \hspace{4cm} &  (C' , C'' ) \longrightarrow  (c'' C' + c' C'' =D )  
\end{array}   \]
only depending on the structure constants $c=(c' , c'' )$.  \\

\noindent
{\bf 7  CONCLUSION}: \\

\noindent
The work of E. Vessiot, motivated by the application of the local theory of Lie pseudogroups ([51]) to the {\it differential Galois theory}, namely the Galois theory for systems of partial differential equations ([52]), has been deliberately ignored by E. Cartan and followers ([8],[26]). As a byproduct, the 
{\it Vessiot structure equations}, introduced as early as in 1903, are still unknown today. Similarly and twenty years later but for other reasons related to his work on general relativity ([6]), in particular his correspondence with A. Einstein, Cartan did not acknowledge the work of M. Janet on the formal theory of systems of partial differential equations (Compare [7] and [23]). Accordingly, the combination of natural bundles and geometric objects with differential sequences, in particular the {\it linear and nonlinear Janet sequences} has been presented for the first time in a rather self-contained manner through this paper. The main idea is to induce from the linear Janet sequence, considered as a linear differential sequence with vector bundles and linear differential operators, a purely algebraic {\it deformation sequence} with finite dimensional vector spaces and linear maps which may not be exact. We have shown that the equivalence problem for structures on manifolds has to do with the local exactness of the Janet sequence at $F_0$ while the deformation problem of the corresponding algebraic structures has to do with the exactness of the deformation sequence at ${\Upsilon}_1=\Upsilon(F_1)$, the space of invariant sections of $F_1$, that is {\it one step further on in the sequence}. This result explains why the many tentatives done in order to link the {\it deformation of algebraic structures} like Lie algebras with the {\it deformation of geometric structures} on manifolds have not been successful. Meanwhile, we have emphasized the part that could be played by computer algebra in any effective computation and hope to have opened a new field of research for the future. Finally, ending with a wink, the reader must not forget that one of the first aplications of computer agebra in 1970 has been done in the deformation theory of Lie algebras where a few counterexamples could only be found in dimension greater than 10, that is with more than 500 structure constants !. \\

\noindent
{\bf 8  BIBLIOGRAPHY}:   \\

\noindent
[1] U. AMALDI: Congr\'{e}s de la Soci\'{e}t\'{e} Italienne pour le Progr\'{e}s des Sciences, Parma, Italy, 1907 (in italian) (Estratto dagli Annali di Matematica, XV, III, 293-327). \\
\noindent
[2] J.M. ANCOCHEA BERMUDEZ, R. CAMPOAMOR-STURSBERG, L. GARCIA VERGNOLLE, M. GOZE: Alg\`{e}bres de Lie r\'{e}solubles r\'{e}elles alg\'{e}briquement rigides, Monatsh. Math., 152, 2007, 187-195.\\
\noindent
[3] M. BARAKAT: Jets. A MAPLE-package for formal differential geometry, Computer algebra in scientific computing, (EACA Konstanz), Springer, Berlin, 2001, 1-12. \\
 http://wwwb.rwth-aachen.de/$\sim$barakat/jets/     \\
 \noindent
[4] E. CARTAN: Sur une g\'{e}n\'{e}ralisation de la notion de courbure de Riemann et les espaces \`{a} torsion, C. R. Acad\'{e}mie des Sciences Paris, 174, 1922, 437-439, 593-595, 734-737, 857-860.\\
\noindent
[5] E. CARTAN: Sur les vari\'{e}t\'{e}s \`{a} connexion affine et la th\'{e}orie de la relativit\'{e} g\'{e}n\'{e}ralis\'{e}e, Ann. Ec. Norm. Sup., 40, 1923, 325-412; 41, 1924, 1-25; 42, 1925, 17-88.\\
\noindent
[6] E. CARTAN, A. EINSTEIN: Letters on Absolute Parallelism 1929-1932 (Original letters with english translation), Princeton University Press and Acad\'{e}mie Royale de Belgique, R. Debever Ed., Princeton University Press, 1979.   \\
\noindent
[7] E. CARTAN: Sur la th\'{e}orie des syst\`{e}mes en involution et ses applications \`{a} la relativit\'{e}, Bull. Soc. Math. France, 59, 193, 88-118.  \\
\noindent
[8] E. CARTAN: La th\'{e}orie de Galois et ses diverses g\'{e}n\'{e}ralisations, Oeuvres Compl\`{e}tes, 1938. \\
\noindent
[9] C. CHEVALLEY, S. EILENBERG: Cohomology theory of Lie groups and Lie algebras, Transactions of the American Mathematical Society, 63, (1), 1948, 85-124.  \\
\noindent
[10] L.P. EISENHART: Riemannian Geometry, Princeton University Press, Princeton, 1926.\\
\noindent
[11] M. FLATO: Deformation view of physical theories, Czech. J. Phys., B, 32, 1982, 472-475.\\
\noindent
[12] J. GASQUI, H. GOLDSCHMIDT: D\'{e}formations Infinit\'{e}simales des Structures Conformes Plates, Birkhauser, 1984. \\
\noindent
[13] M. GERSTENHABER: On the deformation of rings and algebras, Ann. of Math., 79, 1964, (I) 59-103; 84, 1966, (II) 1-99; 88, 1968, (III) 1-34; 99, 1974, (IV) 257-267.  \\
\noindent
[14] H. GOLDSCHMIDT: Sur la structure des \'{e}quations de Lie, J. Differential Geometry, 6, 1972, (I) 357-373 + 7, 1972, (II) 67-95 + 11, 1976, (III) 167-223.  \\
\noindent
[15] H. GOLDSCHMIDT, D.C. SPENCER: On the nonlinear cohomology of Lie equations, Acta. Math., 136, 1973, (I) 103-170, (II) 171-239, (III) J. Diff. Geometry, 13, 1978, 409-453.        \\
\noindent
[15] R.E. GREENE, S.G. KRANTZ: Deformation of complex structures, Adv. Math. 43, 1982, 1-86.  \\
\noindent
[16] V. GUILLEMIN, S. STERNBERG: Deformation Theory of Pseudogroup Structures, Memo 66, Am. of Math. Soc., 1966.  \\
\noindent
[17] V. GUILLEMIN, S. STERNBERG: An algebraic model of transitive differential geometry, Bull. Amer. Math. Soc., 70, 1964, 16-47.  \\
\noindent
[18] V. GUILLEMIN, S. STERNBERG: The Lewy counterexample and the local equivalence problem for G-structures, J. Diff. Geometry, 1, 1967, 127-131.  \\
\noindent
[19] I. HAYASHI: Embedding and existence theorems of infinite Lie algebras, J. Math. Soc. Japan, 22, 1970, 1-14.  \\
\noindent
[20] M. HAZEWINKEL, M. GERSTENHABER: Deformation Theory of Algebras and Structures and Applications, Kluwer, NATO proceedings, 1988.   \\
\noindent
[21] G. HOCHSCHILD: On the cohomology groups of an associative algebra, Ann. of Math., 46, 1945, 58-67.  \\
\noindent
[22] E. INONU, E.P. WIGNER: On the contraction of Lie groups and Lie algebras, Proc. Nat. Acad. Sci. USA, 39, 1953, 510.  \\
\noindent
[23] M. JANET: Sur les syst\`{e}mes aux d\'{e}riv\'{e}es partielles, Journal de Math., 8, (3), 1920, 65-151. \\
\noindent
[24] K. KODAIRA, L. NIRENBERG: On the existence of deformations of complex analytic structures, Ann. of Math., 68, 1958, 450-459.  \\
\noindent
[25] K. KODAIRA, D.C. SPENCER: On the deformation of complex analytic structures, I, II, Ann. of Math., 67, 1968, 328-466.  \\
\noindent
[26] A. KUMPERA, D.C. SPENCER: Lie Equations, Ann. Math. Studies 73, Princeton University Press, Princeton, 1972.\\
\noindent
[27] M. KURANISHI: On the locally complete families of complex analytic structures, Ann. of Math., 75, 1962, 536-577.   \\
\noindent
[28] M. LEVY-NAHAS: Deformation and contraction of Lie algebras, J. Math. Phys., 8, 1967, 1211-1222.   \\
\noindent
[29] M. LEVY-NAHAS: Two simple applications of the deformations of Lie algebras, Ann. Inst. H. Poincar\'{e}, A, 13, 1970, 221-227.   \\
\noindent
[30] A. LORENZ: Jet Groupoids, Natural Bundles and the Vessiot Equivalence Method, Ph.D. Thesis (published electronically), Department of Mathematics,  RWTH Aachen-University, march 18, 2009.See also: On local integrability conditions of jet groupoids, Acta Appl. Math., 01, 2008, 205-213 at :\\
http://dx.doi.org/10.1007/s10440-008-9193-7   \\
\noindent
[31] A. NIJENHUIS, R.W. RICHARDSON: Cohomology and deformations of algebraic structures, Bull. Am. Math. Soc., 79, 1964, 406-411.  \\
\noindent
[32] V. OUGAROV: Th\'{e}orie de la Relativit\'{e} Restreinte, MIR, Moscow, 1969; french translation, 1979.\\
\noindent
[33] W. PAULI: Theory of Relativity, Pergamon Press, London, 1958.\\
\noindent
[34] W.S. PIPER: Algebraic deformation theory, J. Diff. Geometry, 1, 1967, 133-168.  \\
\noindent
[35] W.S. PIPER: Algebras of matrices under deformations, J. Diff. Geometry, 5, 1971, 437-449.  \\
\noindent
[36] J.-F. POMMARET: Systems of Partial Differential Equations and Lie Pseudogroups, Gordon and Breach, New York, 1978; Russian translation: MIR, Moscow, 1983.\\
\noindent
[37] J.-F. POMMARET: Differential Galois Theory, Gordon and Breach, New York, 1983.\\
\noindent
[38] J.-F. POMMARET: Lie Pseudogroups and Mechanics, Gordon and Breach, New York, 1988.\\
\noindent
[39] J.-F. POMMARET: Partial Differential Equations and Group Theory: New Perspectives for Applications, Kluwer, Dordrecht, 1994.\\
\noindent
[40] J.-F. POMMARET: Partial Differential Control Theory, Kluwer, Dordrecht, 2001.\\
\noindent
[41] J.-F. POMMARET: Parametrization of Cosserat equations, Acta Mechanica, 215, 2010, 43-55.\\
\noindent
[42] J.-F. POMMARET: Macaulay inverse systems revisited, Journal of Symbolic Computation, 46, 2011, 1049-1069.\\
\noindent
[43] J.-F. POMMARET: Spencer Operator and Applications: From Continuum Mechanics to Mathematical Physics, in "Continuum Mechanics-Progress in Fundamentals and Engineering Applications", Dr. Yong Gan (Ed.), ISBN: 978-953-51-0447--6, InTech, 2012, Available from: \\
http://www.intechopen.com/books/continuum-mechanics-progress-in-fundamentals-and-engineerin-applications/spencer-operator-and-applications-from-continuum-mechanics-to-mathematical-physics  \\
\noindent
[44] D.S. RIM: Deformations of transitive Lie algebras, Ann. of Math., 83, 1966, 349-357.  \\
\noindent
[45] A.N. RUDAKOV: Deformations of simple Lie algebras, Izv. Akad. Nauk. SSSR, Ser. Math., 35, 1971, 5. \\
\noindent
[46] W. M. SEILER, Involution: The Formal Theory of Differential Equations and its Applications to Computer Algebra, Springer, 2009, 660 pp. (See also doi:10.3842/SIGMA.2009.092 for a recent presentation of involution, in particular sections 3 and 4). \\
\noindent
[47] D.C. SPENCER: Some remarks on homological analysis and structures, Proc. Symp. Pure Math., 3, 1961, 56-86.   \\
\noindent
[48] D.C. SPENCER: deformation of structures on manifolds defined by transitive continuous pseudogroups, Ann. of Math., 76, 1962, (I) 306-445; 81, 1965, (II, III) 389-450.  \\
\noindent
[49] D. C. SPENCER: Overdetermined Systems of Partial Differential Equations, Bull. Am. Math. Soc., 75, 1965, 1-114.\\
\noindent
[50] M. VERGNE: Cohomologie des alg\'{e}bres de Lie nilpotentes, Bull. Soc. Math. France, 98, 1970, 81-116.
\\
\noindent
[51] E. VESSIOT: Sur la th\'{e}orie des groupes infinis, Ann. Ec. Norm. Sup., 20, 1903, 411-451.\\
\noindent
[52] E. VESSIOT: Sur la th\'{e}orie de Galois et ses diverses g\'{e}n\'{e}ralisations, Ann. Ec. Norm. Sup., 21, 1904, 9-85 (M\'{e}moire couronn\'{e} par l'Acad\'{e}mie des Sciences). \\

\end{document}